\documentclass[11pt]{article} 
\usepackage{amsmath,amssymb,amsthm}
\usepackage[usenames,dvipsnames]{color}
\usepackage{enumerate}
\usepackage{geometry}
\usepackage{graphicx}
\usepackage{svg}
\usepackage{subfigure}
\usepackage{epsfig}
\usepackage{graphicx}
\usepackage{colordvi}
\usepackage{graphics}
\usepackage{color}
\usepackage{float}
\usepackage{url}

\newcommand{\inv}{^{-1}}

\newcommand{\norm}[1]{\left|\left|{}#1\right|\right|}

\newcommand{\pare}[1]{\left({}#1\right)}
\newcommand{\brak}[1]{\left[{}#1\right]}

\newcommand{\divergence}{\nabla \cdot}

\numberwithin{equation}{section}

\newtheorem{theorem}{Theorem}[section]

\newtheorem{lemma}[theorem]{Lemma}
\newtheorem{definition}[theorem]{Definition}

\newtheorem{remark}[theorem]{Remark}

\usepackage{mathtools}

\usepackage{fullpage}

\newcommand{\R}{\mathbb{R}}

\def \ba{{\mathbf a}}

\newcommand{\bff}{\boldsymbol{f}}

\newcommand{\LR}[1]{{#1}}
\newcommand{\ti}[1]{{#1}}
\graphicspath{{./backup/EMAC-ROM_1_7_22/}}
\title{
Full and Reduced Order Model Consistency of the Nonlinearity Discretization in Incompressible Flows
}

\author{
Sean Ingimarson
\thanks{Department of Mathematical Sciences, Clemson University, Clemson, SC, 29634;
	email: singama@clemson.edu; partially supported by NSF grant DMS 2011490.}
\and	
Leo G. Rebholz
\thanks{Department of Mathematical Sciences, Clemson University, Clemson, SC, 29634;
	email: rebholz@clemson.edu; partially supported by NSF grant DMS 2011490.}
\and
Traian Iliescu
\thanks{Department of Mathematics, Virginia Tech, Blacksburg, VA, 24061;
    email: iliescu@vt.edu;
    partially supported by NSF grant DMS-2012253.}
}

\begin{document}
\date{\today}
\maketitle

\begin{abstract}
We investigate both theoretically and numerically the consistency between the nonlinear discretization in full order models (FOMs) and reduced order models (ROMs) for incompressible flows.
To this end, we consider two cases:
(i) FOM-ROM consistency, i.e., when we use the same nonlinearity discretization in the FOM and ROM; and 
(ii) FOM-ROM inconsistency, i.e., when we use different nonlinearity discretizations in the FOM and ROM.
Analytically, we prove that while the FOM-ROM consistency yields optimal error bounds, FOM-ROM inconsistency yields additional terms dependent on the FOM divergence error, which prevent the ROM from recovering the FOM as the number of modes increases.  Computationally, we consider channel flow around a cylinder and Kelvin-Helmholtz instability, and show that FOM-ROM consistency yields significantly more accurate results than FOM-ROM inconsistency.  
\end{abstract}

\section{Introduction}
    \label{sec:introduction}

{\it Reduced order models (ROMs)} are computational models whose dimensions are orders of magnitude lower than the dimensions of full order models (FOMs), i.e., models obtained with classical numerical discretizations, e.g., finite element (FE) or finite volume methods.
The main ingredient used to construct ROMs is the available numerical or experimental {\it data}.
ROMs are used in many query engineering applications, e.g., uncertainty quantification, design, and control, where the FOM computational cost is generally prohibitive.
There are many types of ROMs and several classifications.
One such classification is the following:
\begin{enumerate}
    \item {\it Nonintrusive ROMs}.
    Examples in this class include machine learning ROMs~\cite{ahmed2021closures}, sparse identification of nonlinear dynamics (SINDy)~\cite{brunton2016discovering}, or operator inference~\cite{peherstorfer2016data}.
    \item {\it Intrusive ROMs}.
    Examples in this class include {\it Galerkin ROMs (G-ROMs)} (also called projection ROMs) or Petrov-Galerkin ROMs (PG-ROMs)~\cite{carlberg2011efficient}.
    The proper orthogonal decomposition (POD)~\cite{HLB96} and reduced basis method (RBM)~\cite{hesthaven2015certified,quarteroni2015reduced} are classic examples of G-ROMs.
\end{enumerate}
Of course, each class of ROMs has its pros and cons.

For example, one of the main advantages of nonintrusive ROMs is the {\it practical convenience}.
Indeed, the nonintrusive ROMs are generally constructed without knowing too many details about the data used, e.g., the numerical discretization, the computational model, or even the underlying equations.
This makes the nonintrusive ROMs easy to implement in practice.
(We note that, in contrast, the intrusive ROMs generally require the knowledge of the underlying equations and, sometimes, even the corresponding numerical discretizations.)

One of the main advantages of the intrusive ROMs is their {\it principled} character.
For example, G-ROMs (such as POD models or RBMs) are constructed by using the {\it same Galerkin framework} as that used to build classical numerical methods, e.g., the FE method or the spectral method.
Indeed, G-ROMs are constructed by using the following algorithm:
(i) data are used to build a ROM  basis $\{ \psi_1, \ldots, \psi_r \}$;
(ii) the sought ROM solution is written as a linear combination of the ROM basis functions $u_r(x,t) = \sum_{j=1}^{r} a_j(t) \, \psi_j(x)$;
(iii) the ROM solution is used in the underlying equations, which are then projected onto the space spanned by the ROM basis functions.
After applying this three step algorithm \ti{to a PDE that is first order in time (such as the Navier-Stokes equations)}, we obtain the G-ROM:
\begin{eqnarray}
    \overset{\bullet}{\ba}
	= \bff(\ba),
	\label{eqn:g-rom-ds} 
\end{eqnarray}
where $\ba(t) = (a_{i}(t))_{i=1, \ldots,r}$ is the vector of coefficients in the Galerkin expansion in step (ii), $\bff$ comprises the ROM operators, \LR{ and $\overset{\bullet}{\ba}$ represents the time derivative of $\ba$.}
The algorithm used to construct the G-ROM~\eqref{eqn:g-rom-ds} is the same as the algorithm used to build, e.g., the FE model.
The only difference (which is, indeed, fundamental) is the basis construction:
The G-ROM is based on a {\it data-driven basis}, whereas the FE model is based on a universal basis (i.e., piecewise polynomial basis).

The Galerkin framework used to construct G-ROMs allows us to {\it prove fundamental mathematical properties} of the G-ROM, such as {\it stability, consistency, and convergence}.
These properties are proven in, e.g.,~\cite{KV01,singler2014new,koc2021optimal} for POD, and~\cite{hesthaven2015certified,quarteroni2015reduced} for RBM.
We emphasize that, in contrast, these fundamental mathematical properties are generally not proven for \ti{nonintrusive ROMs (e.g., data-driven models based on machine learning)}.

Recently, another step toward laying the {\it mathematical foundations} for G-ROMs has been taken.
Specifically, several studies have addressed the following natural questions:
{\it If we use a Galerkin framework to construct both the FOM and the ROM, should we also use the same numerical discretization for the two models?
Would this approach yield more accurate ROMs?}
ROMs that satisfy this property are called FOM consistent. 
\ti{Next, we propose a possible definition of FOM-ROM consistency.}



\begin{definition}[FOM-ROM Consistency]
A ROM is FOM consistent if it uses the same computational model and the same numerical discretization \ti{(e.g., spatial discretization, time discretization, and nonlinearity and pressure treatment)} as the FOM.
\end{definition}

\ti{
\begin{remark}[Projection]
All the ROMs investigated in this paper are constructed by using the continuous projection~\cite[Chapter 2]{volkwein2013proper}.
Alternatively, one can construct the ROMs by using the discrete projection~\cite[Chapter 1]{volkwein2013proper}.
In that case, the system of ODEs resulting from the FOM spatial discretization is projected onto a ROM basis in $\R^r$ by using a discrete inner product~\cite[Section 4.3.2]{kalashnikova2014reduced}. 
Thus, the ROM constructed by using the discrete projection is consistent with the FOM with respect to the spatial discretization, provided no extra terms are added to the FOM equations before the projection or the ROM equations after the projection.
We note, however, that the ROM can be inconsistent with the FOM with respect to the time discretization.
\end{remark}
}

\begin{remark}[Parameters]
In the definition of the FOM-ROM consistency, we require that the computational model and numerical discretization be the same in the FOM and ROM.
We note, however, that the parameters in the FOM and ROM models and discretizations can be different.
\end{remark}

{\it Computationally}, the FOM-ROM consistency has been investigated with respect to 
the computational model, but not with respect to the numerical discretization:

\begin{enumerate}
    \item {\it Closure Model}: \ 
    The consistency of the closure model~\cite{ahmed2021closures} used in the FOM and ROM was investigated \ti{numerically}  in~\cite{rebollo2017certified}.
    \item {\it Numerical Stabilization}: \ 
    \ti{The consistency with respect to the stabilization was investigated numerically for} 
    residual-based stabilization methods, e.g., SUPG~\cite{ali2020stabilized,giere2015supg, pacciarini2014stabilized}, the variational multiscale method~\cite{stabile2019reduced}, the Leray model~\cite{girfoglio2021pod}, and the evolve-filter-relax model~\cite{strazzullo2021consistency}.
\end{enumerate}

We emphasize that, in all these numerical investigations, FOM-ROM consistency was found to be beneficial.
Specifically, ROMs that were consistent with the corresponding FOMs were more accurate than ROMs that were inconsistent.

\ti{
The \textit{theoretical} investigations of the FOM-ROM consistency are scarce.
Indeed, to our knowledge, the \textit{numerical analysis} of the FOM-ROM consistency has been performed only in two studies~\cite{giere2015supg,pacciarini2014stabilized}. 
Furthermore, in both studies the FOM-ROM consistency was investigated with respect to the computational model (i.e., the SUPG stabilization), but not with respect to the numerical discretization.
Finally, both studies considered \textit{linear} problems (i.e., the convection-diffusion-reaction equations).
Specifically, error bounds for the FOM-ROM inconsistent approach (in which the SUPG stabilization is used at the FOM level but not at the ROM level) were proven in Proposition 3.1 in~
\cite{pacciarini2014stabilized}.
Furthermore, in Section 3.3 of~
\cite{pacciarini2014stabilized} it was argued that using different bilinear forms at the FOM and ROM levels in the inconsistent approach increases the error.
A different theoretical investigation was performed in~\cite{giere2015supg}, where \textit{a priori} error bounds for the consistent approach (in which SUPG is used both at a FOM and at a ROM level) were leveraged to determine two different types of scalings for the SUPG stabilization parameter: a scaling based on the underlying FE discretization, and another scaling based on the ROM discretization.
}

Investigating the FOM-ROM consistency \ti{both computationally and theoretically} is an important research direction in reduced order modeling~\cite{barone2009stable,kalashnikova2014reduced}.
Indeed, despite the fact that several 
\ti{numerical} investigations have shown that the ROMs that are FOM consistent are more accurate than the ROMs that are FOM inconsistent, most of the current ROMs are inconsistent.
As explained in Remark 1 in~\cite{strazzullo2021consistency}, ROMs can be  inconsistent with respect to the computational model~\cite{bergmann2018zonal,carlberg2017galerkin,grimberg2020stability,mou2021data,wang2012proper,wells2017evolve}, stabilization~\cite{ali2020stabilized,pacciarini2014stabilized,stabile2019reduced}, and time discretization~\cite{carlberg2017galerkin,wang2012proper,zucatti2021calibration}.
These inconsistent ROMs are commonly used in practice because they are convenient: to construct the ROM, one is not restricted by the particular modeling and discretization choices made at a FOM level.
Thus, providing both theoretical and computational support for using consistent ROMs could determine more practioners to use them in applications.

In this paper, we take 
\ti{a significant} step in the investigation of the FOM-ROM consistency and investigate both theoretically and numerically the FOM-ROM consistency with respect to an important component of the numerical discretization of incompressible flows, i.e., the {\it nonlinearity}.
To our knowledge, this is the \ti{\textit{first investigation of the FOM-ROM consistency for the nonlinearity discretization and for numerical discretization in general}}.
We show both analytically and numerically that ROM that are inconsistent with respect to the nonlinearity discretization are significantly less accurate than consistent ROMs.
Specifically, we prove that the error bounds for inconsistent ROMs include terms that are not part of the error bounds for consistent ROMs.
Furthermore, these additional error terms can prevent the inconsistent ROM from recovering the corresponding FOM as the number of ROM basis functions increases.

The rest of the paper is organized as follows.
In Section~\ref{sec:fom}, we present the FOM, which is based on a FE spatial discretization.
In Section~\ref{sec:rom}, we outline the ROM and list several results for its numerical analysis.
In Section~\ref{sec:analysis}, we prove an error bound for an inconsistent ROM that is larger than the error bound for the corresponding consistent ROM.
This shows, theoretically, that consistent ROMs are more accurate than inconsistent ROMs.
The same conclusion is reached in the numerical investigation in Section~\ref{sec:numerical-results} for two test problems, channel flow around a cylinder and Kelvin-Helmholtz instability.  
Finally, in Section~\ref{sec:conclusions}, we summarize our findings and outline future research directions.


\section{Full Order Model (FOM)}
    \label{sec:fom}

In our theoretical and numerical investigation, 
we use the accepted physical model for the evolution of the incompressible Newtonian fluid flow, which is the Navier-Stokes equations (NSE), given  in a domain $\Omega$ with unknowns $u$ and $p$ representing velocity and pressure, $f$ the external forcing, and $\nu$ the
inverse of the Reynolds number, $Re$,
by
\begin{eqnarray}
&& u_t + u\cdot\nabla u + \nabla p - \nu\Delta u = f,
\label{eqn:nse-1} \\[0.1cm]
&& \nabla \cdot u = 0,
\label{eqn:nse-2}
\end{eqnarray}
along with an appropriate initial condition $u^0$ and no-slip (for simplicity) velocity boundary conditions. 

To construct the FOM, we use a FE spatial discretization.
We note that, while
the conservation of mass equation of the NSE holds pointwise (see~\eqref{eqn:nse-2}), in numerical simulations generally it does not.  In 
FE simulations, the conservation holds generally only in a global sense.  To better discuss this, consider a 
FE scheme for the NSE using backward Euler time stepping (for simplicity of analysis to follow; any stable time stepping method could be used) and inf-sup stable velocity-pressure element spaces $X_h\subset H^1_0(\Omega) \cap P_k(\tau_h)$,  $Q_h\subset L^2_0(\Omega) \cap P_{k-1}(\tau_h)$, for a regular triangulation $\tau_h$ of $\Omega$, and $k\ge 2$ (i.e.,  Taylor-Hood elements for simplicity of analysis to follow; any other inf-sup stable pair could also be used).  Denoting the $L^2(\Omega)$ inner product by $(\cdot,\cdot)$, the scheme reads: given $u_h^0\in X_h$, for $n=0,1,2,...$, find $(u_h^{n+1},p_h^{n+1}) \in (X_h,Q_h)$ satisfying
\begin{align}
\frac{1}{\Delta t}\left( u_h^{n+1} - u_h^n,v_h\right) + b\left(u_h^{n+1},u_h^{n+1},v_h\right) - \left(p_h^{n+1},\nabla \cdot v_h\right) + \nu\left(\nabla u_h^{n+1},\nabla v_h\right) & = \left(f(t^{n+1}),v_h\right), \label{ns1} \\
\left(\nabla \cdot u_h^{n+1},q_h\right) &= 0, \label{ns2}
\end{align}
for all $v_h\in X_h,\ q_h \in Q_h$.  The $b$ operator can represent any of several nonlinear forms, including:
\begin{align}
b_{c}(u,u,v) & = (u\cdot\nabla u,v) & \mbox{Convective form} \label{eqn:bc} \\
b_{s}(u,u,v) & = (u\cdot\nabla u,v) + \frac12 ((\nabla \cdot u)u,v) & \mbox{Skew-symmetric form} \label{eqn:bs} \\\
b_{r}(u,u,v) & = ((\nabla \times u) \times u,v) & \mbox{Rotational form} \label{eqn:br} \\
b_e(u,u,v) & = 2(D(u)u,v) + ((\nabla \cdot u)u,v) 
& \mbox{EMAC form,}
\label{eqn:be}
\end{align}
\LR{where EMAC stands for Energy, Momentum, and Angular Momentum Conserving.  See \cite{CHOR17} and \cite{GS98} for more on these and other NSE formulations.}

We note that with Rotational and EMAC forms the pressure variable will represent a dynamic pressure that adds $\frac12 |u|^2$ to usual pressure for Rotational form and subtracts it for EMAC.  If $\nabla\cdot u=0$, then all these nonlinear forms are equivalent.  However, due to the weak FE enforcement in \eqref{ns2}, these four forms can give very different behavior in simulations \cite{CHOR17,LMNOR09,OR20} and even have different physical properties such as conservation (or not) of energy, momentum, and angular momentum - for example the EMAC nonlinear form preserves all of these quantities while the convective nonlinear form preserves none of them \cite{CHOR17}.  Numerical tests in \cite{CHOR19,CHOR17} show how these conservation properties correspond to better stability and accuracy properties, and it was also shown in \cite{OR20} that EMAC reduces the Gronwall constant compared to other schemes (reducing or eliminating its dependence on $Re$) which provides for better long time accuracy.




There is extensive literature on the accuracy of \eqref{ns1}-\eqref{ns2}, see e.g. \cite{Laytonbook,OR20,temam}, which proves the following for $(X_h,Q_h)=(P_k,P_{k-1})$ Taylor-Hood elements if the Skew-symmetric, Rotational or EMAC forms are used.  For a sufficiently smooth NSE solution and sufficiently small time step size,  
\begin{equation}
    \| u(T) - u_h^M \|^2 + \nu \Delta t \sum_{n=1}^M \| \nabla (u(t^n) - u_h^n) \|^2 \le C (\Delta t^2 + h^{2k}). \label{feerr1}
\end{equation}
If BDF2 time stepping is used, then the $\Delta t^2$ becomes a $\Delta t^4$, and for EMAC the right hand side constant $C$ does not depend exponentially on $Re$ while for Rotational or Skew-symmetric forms it does \cite{OR20}.  For Convective form, since this scheme lacks unconditional stability, similar results may be possible but would require a much more technical analysis and 
a time step restriction.

\section{Reduced Order Model (ROM)}
    \label{sec:rom}

To construct the Galerkin ROM, we need to build the ROM basis.
To this end, we use the proper orthogonal decomposition (POD)~\cite{HLB96}.
First, we take FOM snapshots, i.e., we construct  
a matrix with columns that are solutions to the FE model  \eqref{ns1}-\eqref{ns2}:
\[
A = \{ u_h^0,\ u_h^1,\ ..., u_h^M \}.
\]
Then, we perform an eigenvalue decomposition on $A$ 
(using an appropriate inner product, in this case $L^2$) and obtain a set of orthonormal functions.
Next, we choose an $r$ value that is much smaller than the FE space dimension, 
and use 
the top $r$ 
modes to create the ROM space~\cite{HLB96}
\[
V_r = \mbox{span}\left\{ \psi_1,\ \psi_2,...,\ \psi_r \right\}.
\]
The next step in the Galerkin ROM construction is to replace $u$ with $w_r \in V_r$ in the NSE~\eqref{eqn:nse-1}--\eqref{eqn:nse-2} and project the resulting equations onto $V_r$.
We note that, by construction, functions in $V_r$ automatically satisfy the weak FE conservation of mass constraint \eqref{ns2}.  
Thus, the pressure term in the resulting equations vanishes.
The Galerkin ROM (G-ROM) now reads: Given $w_r^0\in V_r$, for $n=0,1,2,...$, find $w_h^{n+1}\in V_r$ satisfying
\begin{equation}
\frac{1}{\Delta t}(w_r^{n+1}-w_r^n,v_r) + b(w_r^{n+1},w_r^{n+1},v_r) + \nu(\nabla w_r^{n+1},v_r)  = (f(t^{n+1}),v_r) \quad \forall v_r \in V_r. \label{rom1}
\end{equation}
As in the FE model~\eqref{ns1}-\eqref{ns2}, 
the nonlinear form $b$ can take one of 
the forms \eqref{eqn:bc}--\eqref{eqn:be}.


One can prove that the POD projection error satisfies the following equality~\cite{singler2014new}:
\begin{eqnarray}
\frac{1}{M+1} \, \sum_{j=0}^{M} \left\| \nabla \left( u_h(\cdot, t_j) - \sum_{i=1}^{r} \bigl( u_h(\cdot, t_j) , \psi_i(\cdot) \bigr) \, \psi_i(\cdot) \right) \right\|^2
= \sum_{k=r+1}^{d} \| \nabla \psi_k \|^{2} \, \lambda_k,
    \label{eqn:projection-error}
\end{eqnarray}
where $\| \cdot \|$ is the $L^2(\Omega)$ norm, $d$ is the rank of the snapshot matrix, $A$, and $\lambda_j$ is the eigenvalue corresponding to the ROM basis function $\psi_j$.

 
 
\section{Analysis of FOM-ROM inconsistency}
    \label{sec:analysis}

In this section, we perform an error analysis that reveals how FOM-ROM inconsistency of the nonlinearity discretization
creates additional error.  To this end, we consider a ROM that is inconsistent with the FOM with respect to the nonlinearity discretization.
Specifically, we assume that the FEM simulation was done with the Skew-symmetric form of the nonlinear term ($b=b_s$ in \eqref{ns1}) and EMAC is used in the ROM ($b=b_e$ in \eqref{rom1}).  We note that the analysis that follows would hold in essentially the same way with any other combination of nonlinear forms that are not the same in the FEM and ROM schemes (although with additional difficulties if $b_c$ is used).

 In the analysis to follow, 
 norms other than $L^2$ 
 are labeled with subscripts.
We will also use the following discrete Gronwall lemma.

\begin{lemma}[\LR{Lemma 27 in~\cite{layton2008introduction}}]
    Let $\Delta t, H,$ and $a_n, b_n, c_n, d_n$ (for integers $n \geq 0$) be nonnegative numbers such that
    \begin{eqnarray}
        a_l
        + \Delta t  \sum_{n=0}^{l} b_n
        \leq \Delta t  \sum_{n=0}^{l} d_n a_n
        + \Delta t  \sum_{n=0}^{l} c_n
        + H
        \qquad \text{for } l \geq 0.
        \label{eqn:lemma-gronwall-1}
    \end{eqnarray}
    Suppose that $\Delta t \, d_n < 1$.
    Then,
    \begin{eqnarray}
        a_l
        + \Delta t  \sum_{n=0}^{l} b_n
        \leq \exp \left( \Delta t  \sum_{n=0}^{l} \frac{d_n}{1 - \Delta t \, d_n} \right)
        \left( \Delta t  \sum_{n=0}^{l} c_n
        + H \right)
        \qquad \text{for } l \geq 0.
        \label{eqn:lemma-gronwall-2}
    \end{eqnarray}
    \label{lemma:gronwall}
\end{lemma}

\begin{theorem}\label{FOMROMerr}
Let $\{ u_h^n\}_{n=1}^{M}$ be the FE  solution from \eqref{ns1} with $b=b_s$, and $\{ w_r^n \}_{n=1}^{M}$ be the ROM solution from \eqref{rom1} with $b=b_e$ and with initial condition $w_r^0 = P_{L_2}^{V_r}(u_h^0)$, where $P_{L^2}^{V_r}$ is the $L^2$ projection onto the ROM space.
Assuming that the FE  solution satisfies $C_u=\max_n \{ \| \nabla u_h^n \| \}_{n=1}^{M}$ is bounded independent of any discretization parameters and $\Delta t$ is sufficiently small, we have the following error bound for the ROM as an approximation of the FOM: 
\begin{multline}
\norm{w_r^M - u_h^M}^2
 +\nu \, \Delta t \sum_{n=1}^{M}\norm{\nabla (w_r^n - u_h^n)} ^2 
 \\ \leq C e^{\pare{\LR{ (\nu^{-3}C_u^4+\nu^{-1})} T}} \bigg( 
\nu^{-1}(C_u^2+1)  \sum_{k=r+1}^d \| \nabla \psi_k \|^2 \lambda_k  
+ \underbrace{C_u^2 \|\nabla \cdot 
\left(P_{L^2}^{V_r}(u_h^n) - w_r^n \right) \|_{1,0} 
+ C_u^2 \|\nabla \cdot u_h^{n} \|^2_{2,0} }_{\text{arises from FOM-ROM inconsistency}}
\bigg) ,
\label{eqn:errthm}
\end{multline}
where $\lambda_k$ are the eigenvalues in the eigenvalue problem used to construct the ROM basis~\cite{KV01,singler2014new}, 
$T = 
M \Delta t$ is the final time, 
and $\| v \|_{p,k} := \left( \Delta t \, \sum_{n=1}^{M} \| v^{n} \|_{k}^{p} \right)^{1/p}$
is the discrete time norm~\cite{layton2008introduction}.
\end{theorem}

\begin{proof}[Proof of Theorem \ref{FOMROMerr}]

  Denote $e^n = u_h^n - w_r^n = u_h^n - P_{L^2}^{V_r}(u_h^n) + P_{L^2}^{V_r}(u_h^n) - w_r^n =: \eta^n + \phi_r^n$. 
  Subtracting \eqref{rom1} from \eqref{ns1} with $v_h=v_r\in V_r$ yields
\begin{equation}
\frac{1}{\Delta t}(e_r^{n+1}-e_r^n,v_r) + b_s\left(u_h^{n+1},u_h^{n+1},v_r\right)  -  b_e(w_r^{n+1},w_r^{n+1},v_r) + \nu(\nabla e^{n+1},\nabla v_r)  = 0 \quad \forall v_r \in V_r, \label{err1}
\end{equation}
with the pressure term vanishing thanks to $v_r\in V_r$ being discretely divergence free.  Choosing  $v_r=\phi_r^{n+1}$, we obtain 

\begin{align}
\begin{split}
& \frac{1}{2\Delta t}\pare{\norm{\phi_r^{n+1}}^2-\norm{\phi_r^{n}}^2}+\nu\norm{\nabla \phi_r^{n+1}}^2 \le \\
  & - \left(b_s\left(u_h^{n+1},u_h^{n+1},\phi_r^{n+1}\right)  -  b_e(w_r^{n+1},w_r^{n+1},\phi_r^{n+1}) \right)-\nu(\nabla \eta^{n+1},\nabla \phi_r^{n+1}),\label{err3}
\end{split}
\end{align}
using that $(\eta^{n+1}-\eta^n,\phi_r^{n+1})=0$ since $\eta^n, \: \eta^{n+1}$ are $L^2$-orthogonal to $\phi_r^{n+1}$ by construction.  The last term on the right-hand side of~\eqref{err3} is bounded in the usual way with Cauchy-Schwarz and Young's inequalities via
\begin{align*}
-\nu(\nabla \eta^{n+1},\nabla \phi_r^{n+1}) & \le \nu \| \nabla \eta^{n+1} \| \| \nabla \phi_r^{n+1} \| \\
& \le 2\nu \| \nabla \eta^{n+1} \|^2 + \frac{\nu}{8} \| \nabla \phi_r^{n+1} \|^2.
\end{align*}
We begin the analysis of the nonlinear terms by adding and subtracting $b_e(u_h^{n+1},u_h^{n+1},\phi_r^{n+1})$ and then applying the triangle inequality to get
\begin{align}
|b_s &\left(u_h^{n+1}, u_h^{n+1},\phi_r^{n+1}\right)  -  b_e(w_r^{n+1},w_r^{n+1},\phi_r^{n+1})| \nonumber \\ 
& \ \leq 
|b_s(u_h^{n+1},u_h^{n+1},\phi_r^{n+1})-b_e(u_h^{n+1},u_h^{n+1},\phi_r^{n+1})|
 +|b_e(u_h^{n+1},u_h^{n+1},\phi_r^{n+1})-b_e(w_r^{n+1},w_r^{n+1},\phi_r^{n+1})|.
 \label{trilin2}
\end{align}
For the second term on the right-hand side of 
\eqref{trilin2}, we follow analysis from \cite{OR20} to obtain
\begin{multline}
|b_e(u_h^{n+1}, u_h^{n+1},\phi_r^{n+1})-b_e(w_r^{n+1},w_r^{n+1},\phi_r^{n+1})| \\
\leq\  C\nu\inv\brak{\norm{\nabla u_h^{n+1}}^2\norm{\eta^{n+1}}\norm{\nabla \eta^{n+1}}+\norm{u_h^{n+1}}\norm{\nabla u_h^{n+1}}\norm{\nabla \eta^{n+1}}^2}
\\
+M_0\nu^{-3}\norm{\nabla u_h^{n+1}}^4\norm{\phi_r^{n+1}}^2+\frac{\nu}{8}\norm{\nabla \phi_r^{n+1}}^2,\label{trilin6}
\end{multline}
where $M_0$ is a constant depending only on the size of $\Omega$.

%
%

To bound the first 
term on the right-hand side of \eqref{trilin2}, we expand the definitions of $b_s$ and $b_e$ and use that
\begin{multline*}
2(D(u_h^{n+1})u_h^{n+1},\phi_r^{n+1}) =
((\nabla u_h^{n+1}) u_h^{n+1} ,\phi_r^{n+1})
+((\nabla u_h^{n+1})^T u_h^{n+1} ,\phi_r^{n+1})
= \\
(u_h^{n+1}\cdot\nabla u_h^{n+1} ,\phi_r^{n+1})
+(\phi_r^{n+1}\cdot\nabla u_h^{n+1} ,u_h^{n+1})
\end{multline*}
to get
\begin{align}
|b_s &\left(u_h^{n+1}, u_h^{n+1},\phi_r^{n+1}\right)  -  b_e(u_h^{n+1},u_h^{n+1},\phi_r^{n+1})|\notag\\
=&|(u_h^{n+1}\cdot\nabla u_h^{n+1},\phi_r^{n+1})+\frac{1}{2}((\divergence u_h^{n+1})u_h^{n+1},\phi_r^{n+1})-(u_h^{n+1}\cdot \nabla u_h^{n+1},\phi_r^{n+1})\notag \\
& \ \ \ \ \ -(\phi_{r}^{n+1}\cdot\nabla u_h^{n+1},u_h^{n+1})
 -((\divergence u_h^{n+1})u_h^{n+1},\phi_r^{n+1})|\notag\\
=&|(\phi_r^{n+1}\cdot\nabla u_h^{n+1},u_h^{n+1})+\frac{1}{2}((\divergence u_h^{n+1})u_h^{n+1},\phi_r^{n+1})|.\label{trilin3}
\end{align}
Note that these terms arise only because of the inconsistency between the FOM and ROM, and we now bound both of them.  For the first 
term, after using 
the divergence theorem \cite{Laytonbook} and then H\"older and Sobolev inequalities, we obtain
\begin{align}
(\phi_r^{n+1} \cdot \nabla u_h^{n+1},u_h^{n+1})
& =-\frac{1}{2}((\divergence \phi_r^{n+1})u_h^{n+1},u_h^{n+1}) \notag \\
& \leq \frac12 \| \nabla \cdot \phi_r^{n+1}\| \| u_h^{n+1} \|_{L^6} \| u_h^{n+1} \|_{L^3} \notag \\
& \le \frac12 M_0\norm{\divergence \phi_r^{n+1}}\norm{\nabla u_h^{n+1}}^{\frac{3}{2}}\norm{u_h^{n+1}}^{\frac{1}{2}}.\label{trilin4}
\end{align}

\LR{
For the second term from \eqref{trilin3}, we 
\ti{denote $C$ to be a generic constant depending only on the domain size, and} 
apply H\"older's inequality, Sobolev inequalities, and Young's inequality to get
\begin{align}
\frac{1}{2}((\divergence u_h^{n+1})u_h^{n+1},\phi_r^{n+1})
&
\leq \frac12 \norm{\divergence u_h^{n+1}}\norm{u_h^{n+1}}_{L^6}\norm{\phi_r^{n+1}}_{L^3} \nonumber \\
& \leq \ti{C^{\frac{1}{2}}} \norm{\divergence u_h^{n+1}}\norm{\nabla u_h^{n+1}} \norm{\phi_r^{n+1}}^{1/2} \norm{\nabla \phi_r^{n+1}}^{1/2} \nonumber \\
& \leq C C_u^2 \norm{\divergence u_h^{n+1}}^2 + \norm{\phi_r^{n+1}} \norm{\nabla \phi_r^{n+1}} \nonumber \\
& \leq C C_u^2 \norm{\divergence u_h^{n+1}}^2 + 2\nu^{-1} \norm{\phi_r^{n+1}}^{\ti{2}} + \frac{\nu}{8} \norm{\nabla \phi_r^{n+1}}^{\ti{2}} . \label{trilin5}
\end{align}

}
%
Combining the above bounds 
and applying the Poincar\'e inequality, 
we get
\begin{multline*}
\frac{1}{2\Delta t}\pare{\norm{\phi_r^{n+1}}^2-\norm{\phi_r^{n}}^2}
+\frac{\nu}{2}\norm{\nabla \phi_r^{n+1}}^2 \le M_0\nu^{-3}\norm{\nabla u_h^{n+1}}^4\norm{\phi_r^{n+1}}^2 \LR{  + 2\nu^{-1} \norm{\phi_r^{n+1}}} \\
+  C\nu\inv\pare{ \norm{\nabla u_h^{n+1}}^2\norm{\nabla \eta^{n+1}}^2} + 2\nu \norm{\nabla \eta^{n+1}}^2\\
 \underbrace{ +\frac{1}{2}M_0\norm{\divergence \phi_r^{n+1}}\norm{\nabla u_h^{n+1}}^{\frac{3}{2}}\norm{u_h^{n+1}}^{\frac{1}{2}} + \LR{C C_u^2} \norm{\divergence u_h^{n+1}}^2}_{\text{FOM-ROM inconsistency}}.
\end{multline*}
Using the smoothness 
assumption on $u_h$ and summing over time steps provides us with
\begin{multline}
\| \phi_r^M \|^2 + \nu \Delta t \sum_{n=1}^{M} \norm{\nabla \phi_r^{n}}^2 
\le \| \phi_r^0 \|^2  + \Delta t \sum_{n=1}^{M}  C \LR{ (C_u^4 \nu^{-3}+\nu^{-1})} \norm{\phi_r^{n}}^2 
+ 4\nu \Delta t \sum_{n=1}^{M} \norm{\nabla \eta^{n}}^2\\
+ C \nu\inv C_u^2 \Delta t \sum_{n=1}^{M}  \| \nabla \eta^{n} \|^2  
 + \underbrace{ C C_u^2 \Delta t \sum_{n=1}^{M} \left(  \norm{\divergence \phi_r^{n}} +\norm{\divergence u_h^{n}}^2 \right) }_{\text{FOM-ROM inconsistency}}.
\end{multline}
Note that $\phi_r^0=0$ by definition of the initial condition for the ROM being $w_r^0 = P_{L_2}^{V_r}(u_h^0)$.  Hence, we can now apply the discrete Gronwall inequality in Lemma~\ref{lemma:gronwall}, the POD projection error formula~\eqref{eqn:projection-error}~\cite{singler2014new}, and the triangle inequality to find that, for sufficiently small $\Delta t$, 
we get a bound for the ROM as an approximation of the FOM: 
\begin{multline}
\norm{w_r^M - u_h^M}^2
 + \nu \Delta t \sum_{n=1}^{M}\norm{\nabla (w_r^n - u_h^n)} ^2 
\\ \leq C e^{\pare{ \LR{(\nu^{-3}C_u^4+\nu^{-1})} T}} \bigg( 
\nu^{-1}(C_u^2+1)  \sum_{k=r+1}^d \| \nabla \varphi_k \|^2 \lambda_k  
+ \underbrace{C_u^2 \|\nabla \cdot \phi_r \|_{1,0} 
+  C_u^2 \|\nabla \cdot u_h \|^2_{2,0} }_{\text{FOM-ROM inconsistency}}
\bigg) ,
\label{eqn:gronwall}
\end{multline}
where $\lambda_k$ are the eigenvalues in the eigenvalue problem used to construct the ROM basis~\cite{KV01}.  This completes the proof.

\end{proof}

\begin{remark}
For the smoothness assumption on $u_h$ in the theorem, if for example $(P_2,P_1)$ Taylor-Hood elements are used for the FE spatial discretization and a BDF2 
scheme with 
$\Delta t \le O(h^{1/2})$
is used for the time discretization, then using the inverse inequality, standard interpolation estimates, convergence of the FE \eqref{feerr1}, and the $H^1$  stability of the $L^2$ projection 
into the FE space \cite{BPS02}, $P_{L^2}^{X_h}$, we obtain
\begin{align*}
\| \nabla u_h^n \| & \le \| \nabla (u_h^n - u^n) \| + \| \nabla u^n \|  \\
&\le \| \nabla (u_h^n - P_{L^2}^{X_h}(u^n) ) \| + \| \nabla (u^n - P_{L^2}^{X_h}(u^n) ) \| +  \| \nabla u^n \| \\
& \le Ch^{-1} \| u_h^n - P_{L^2}^{X_h}(u^n)  \| + (1+C_{\Omega}) \| \nabla u^n \| +  \| \nabla u^n \| \\
& \le Ch^{-1} \left(   \| u_h^n - u^n  \| +  \| u^n - P_{L^2}^{X_h}(u^n)  \| \right) + (2+C_{\Omega}) \| \nabla u^n \| \\
& \le Ch^{-1} \left(   h^2 + 
\Delta t^2 \right)    + (2+C_{\Omega}) \| \nabla u^n \| \\
& \le C,
\end{align*}
where $C$ depends on the true solution and problem data but is independent of discretization parameters.
\end{remark}

\begin{remark}
The theorem reveals how the FOM-ROM inconsistency plays a role in the error: if there is FOM-ROM consistency, the last two terms on the right hand side of~\eqref{eqn:errthm} vanish, and thus as $r$ is increased the ROM solution recovers the FOM solution.  However, if there is FOM-ROM inconsistency in the nonlinear formulation, then no matter how many modes are used, the ROM error will not be better than the divergence error from the FOM (which is generally not very small \cite{JLMNR17}).
{\it Thus, ROMs that are FOM inconsistent do {\it not} converge with respect to the ROM dimension, $r$.}
\end{remark}

\begin{remark}
The theorem shows that if the FOM gives divergence-free solutions (e.g., using 
Scott-Vogelius FEs or others, see e.g. \cite{JLMNR17}), then the FOM-ROM inconsistency error is eliminated.  This is perhaps not surprising, since with divergence-free elements all of the nonlinearity formulations 
\eqref{eqn:bc}--\eqref{eqn:be}
are equivalent, up to a potential term.  
We note, however, that implementing divergence-free FEs can be challenging, and so the FOM-ROM inconsistency is relevant for the majority of FEs used in practical computations.
We also note that, 
if grad-div stabilization is used in the FEM, then the FOM-ROM inconsistency error can be reduced since the divergence error of the FEM solution is reduced.
\end{remark}

\section{Numerical Tests}
    \label{sec:numerical-results}

We give results using two numerical experiments, channel flow around a cylinder and Kelvin-Helmholtz instability.  Both experiments agree with the above analysis: better results are obtained when the ROM formulation matches the FEM formulation whose snapshots created the ROM space, i.e., the FOM-ROM consistency of the nonlinearity discretization is beneficial.  

\LR{The velocity-pressure spaces for the FOM (FEM simulations) are chosen to be $(P_2,P_1)$ Taylor-Hood elements in all of our tests, along with BDF2 time stepping.  The nonlinearity is resolved with Newton's method at each time step, in the FOMs and the ROMs.  The FOM simulations create the snapshots, and we follow the procedure in \cite{CIJS15} to construct the ROM modes and ROM space.}
\ti{
We assemble the ROM operator  corresponding to the nonlinearity in the NSE without using hyperreduction~\cite{barrault2004eim,yano2019discontinuous}.
We note, however, that hyperreduction may be needed for other types of PDEs.
}

In this section, we 
denote FOMs by EMAC-FEM, SKEW-FEM and CONV-FEM, and ROMs by EMAC-ROM, SKEW-ROM, and CONV-ROM.

\subsection{Channel flow around a cylinder}
    \label{sec:channel}

\begin{figure}[h!]
\begin{center}
\includegraphics[width=0.65\textwidth,height=0.25\textwidth, trim=0 0 0 0, clip]{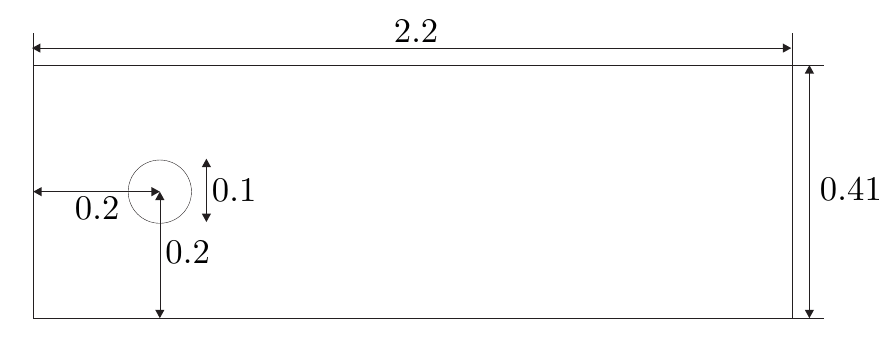}
\end{center}
\caption{The domain for the channel flow past a cylinder numerical experiment.\label{cyldomain} }
\end{figure}

\begin{figure}[h!]
\begin{center}
\includegraphics[width=0.65\textwidth,height=0.17\textwidth, trim=230 25 130 0, clip]{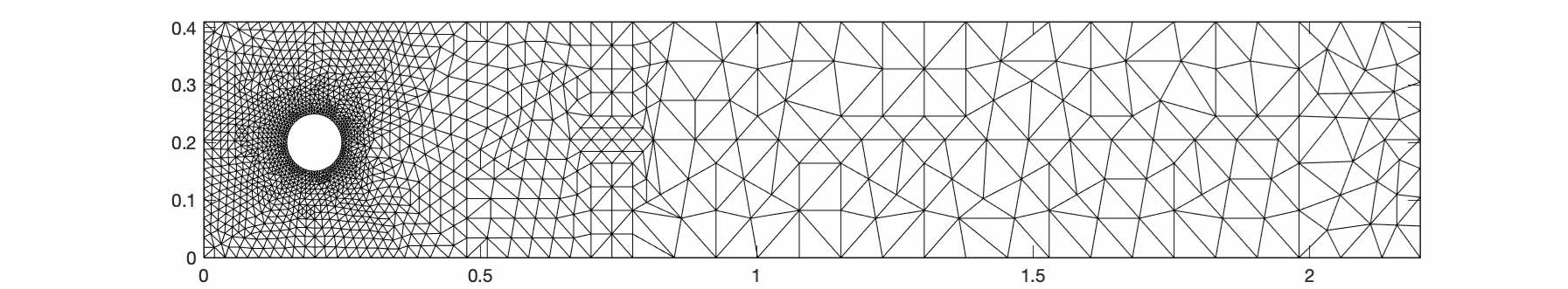}
\end{center}
\caption{The mesh used for the FEM computations for channel flow past a cylinder numerical experiment.\label{cylmesh} }
\end{figure}

For our first set of tests, we consider ROM simulations for 2D channel flow past a cylinder at $Re$=200.   We follow the setup from \cite{J04,ST96}, and use a $2.2\times0.41$ rectangular channel domain containing a cylinder (circle) of radius $0.05$ centered at $(0.2,0.2)$; the domain is shown in Figure \ref{cyldomain} \LR{ and the Delaunay generated 2886 element triangular mesh is shown in Figure \ref{cylmesh}.}   There is no external forcing ($f=0$), the kinematic viscosity is taken to be $\nu=0.0005$, no-slip boundary conditions are prescribed for the walls and the cylinder, and the inflow and outflow velocity profiles are given by
\begin{align*}
u_1(0,y,t)  = u_1(2.2,y,t) = \frac{6}{0.41^2}y(0.41-y), \ \ \
u_2(0,y,t)  = u_2(2.2,y,t) = 0.
\end{align*}

Snapshots to create the ROMs are generated by computing with EMAC-FEM and CONV-FEM as the FOMs, both using $(P_2,P_1)$ Taylor-Hood elements and 
BDF2 time stepping with $\Delta t=0.001$.  
We use a relatively coarse mesh having 12K velocity degrees of freedom (dof), noting 
\ti{that, for many practical applications (e.g., engineering and geophysical turbulent flows), ROMs}
are created from FOMs that are not fully resolved~\cite{ahmed2021closures}.  Each FOM is started from rest, and run to $t=15$.  By $t=4$, a periodic-in-time/statistically steady state is reached.  ROM spaces are created from a full period of snapshots.
Each ROM is run starting from the FOM solution at $t=5$ (projected into the ROM space), and run for 10 time units (which is $t=15$ for the FOM) using the same time step size as the FOM, $\Delta t=0.001$.  Comparisons of the ROMs and FOMs are done using energy and drag prediction, \LR{ with energy at time $t$ defined as $\frac12 \| u(t) \|^2$ and the} drag coefficient defined by \cite{J04}
\begin{align*}
c_d(t) &= 20\int_S \left( \nu \frac{\partial u_{t_S}(t)}{\partial n}n_y - p(t)n_x \right)dS,
\end{align*}
where $S$ is the cylinder, $p(t)$ is pressure, $u_{t_S}$ is the tangential velocity, the length scale is selected to be the diameter of the cylinder, and $n = \langle n_x, n_y\rangle$ is the outward unit normal.  

We run multiple experiments with this test problem.  
In section~\ref{sec:channel-emac-fem}, we use EMAC-FEM as the FOM and then run EMAC-ROM, SKEW-ROM and CONV-ROM and compare results (EMAC-ROM is superior).  
In section~\ref{sec:channel-conv-fem}, we repeat this test but using CONV-FEM as the FOM (CONV-ROM).  
Finally, in section~\ref{sec:channel-consistency}, 
we study convergence of the FOM to the ROM, as the number of modes $N$ increases.  We show that for EMAC-FEM as the FOM, EMAC-ROM converges to the FOM as $N$ increases, but SKEW-ROM does not converge for increasing $N$ and appears to lock (as predicted by our theory above).  

\subsubsection{Comparison of ROMs using EMAC-FEM as FOM}
    \label{sec:channel-emac-fem}

\begin{figure}
\centering
\includegraphics[width=.6\textwidth,height=0.22\textwidth]{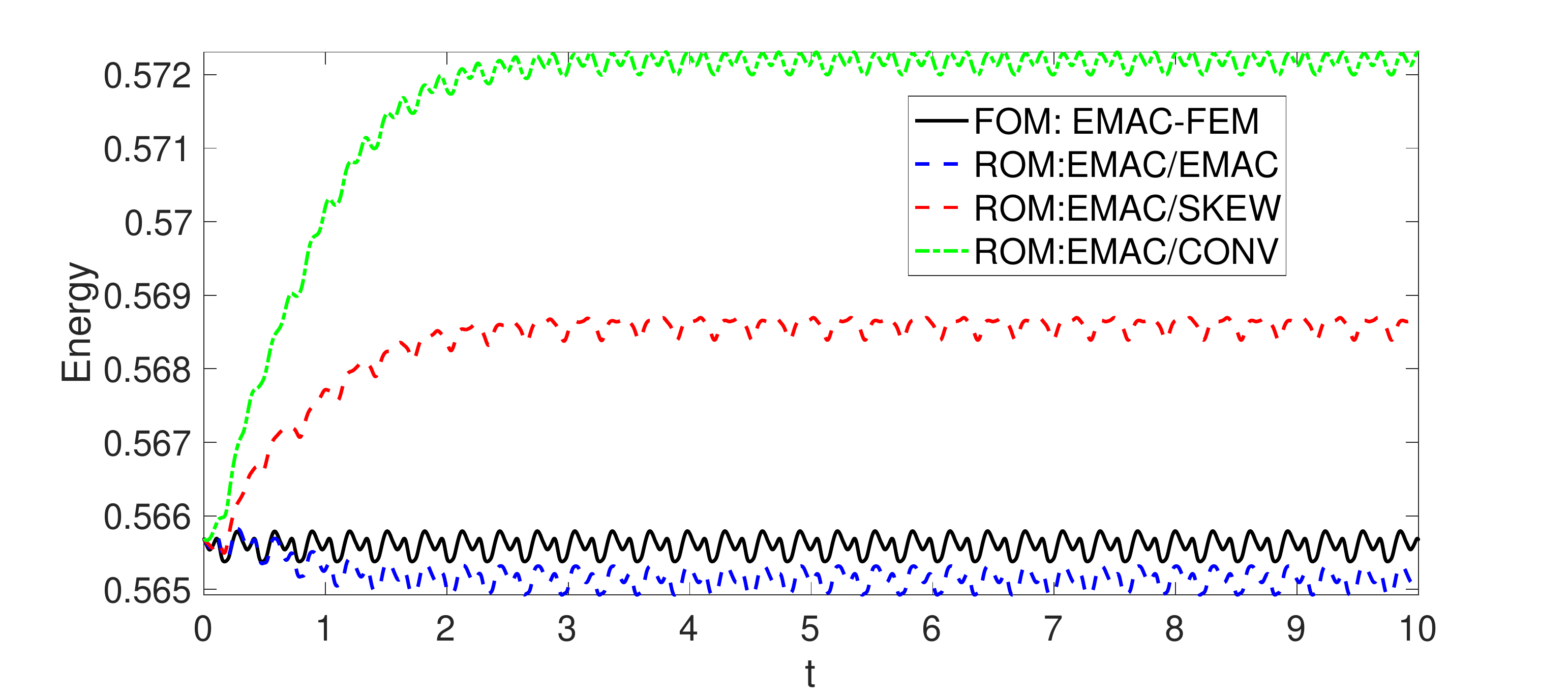}
\includegraphics[width=.3\textwidth,height=0.22\textwidth]{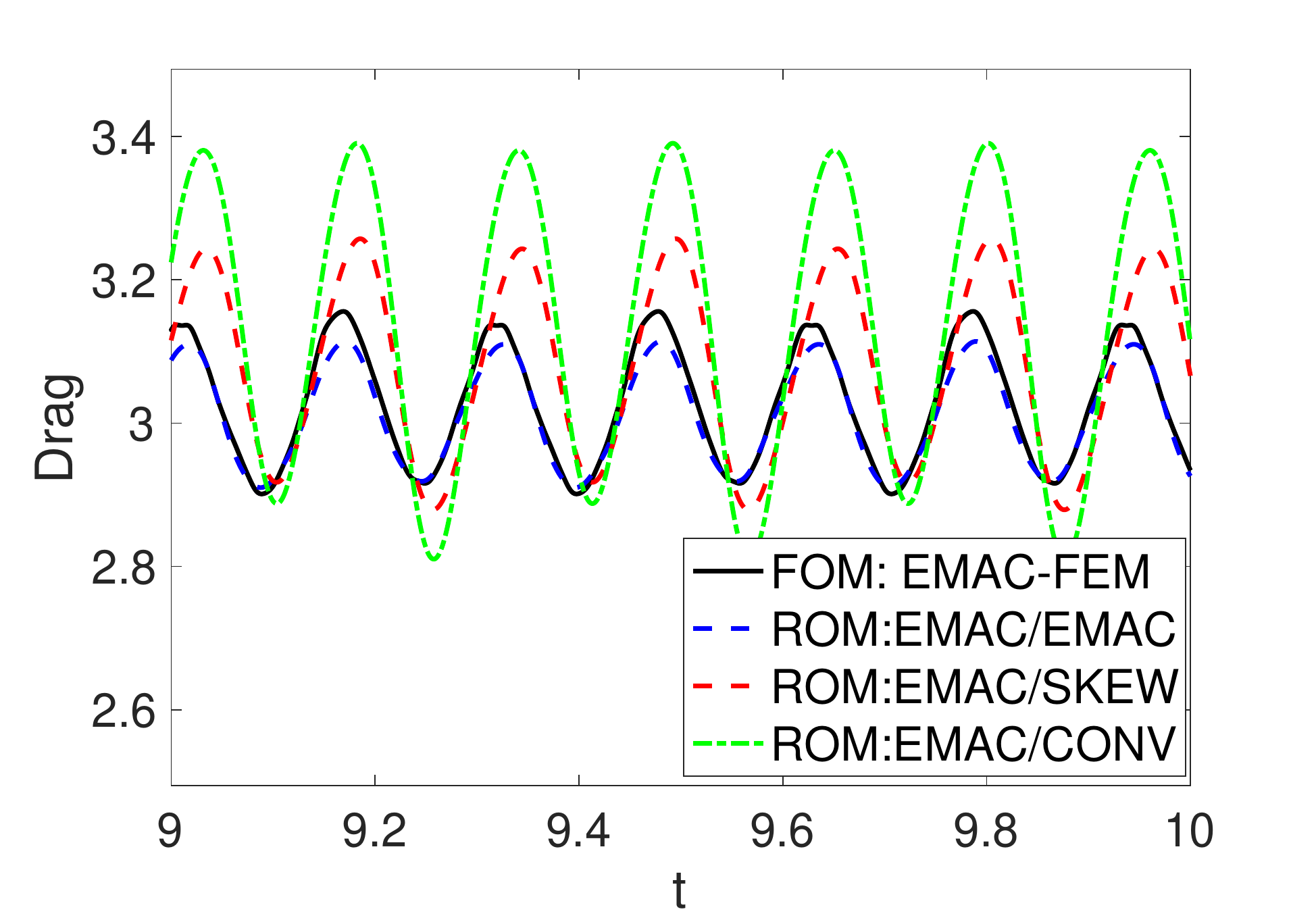}
\caption{Shown above are results of ROM simulations built from EMAC-FEM as the FOM, 
using EMAC-ROM, SKEW-ROM, and CONV-ROM with 13 modes.  
\label{Cyl200EMAC}}
\end{figure}

For our first test, we choose EMAC-FEM as the FOM to build the ROM space, and then run EMAC-ROM, SKEW-ROM and CONV-ROM, each with $N=13$ modes.  Results of each run's energy versus time and drag versus time are shown in Figure \ref{Cyl200EMAC} (drag is shown only on $9\le t\le 10$ so it is easier to see and compare, but differences are similar on $0\le t\le 9$).  From the plots, it is clear the EMAC-ROM (labeled ROM: EMAC/EMAC) is the most accurate for both energy and drag prediction, in terms of best matching the FOM energy and drag predictions.  Thus, this is an excellent illustration of the theory above that shows that ROMs that are consistent with FOMs give better results.

\subsubsection{Comparison of ROMs using CONV-FEM as FOM}
    \label{sec:channel-conv-fem}

\begin{figure}[h]
\centering
\includegraphics[width=.6\textwidth,height=0.22\textwidth]{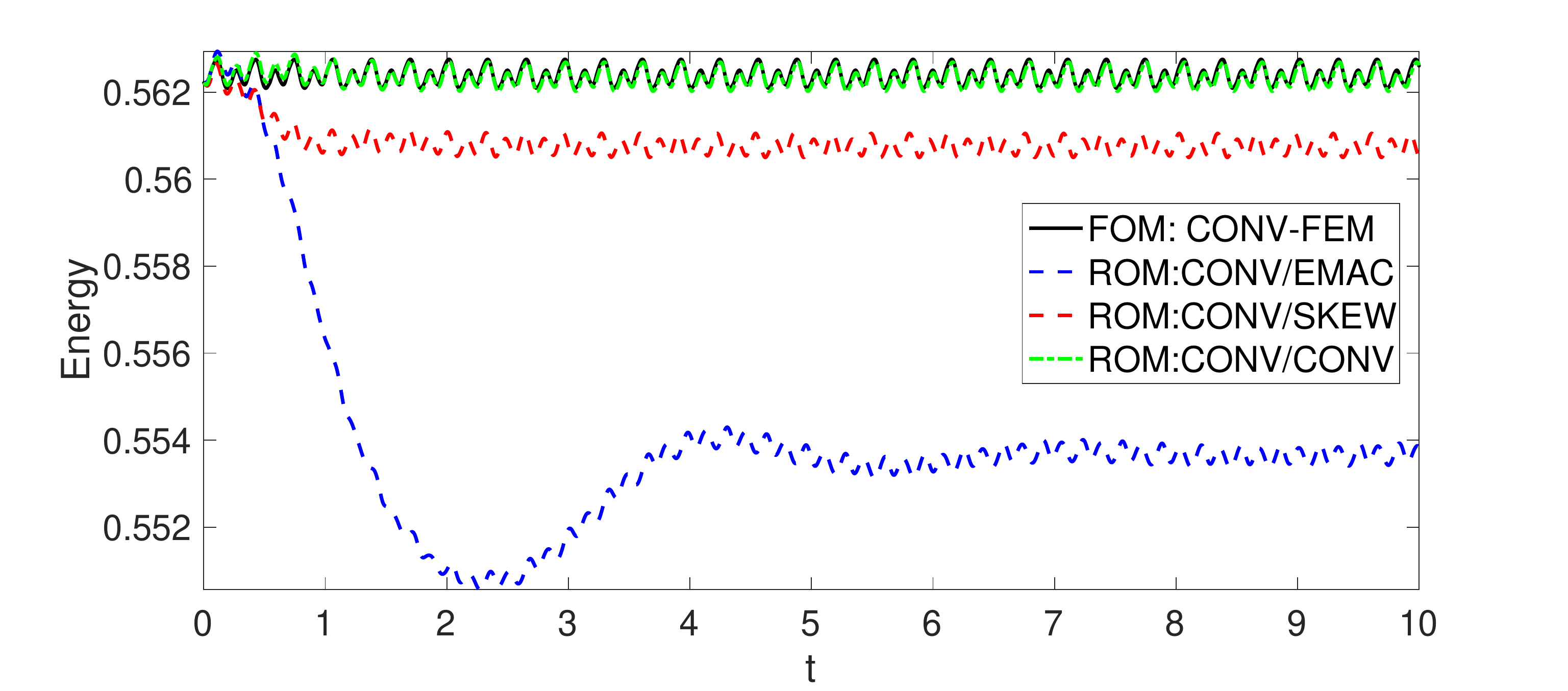}
\includegraphics[width=.3\textwidth,height=0.22\textwidth]{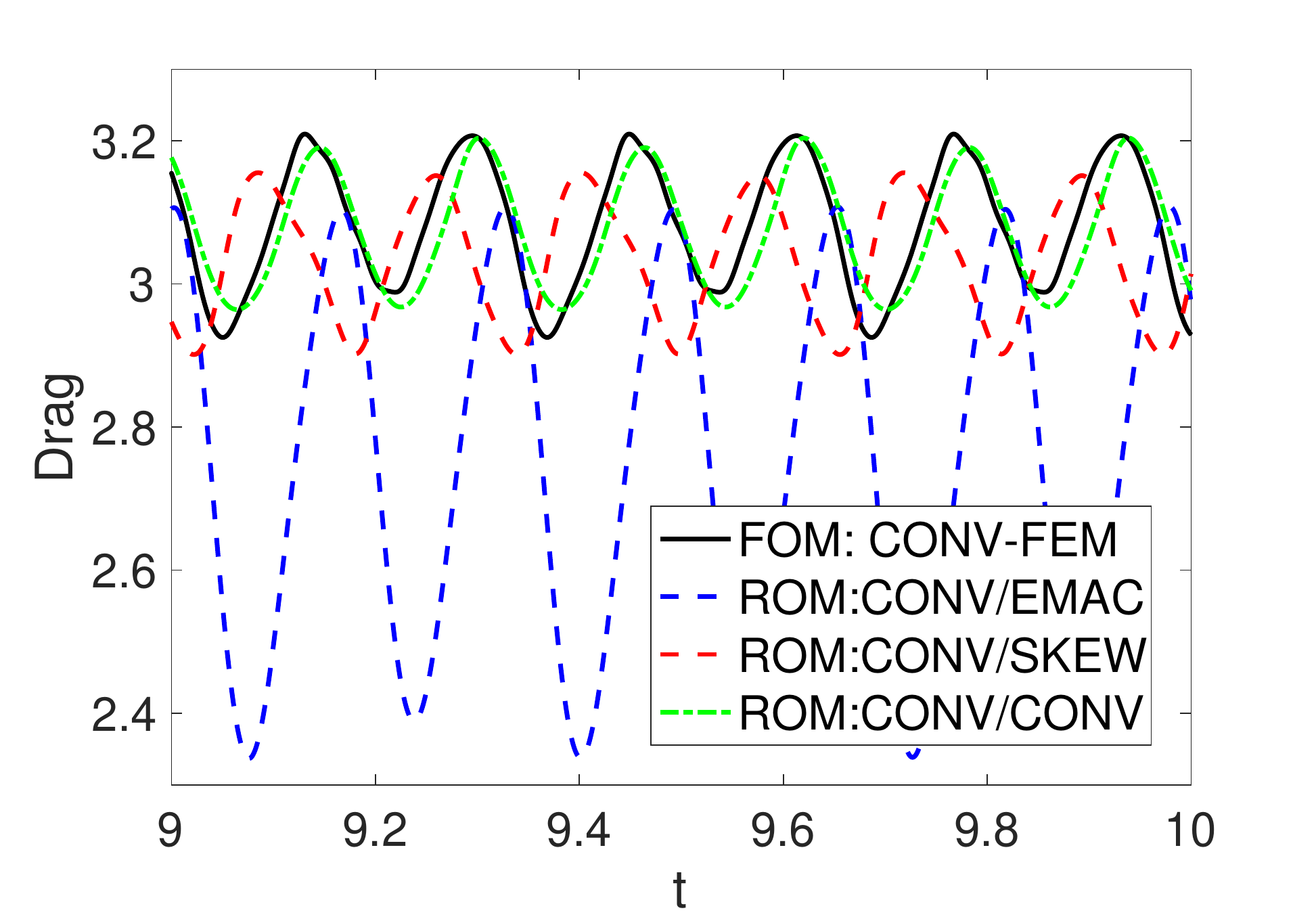}
\caption{Shown above are results of ROM simulations built from 
CONV-FEM as the FOM, 
using EMAC-ROM, SKEW-ROM, and CONV-ROM with 13 modes. }
\label{Cyl200CONV}
\end{figure}

For our second test, we choose CONV-FEM as the FOM to build the ROM space, and then run EMAC-ROM, SKEW-ROM and CONV-ROM, each with $N=13$ modes.  Results of each run's energy versus time and drag versus time are shown in Figure \ref{Cyl200CONV}.  From the plots, it is clear the CONV-ROM is the most accurate for both energy and drag prediction, in terms of best matching the FOM energy and drag predictions.  This is another illustration of the theory above that shows that ROMs that are consistent with FOMs give better results.

\subsubsection{Convergence of ROM to FOM depends on FOM/ROM consistency}
    \label{sec:channel-consistency}

Theorem \ref{FOMROMerr} gives a bound that guarantees that a {\it consistent} ROM will converge to the FOM as the number of modes $N$ increases.  \LR{ The bound also suggests that an inconsistent ROM may only converge to the FOM up to the terms depending on the divergence error in the FOM. This is analogous to locking phenomena in convergence theory for numerical methods, where as discretization parameters get small (and number of degrees of freedom gets large), the error cannot decrease below a certain level.}  We now illustrate this theory.  

\begin{figure}
\centering
\includegraphics[width=.6\textwidth,height=0.22\textwidth]{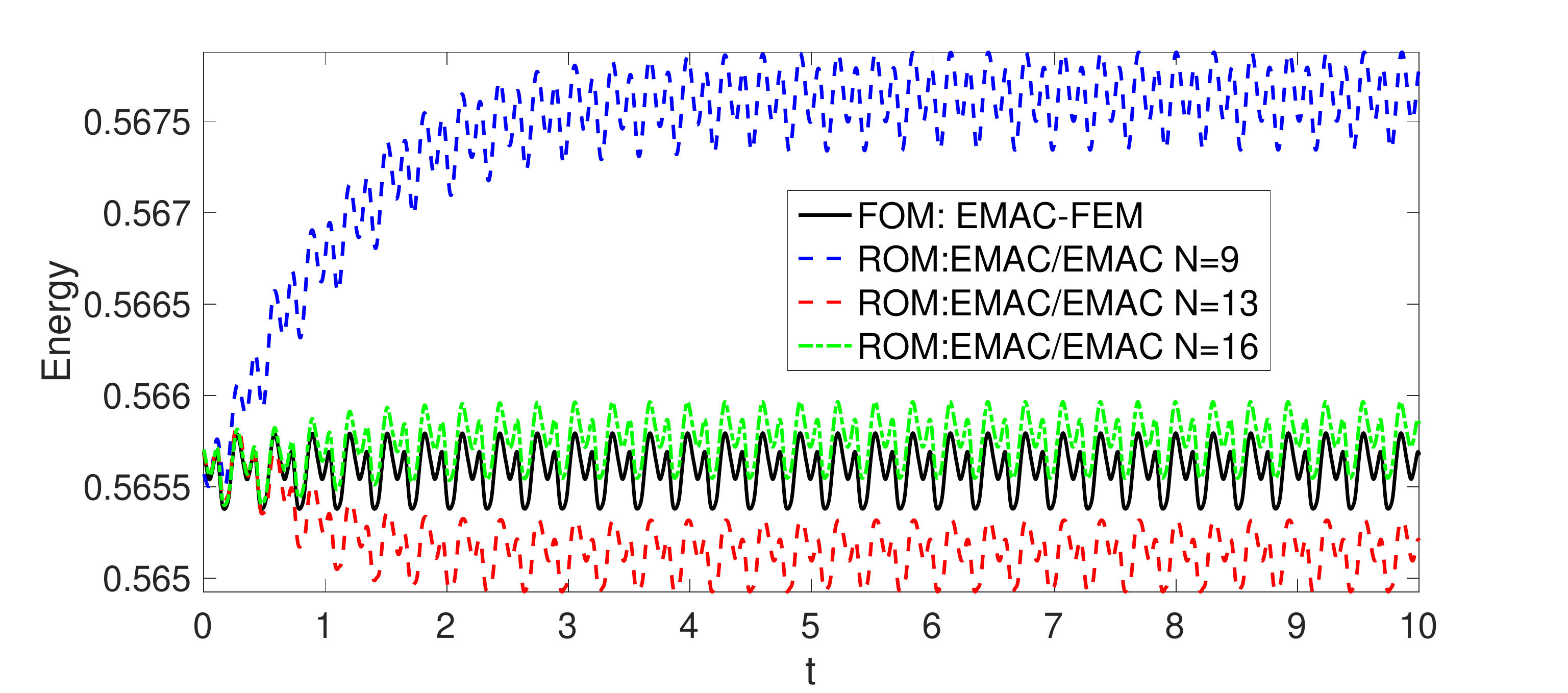}
\includegraphics[width=.3\textwidth,height=0.22\textwidth]{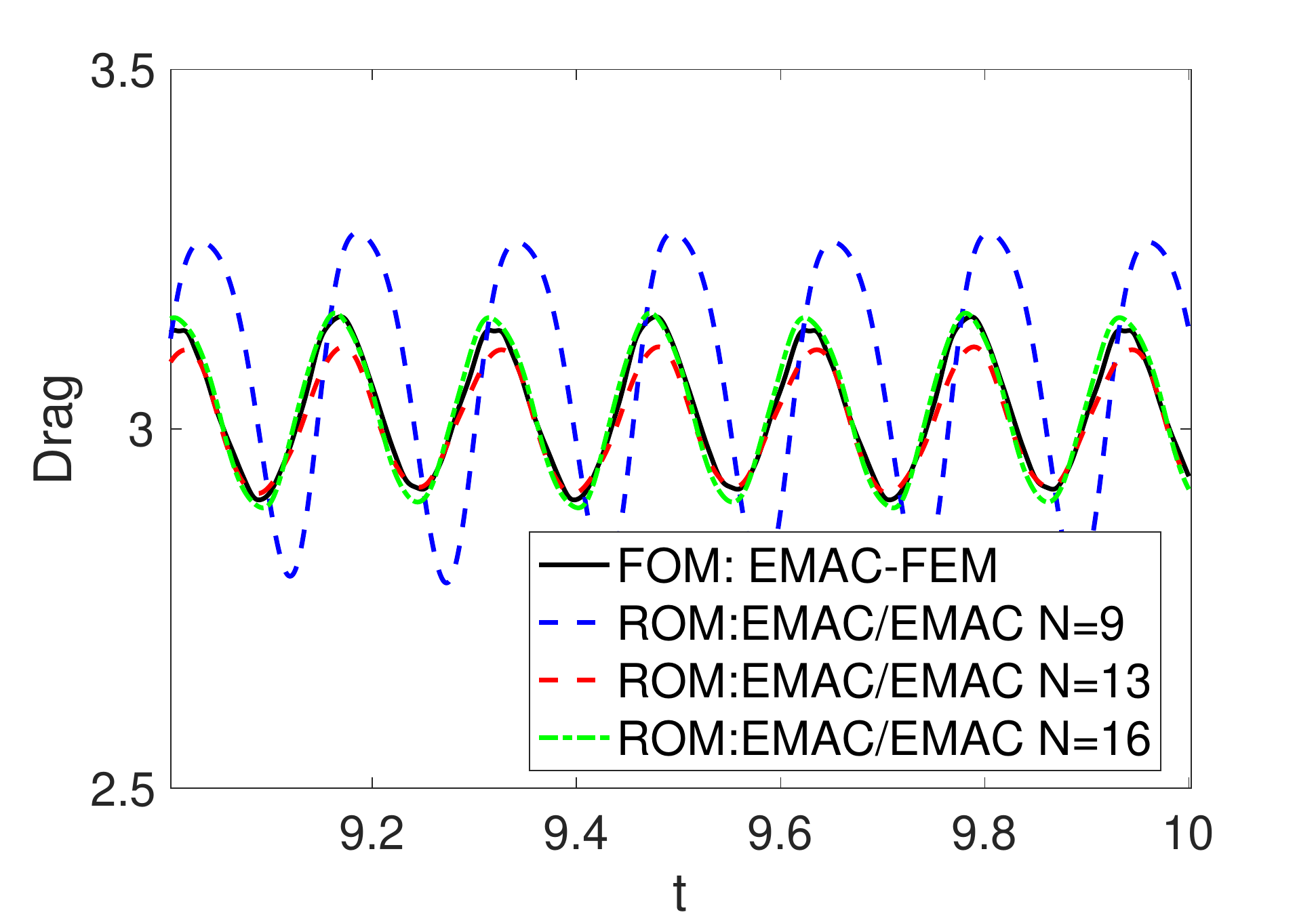} \\ \ \\
\includegraphics[width=.6\textwidth,height=0.22\textwidth]{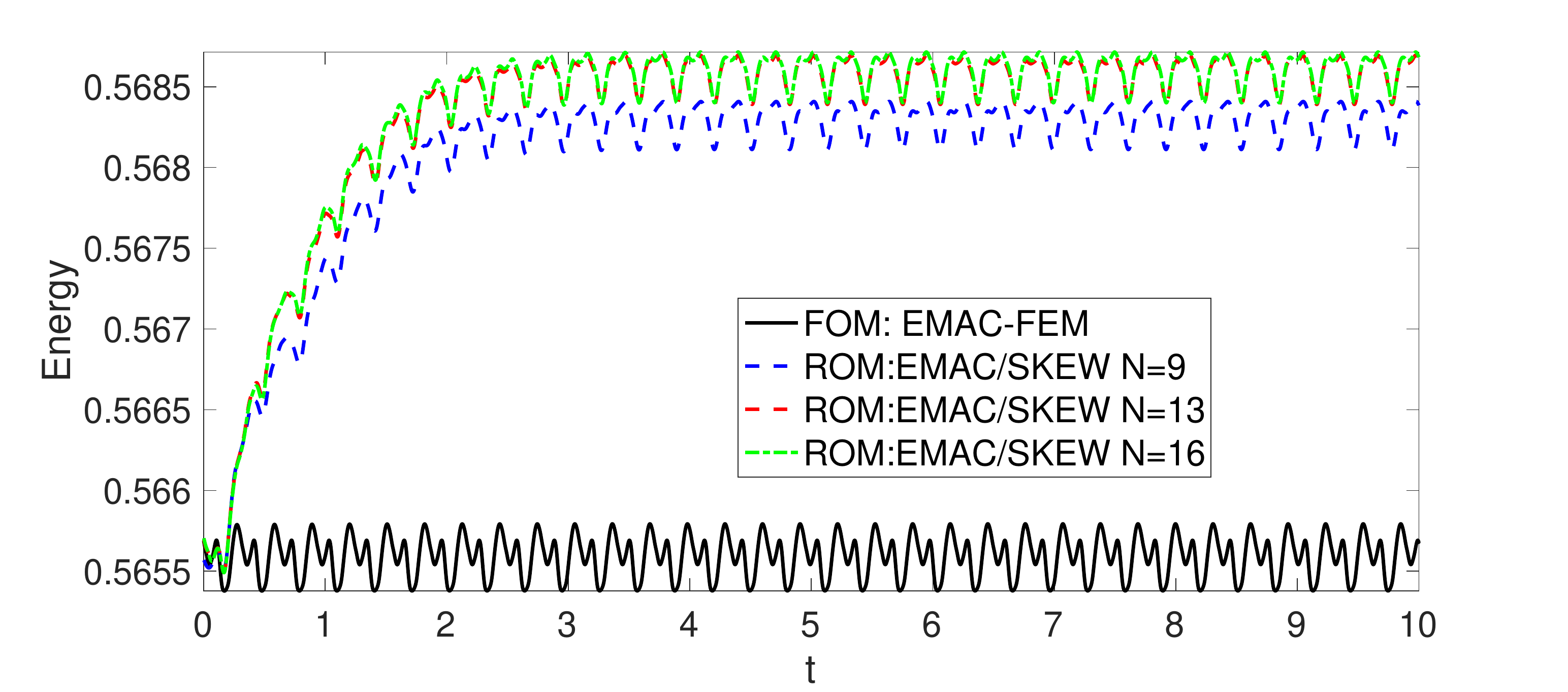}
\includegraphics[width=.3\textwidth,height=0.22\textwidth]{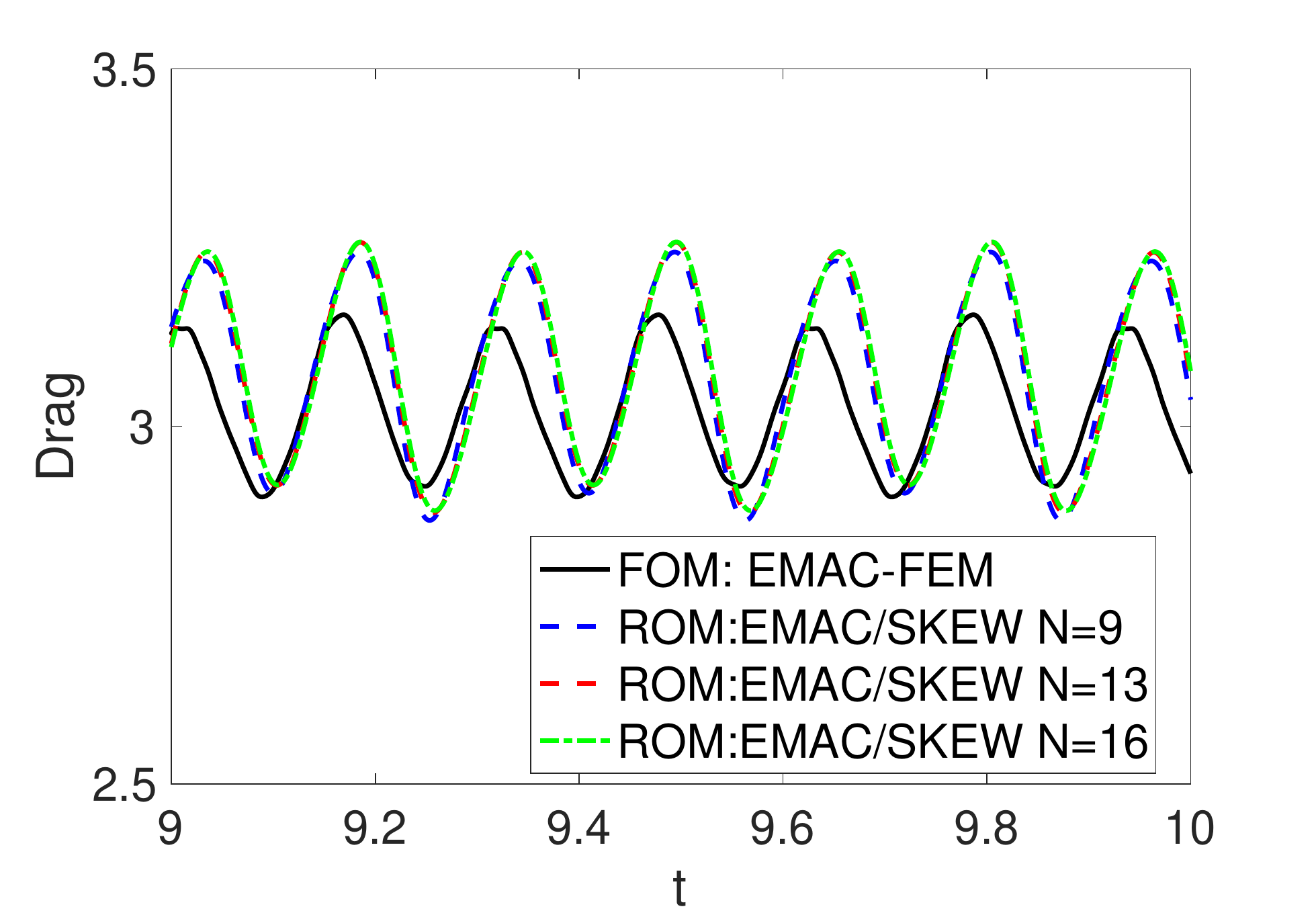}
\caption{Shown above are results of ROM simulations built from EMAC-FEM as the FOM, 
and (top) EMAC-ROM and (bottom) SKEW-ROM with $N$=9, 13, and 16 modes. }
\label{Cyl200EMACcv}
\end{figure}

First, we consider the consistent case, and choose EMAC-FEM for the FOM, and EMAC-ROM as the ROM, and run the ROM using $N$=9, 13, and 16 modes.  From the plots in Figure \ref{Cyl200EMACcv} at the top row, we observe convergence of the energy and drag to the FOM as $N$ increases.

For the inconsistent case, we repeat the same test, again choosing EMAC-FEM as the FOM, but now 
choosing SKEW-ROM for the ROMs.  The plots in Figure \ref{Cyl200EMACcv} at the bottom row
do not show convergence to the FOM, but seem to show convergence to a different limit (since $N$=13 and 16 plots are nearly identical).  This is in agreement with the theory that FOM/ROM inconsistency can prevent convergence of the ROM to the FOM as $N$ increases.

\subsection{Kelvin-Helmholtz instability}
    \label{sec:kh}

For our second numerical experiment, we consider the 2D Kelvin-Helmholtz instability problem from \cite{SJLLLS18} with $Re$=100.  The domain is $(0,1)^2$, with periodic
boundary conditions on the sides, and on the top and bottom we enforce a no penetration condition along with a natural weak enforcement of the free-slip condition $(-\nu\nabla u \cdot n)\times n=0$.  The initial condition is given by

\[
u_0(x,y) = \left( \begin{array}{c} \tanh\left( 28 (2y-1) \right) \\ 0 \end{array} \right) + 10^{-3} \left( \begin{array}{c} \partial_y \psi(x,y) \\ -\partial_x \psi(x,y) \end{array} \right)
\]
with 
\[
\psi =  \exp \left( -28^2 (y-0.5)^2 \right) \left( \cos(8\pi x) + \cos(20\pi x) \right),
\]
where $\frac{1}{28}$ is the initial vorticity thickness and $10^{-3}$ is a noise/scaling factor.  The Reynolds number is defined by $Re=\frac{1}{28 \nu}$, there is no forcing ($f=0$), and $\nu$ is defined by selecting $Re$=100.  As discussed in  \cite{OR20,SJLLLS18}, this is a very sensitive problem where small changes in parameters can lead to differences in how the eddies combine, how long they take to combine, where the final eddy is located, and so on.  This is a particularly 
challenging problem for reduced order modeling, since the solution evolves in time and is never statistically steady.

We proceed with similar tests as in the 2D cylinder problem in section~\ref{sec:channel}.  That is, we use certain formulations to create the FOM, and the run varying ROMs and test for accuracy.  
\subsubsection{Comparison of ROMs using EMAC-FEM as FOM}
    \label{sec:kh-emac-fem}

We first consider using EMAC-FEM as the FOM, and run EMAC-ROM, SKEW-ROM, and CONV-ROM.  EMAC-FEM \LR{ is calculated up to $t=5$ using $\Delta t=0.01$ and $h=1/96$ on a uniform triangulation of the domain}.  Snapshots are taken at each time step, and the ROM space is created using the top $N$=50 modes.  The ROMs are run starting from the same $u_0$ (projected into $V_r$) up to $t=5$ also using $\Delta t=0.01$.

\begin{figure}[!h]
\begin{center}
\includegraphics[width=.45\textwidth, height=.32\textwidth,viewport=0 0 520 400, clip]{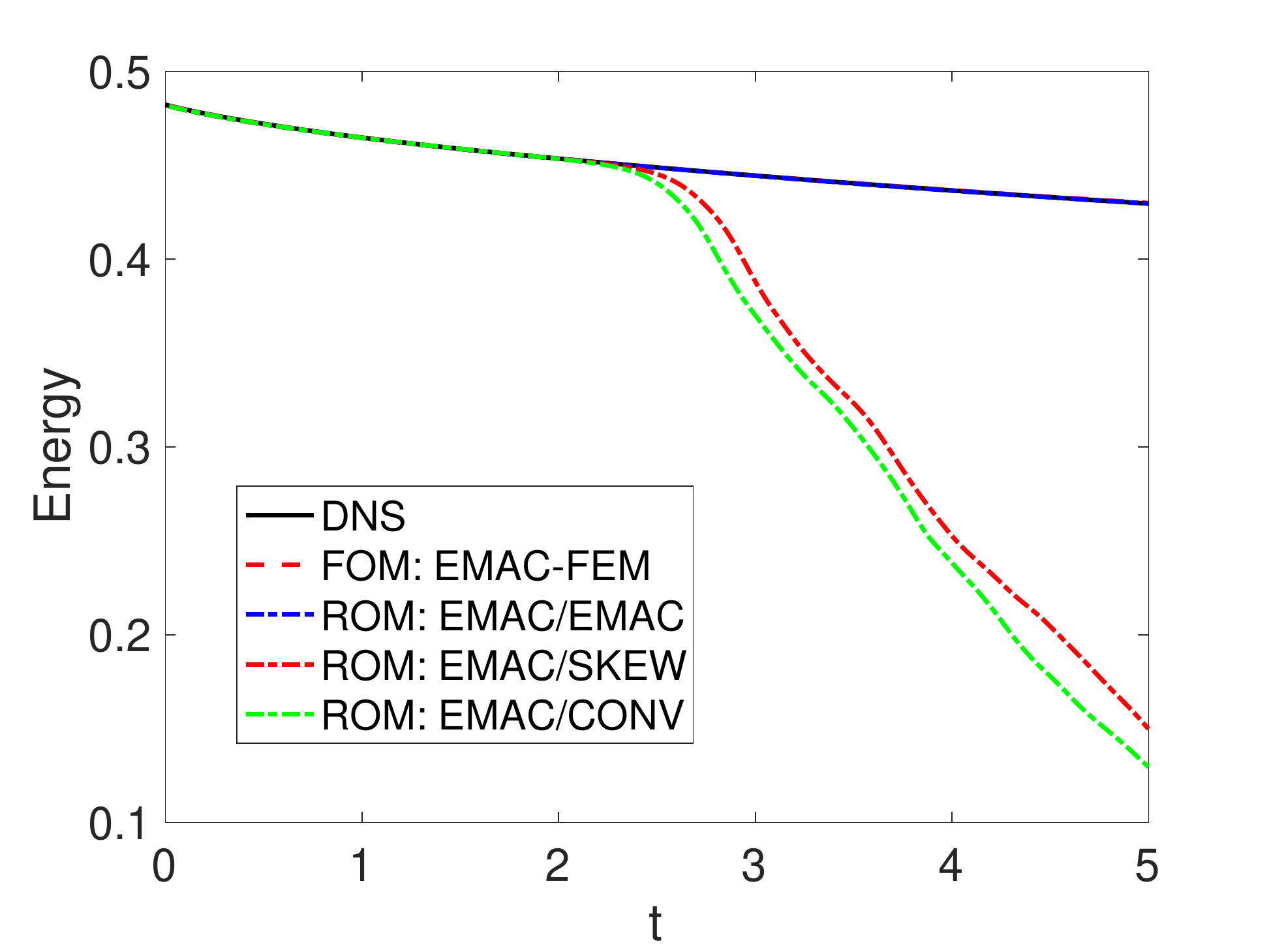}
\includegraphics[width=.45\textwidth, height=.32\textwidth,viewport=0 0 520 400, clip]{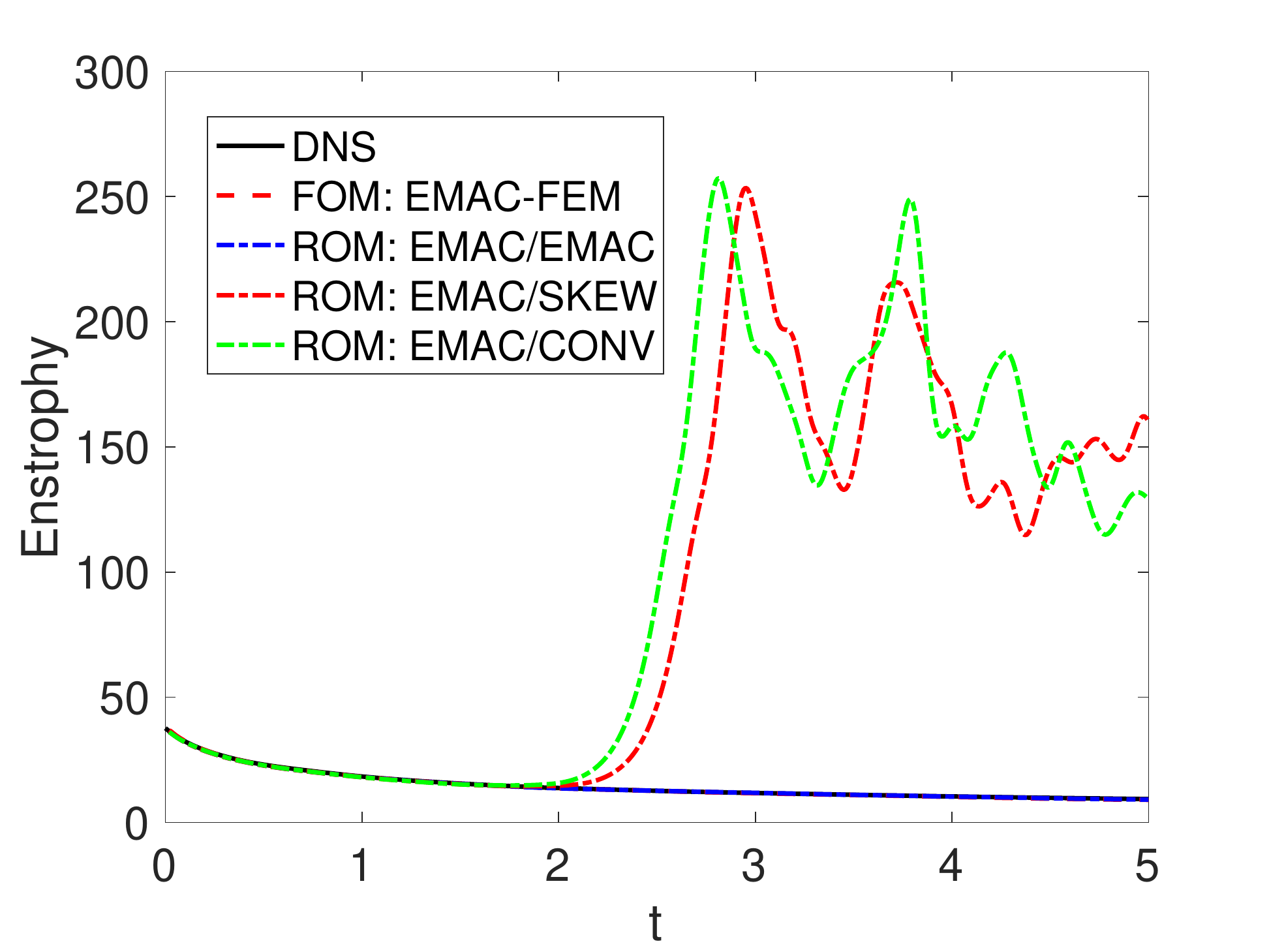} 
\end{center}
\caption{\label{emacEE}
Shown above are the energy and enstrophy plots versus time for the $Re=100$ KH tests for varying ROMs (that are built from EMAC-FEM), DNS, the FOM (EMAC-FEM), and a DNS from  \cite{SJLLLS18}.
}
\end{figure}

\begin{figure}[!h]
\hspace{.2in} SKEW-ROM \hspace{.6in} EMAC-ROM \hspace{.7in} CONV-ROM \hspace{.4in} FOM (EMAC-FEM) 
\begin{center}
\includegraphics[width=.24\textwidth, height=.12\textwidth,viewport=0 0 520 400, clip]{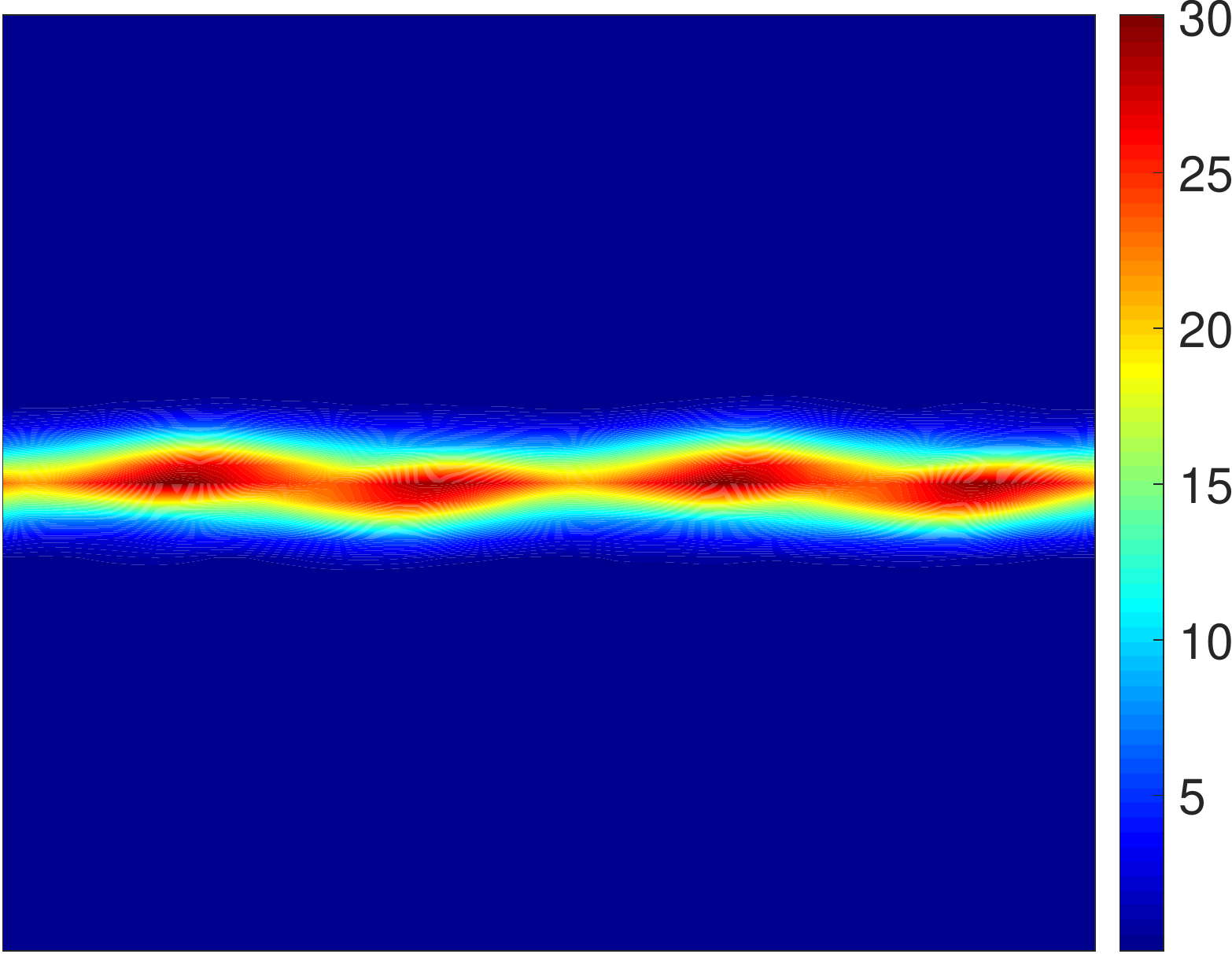}
\includegraphics[width=.24\textwidth, height=.12\textwidth,viewport=0 0 520 400, clip]{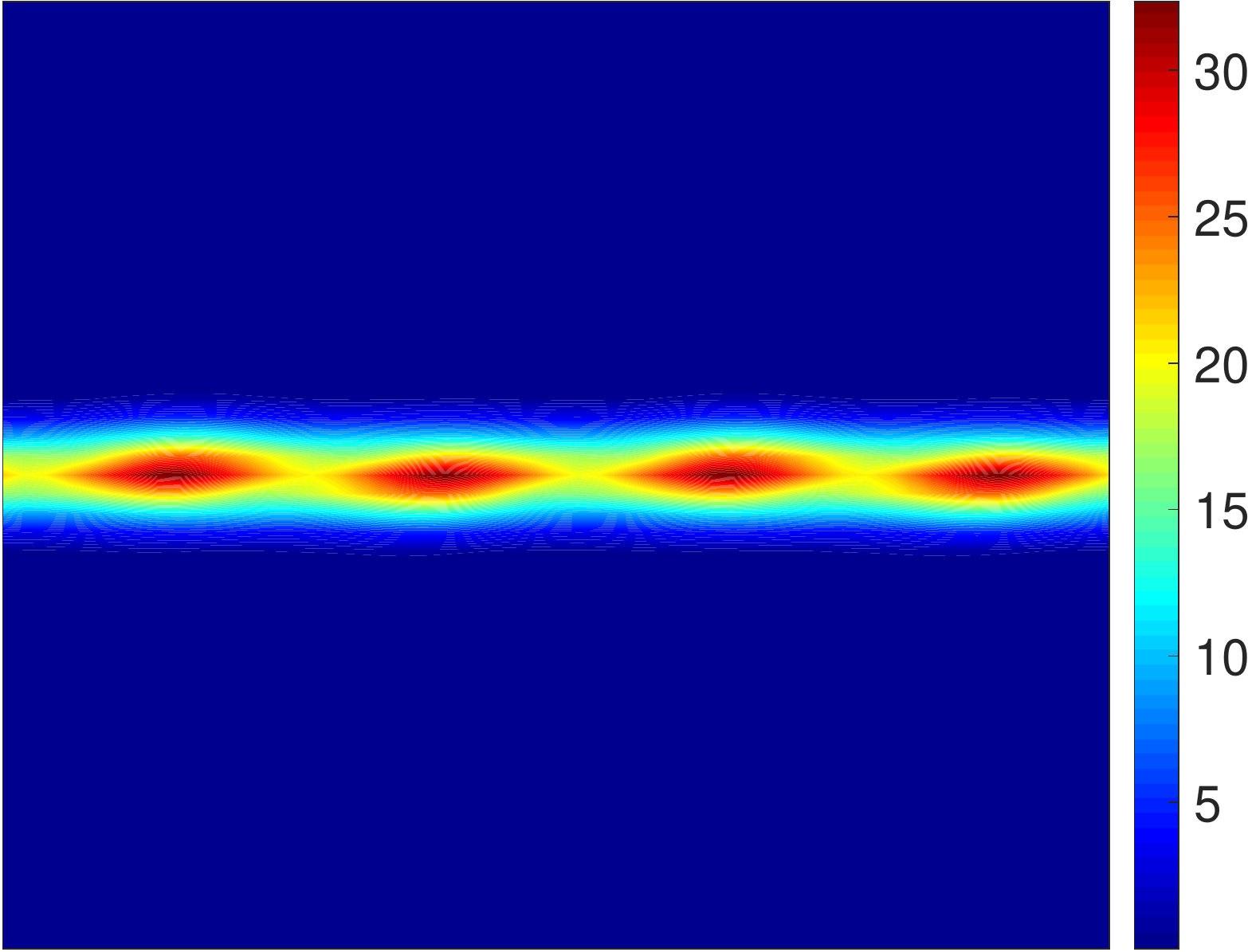} 
\includegraphics[width=.24\textwidth, height=.12\textwidth,viewport=0 0 520 400, clip]{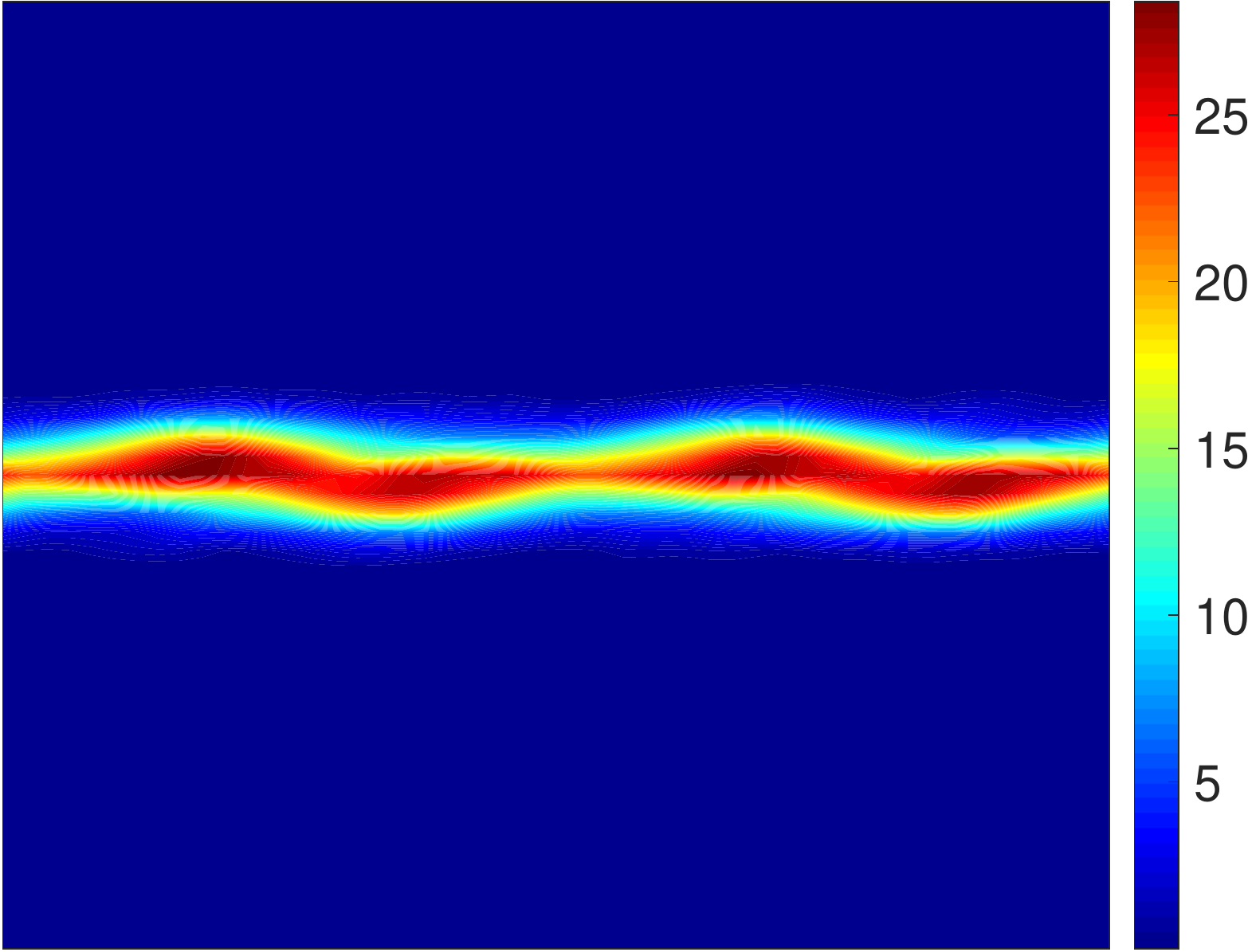}
\includegraphics[width=.24\textwidth, height=.12\textwidth,viewport=0 0 520 400, clip]{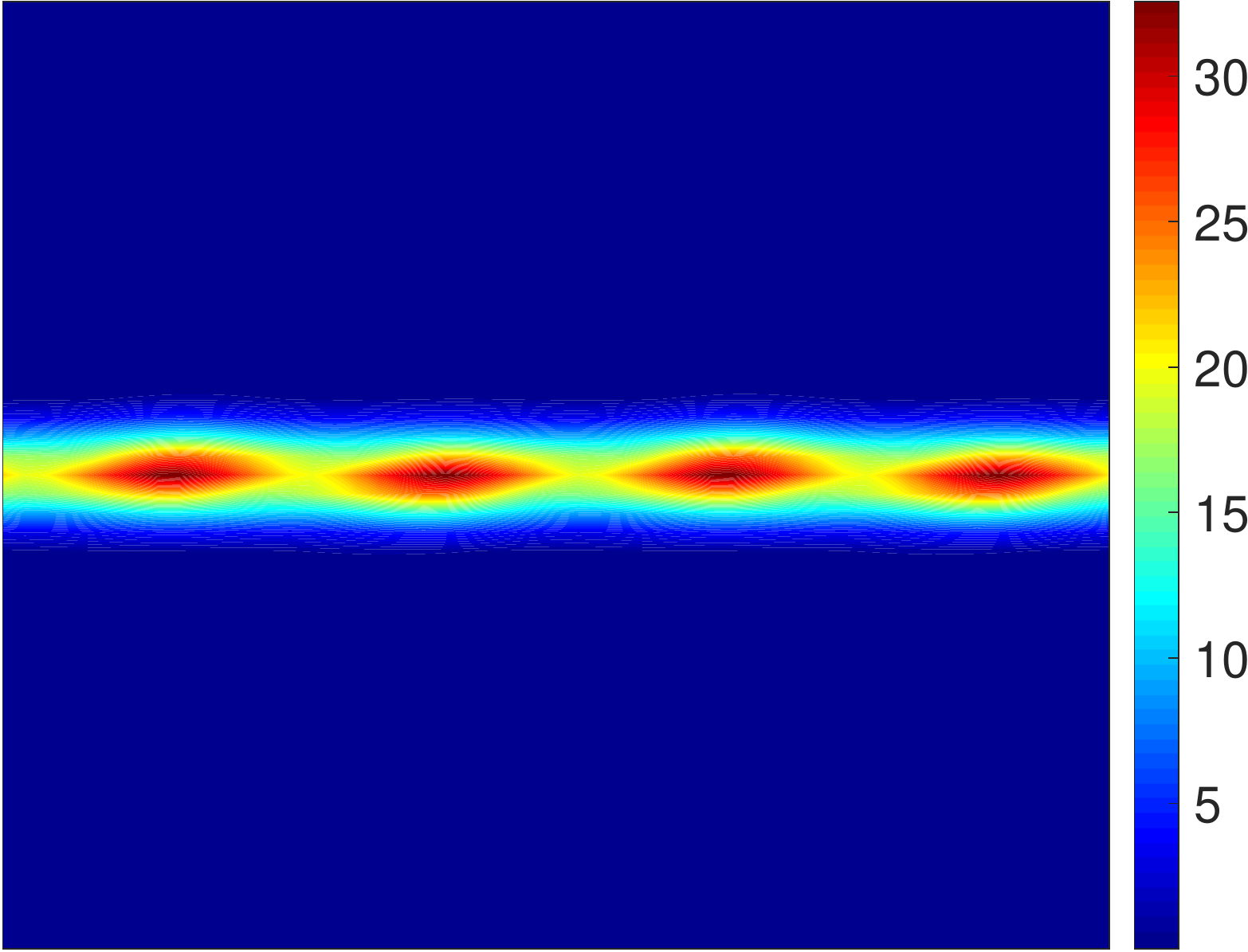}\\
\includegraphics[width=.24\textwidth, height=.12\textwidth,viewport=0 0 520 400, clip]{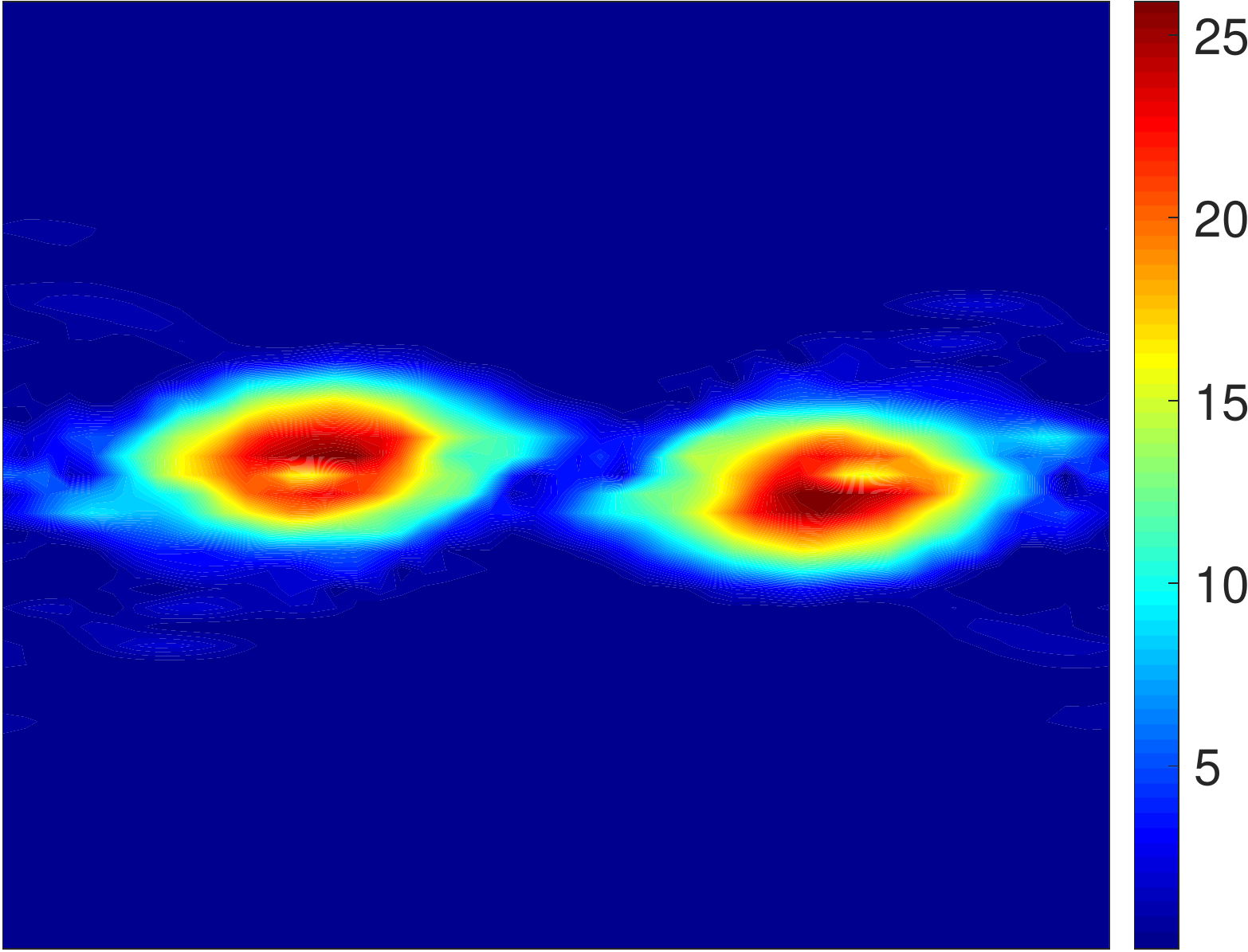}
\includegraphics[width=.24\textwidth, height=.12\textwidth,viewport=0 0 520 400, clip]{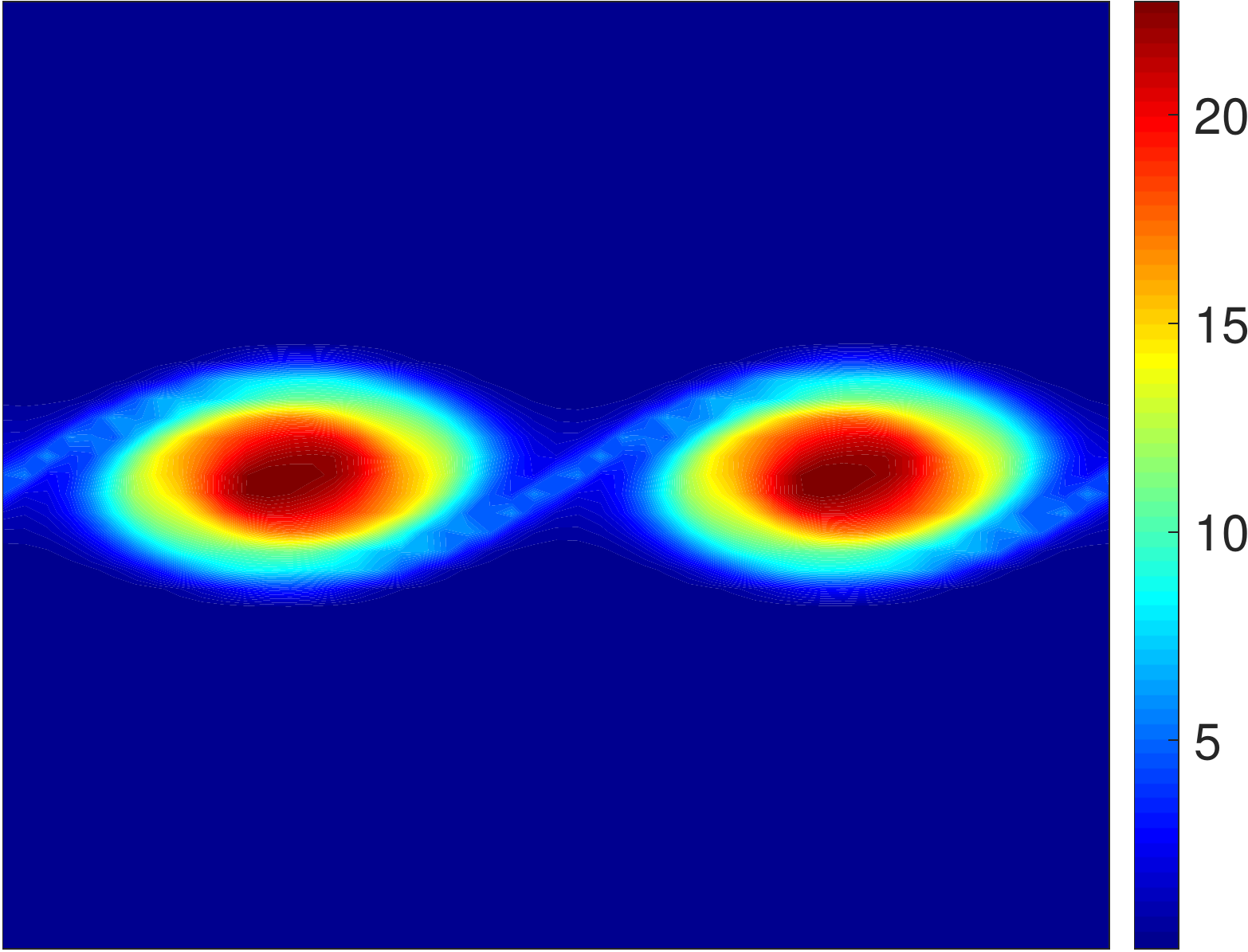} 
\includegraphics[width=.24\textwidth, height=.12\textwidth,viewport=0 0 520 400, clip]{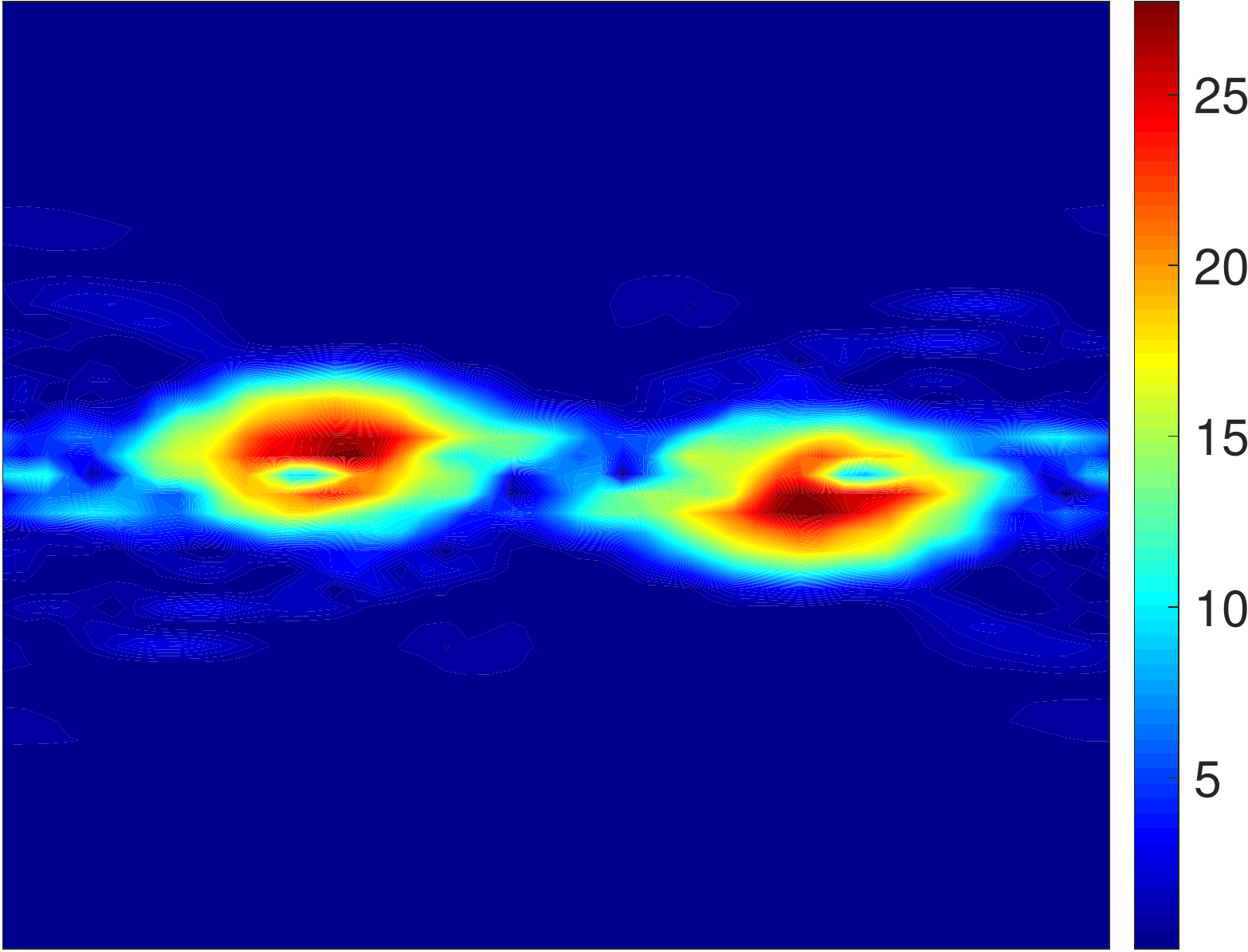}
\includegraphics[width=.24\textwidth, height=.12\textwidth,viewport=0 0 520 400, clip]{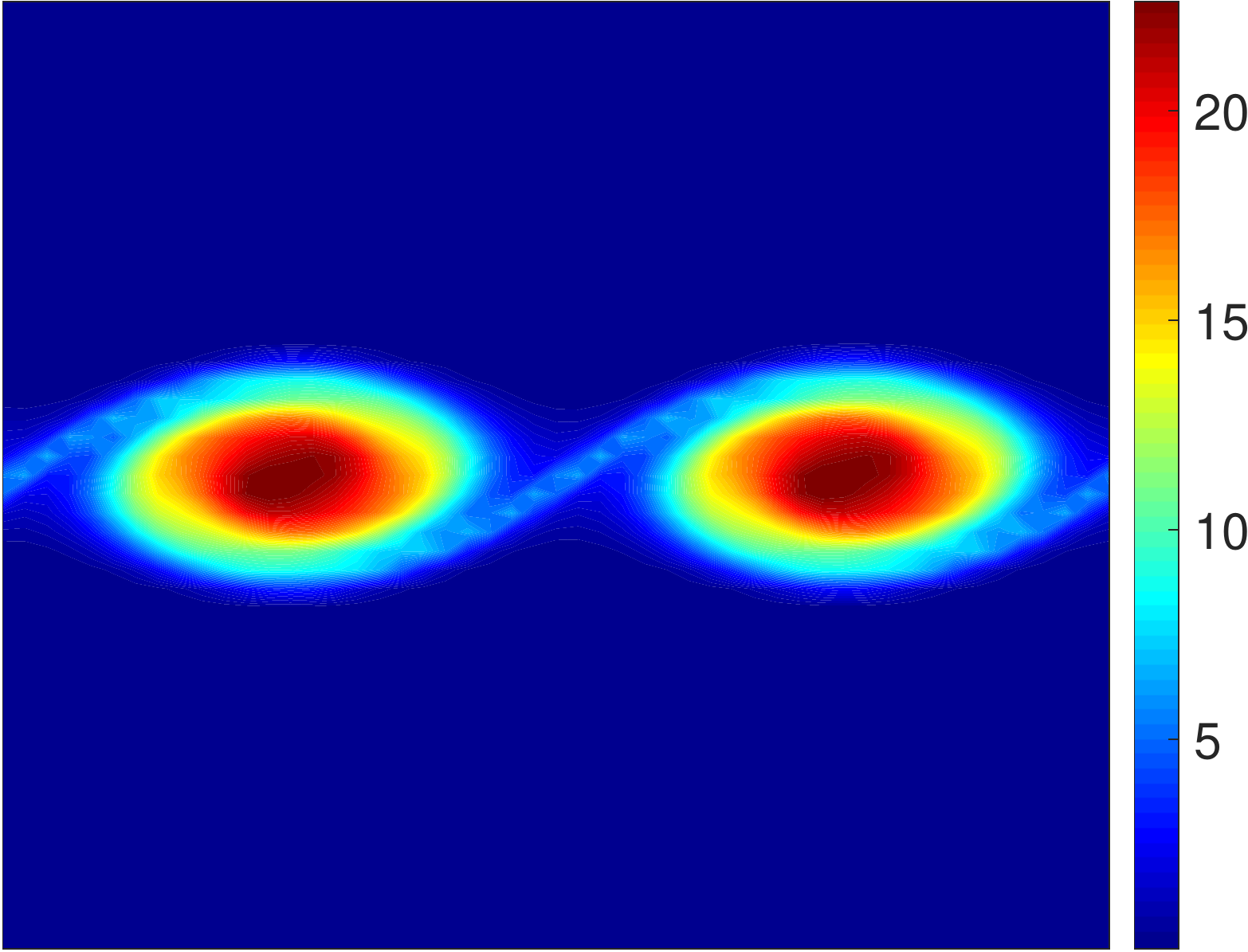}\\
\includegraphics[width=.24\textwidth, height=.12\textwidth,viewport=0 0 520 400, clip]{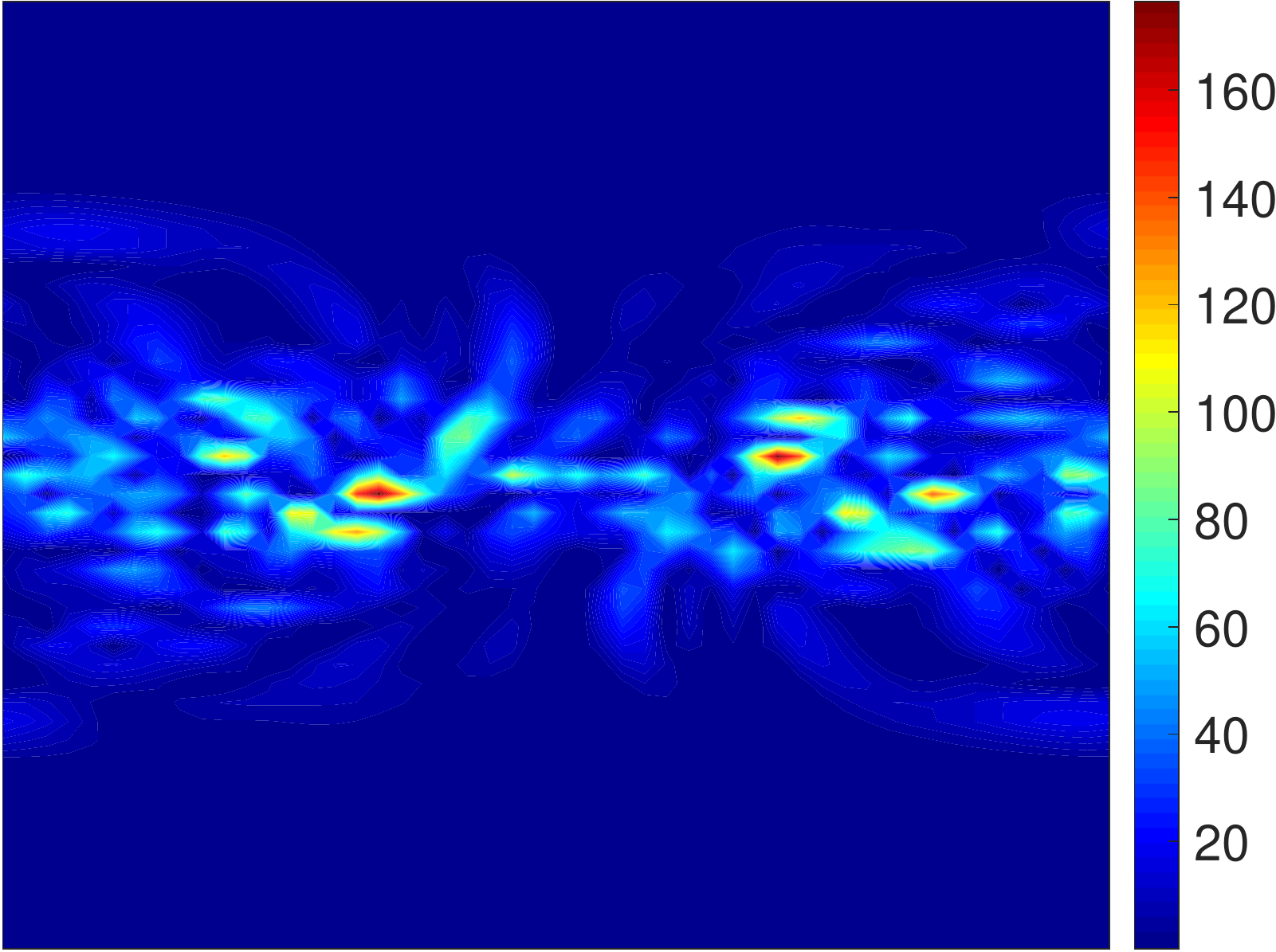}
\includegraphics[width=.24\textwidth, height=.12\textwidth,viewport=0 0 520 400, clip]{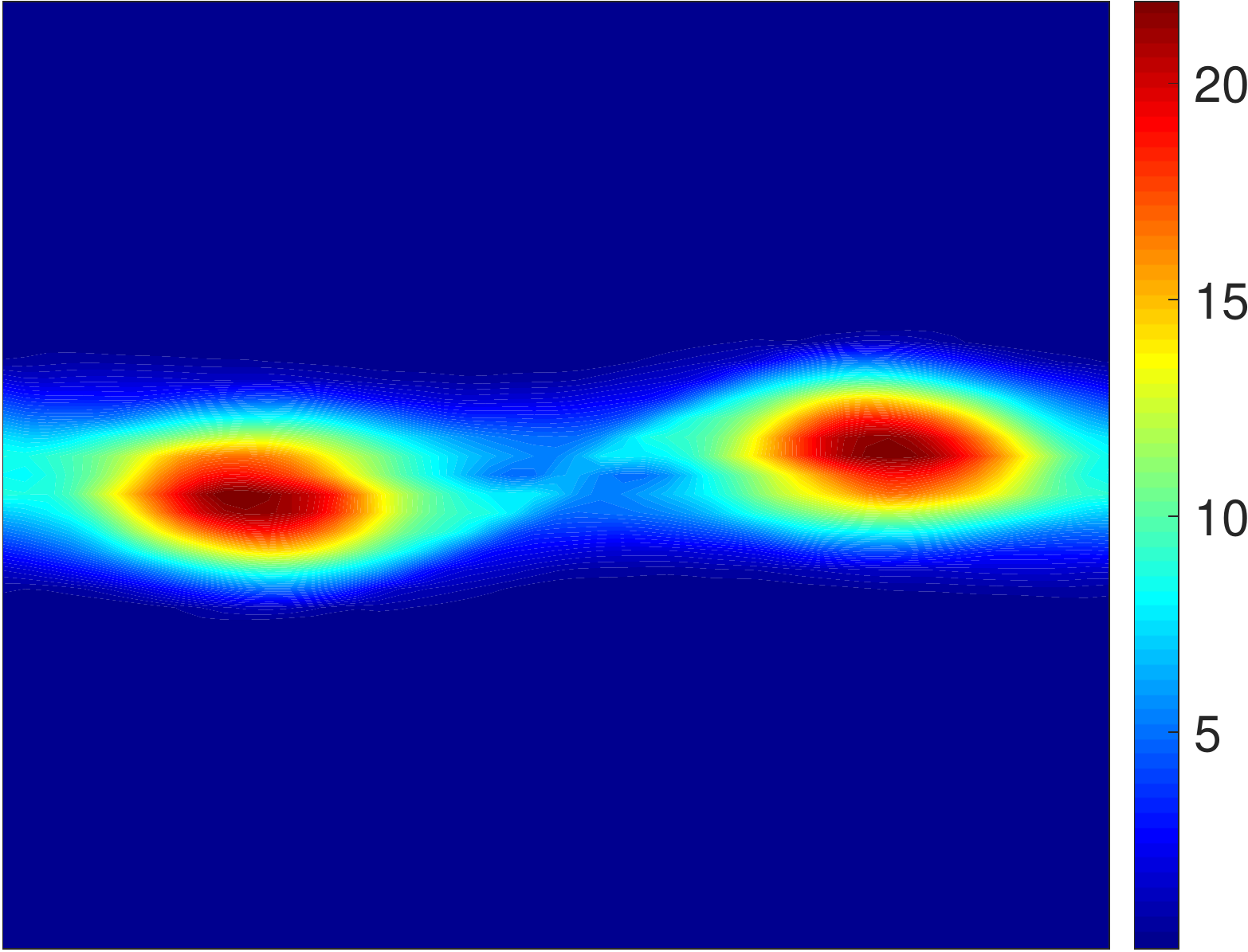} 
\includegraphics[width=.24\textwidth, height=.12\textwidth,viewport=0 0 520 400, clip]{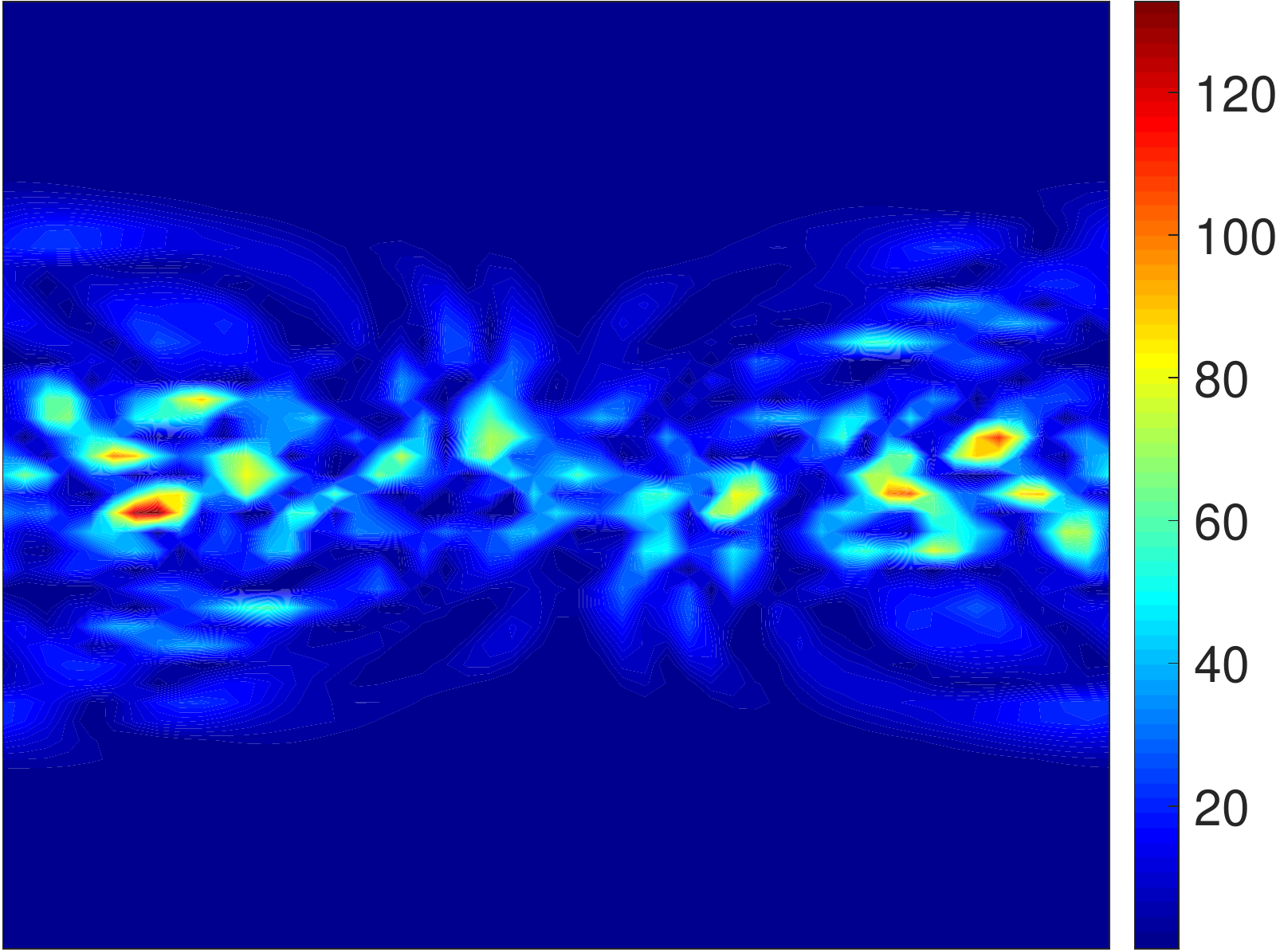}
\includegraphics[width=.24\textwidth, height=.12\textwidth,viewport=0 0 520 400, clip]{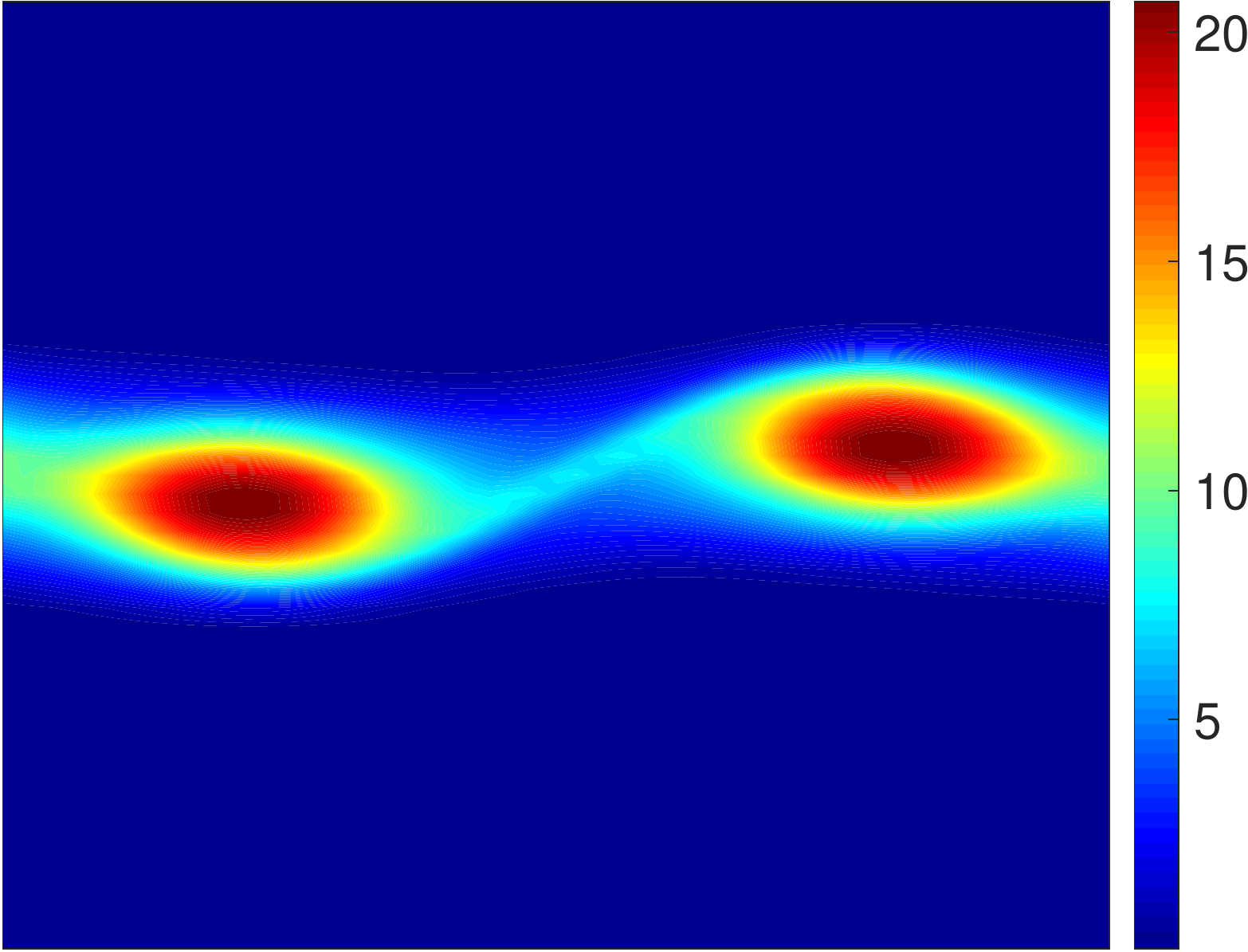}\\
\includegraphics[width=.24\textwidth, height=.12\textwidth,viewport=0 0 520 400, clip]{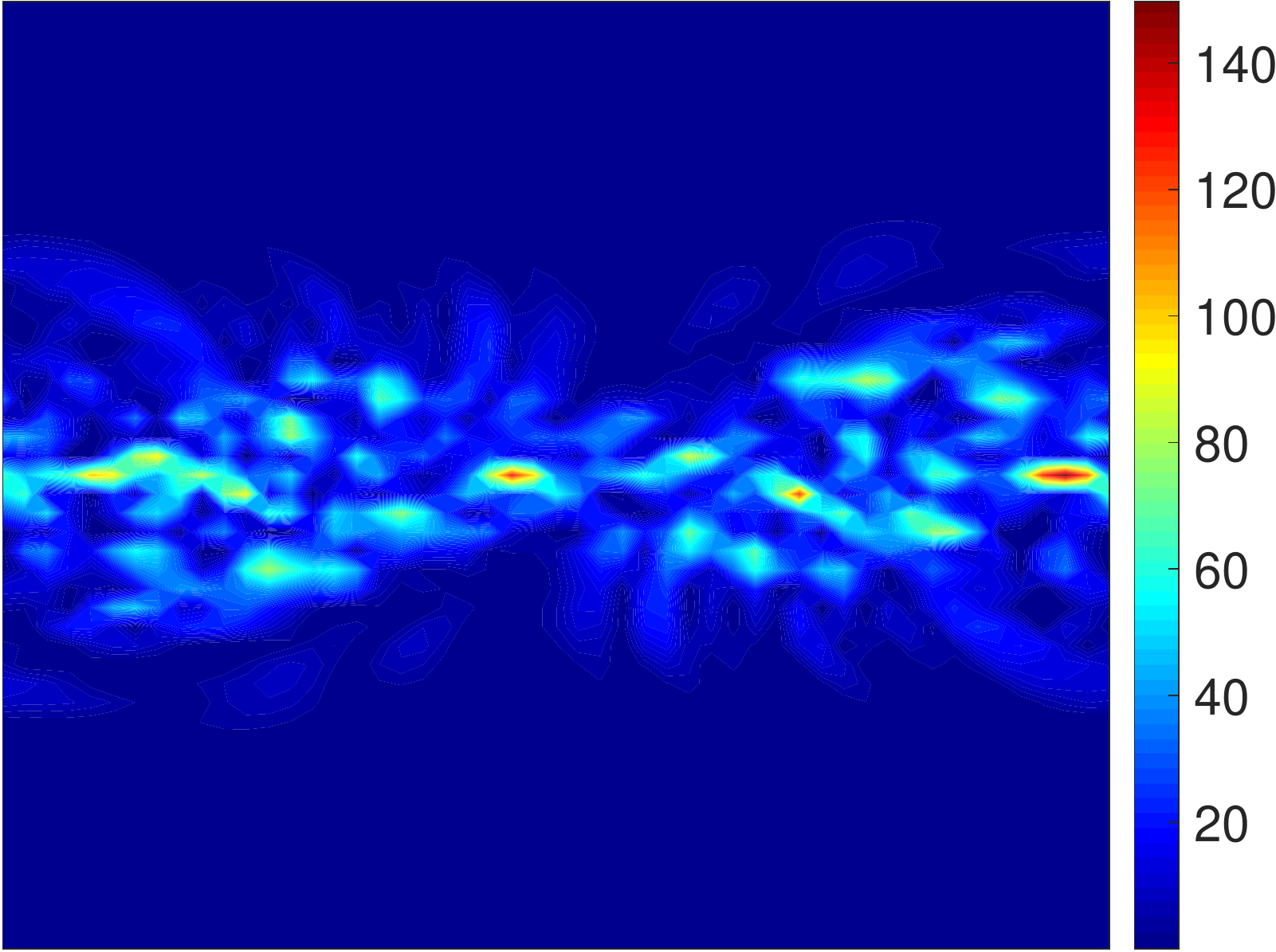}
\includegraphics[width=.24\textwidth, height=.12\textwidth,viewport=0 0 520 400, clip]{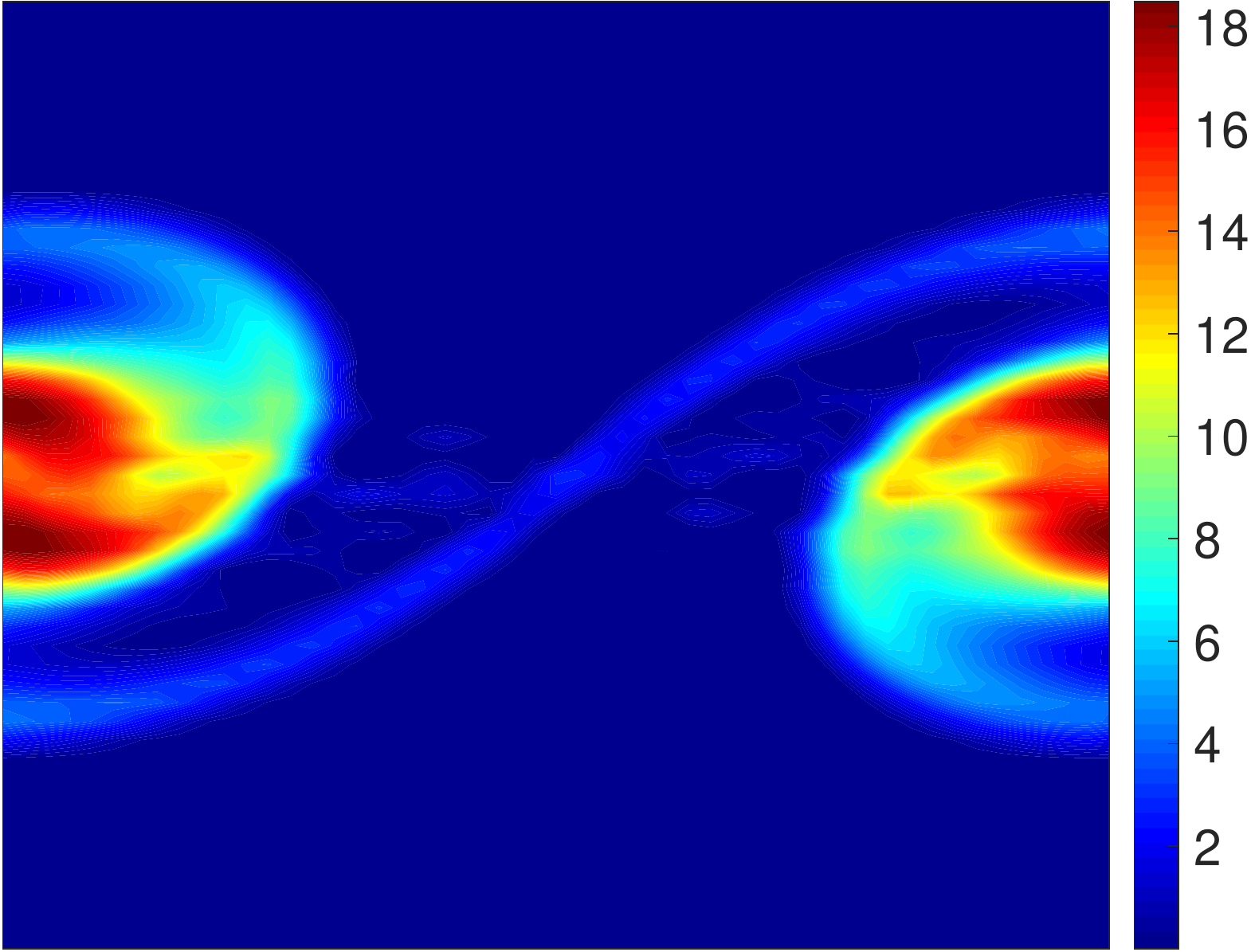} 
\includegraphics[width=.24\textwidth, height=.12\textwidth,viewport=0 0 520 400, clip]{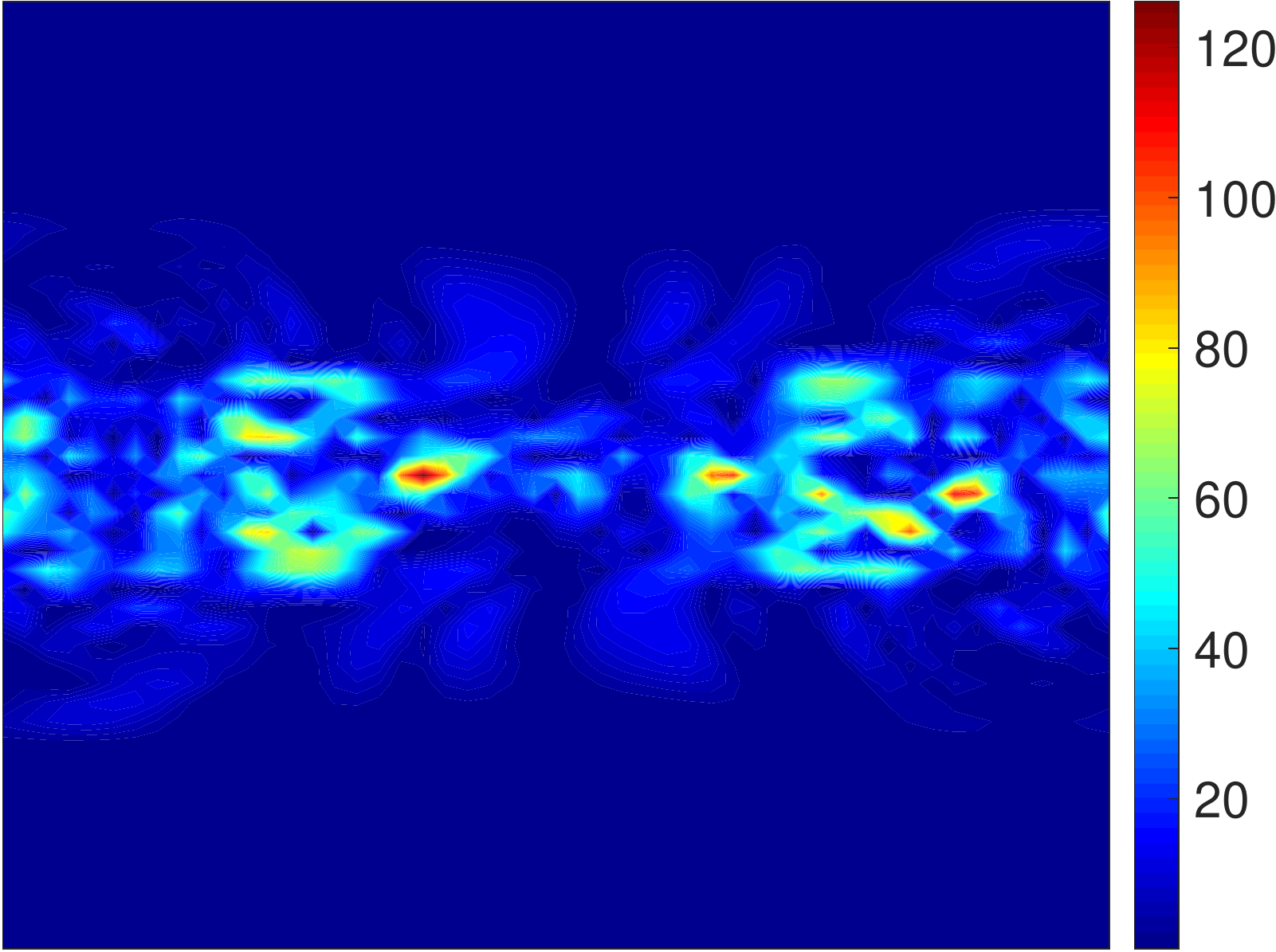}
\includegraphics[width=.24\textwidth, height=.12\textwidth,viewport=0 0 520 400, clip]{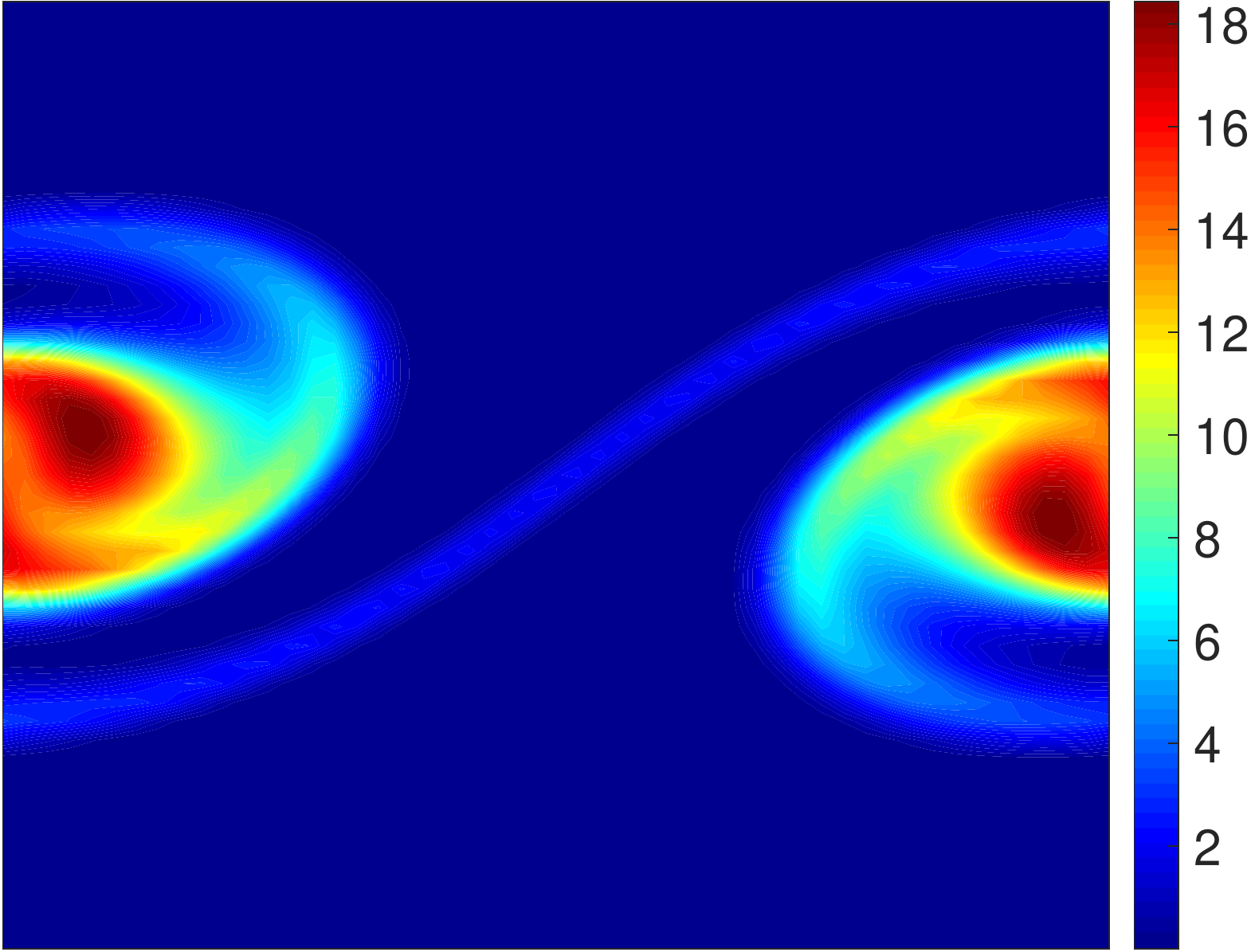}\\
\includegraphics[width=.24\textwidth, height=.12\textwidth,viewport=0 0 520 400, clip]{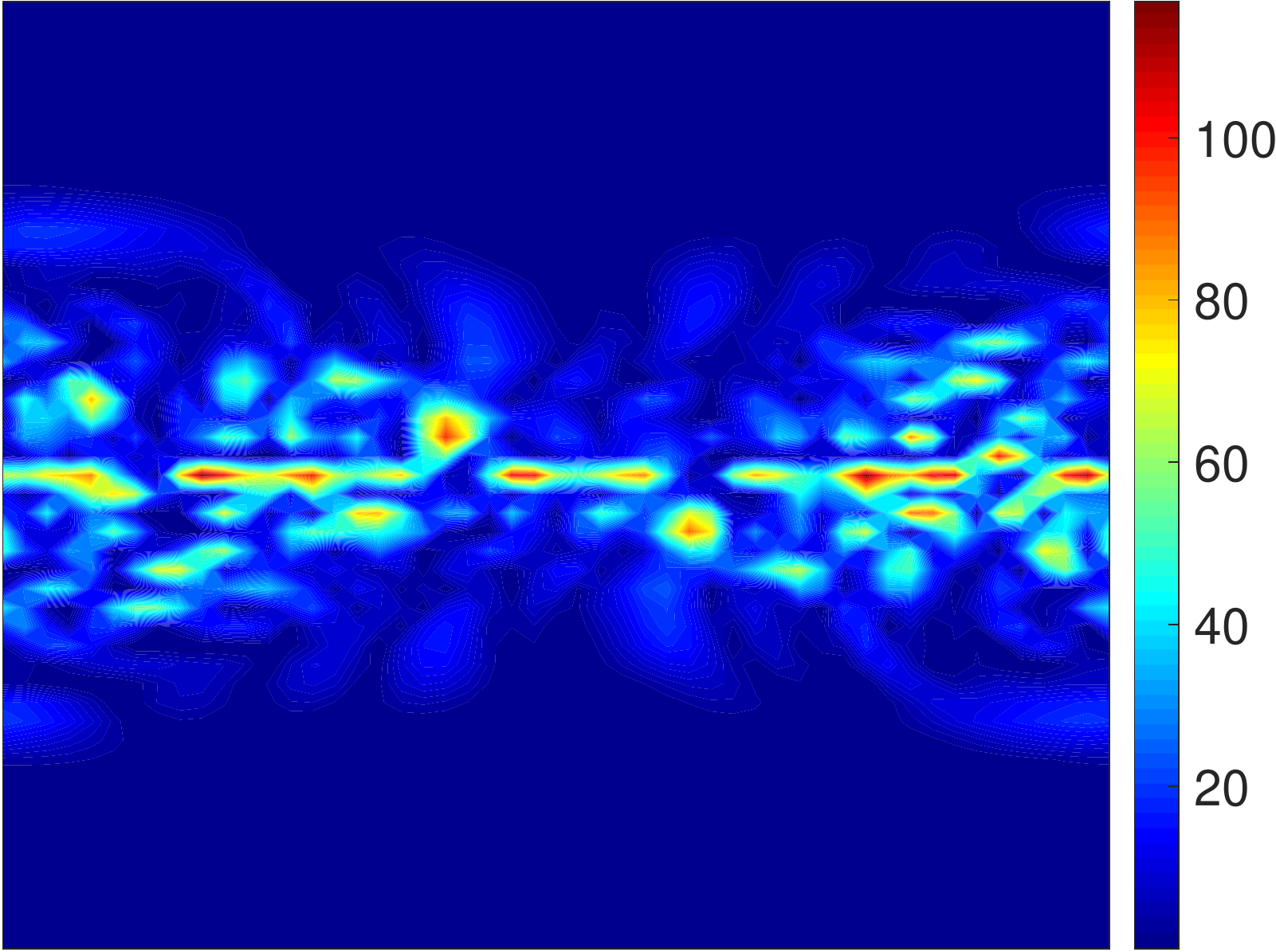}
\includegraphics[width=.24\textwidth, height=.12\textwidth,viewport=0 0 520 400, clip]{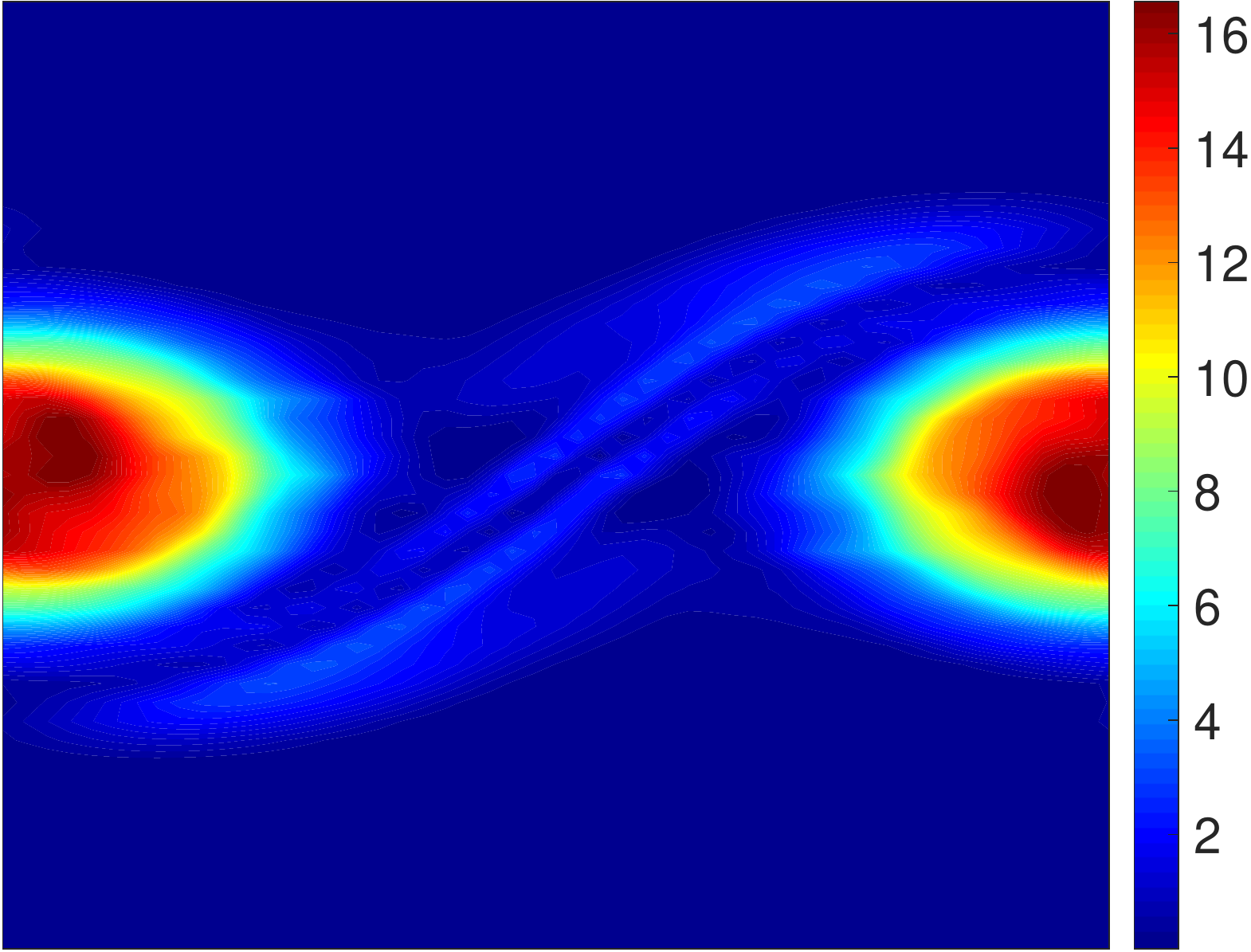} 
\includegraphics[width=.24\textwidth, height=.12\textwidth,viewport=0 0 520 400, clip]{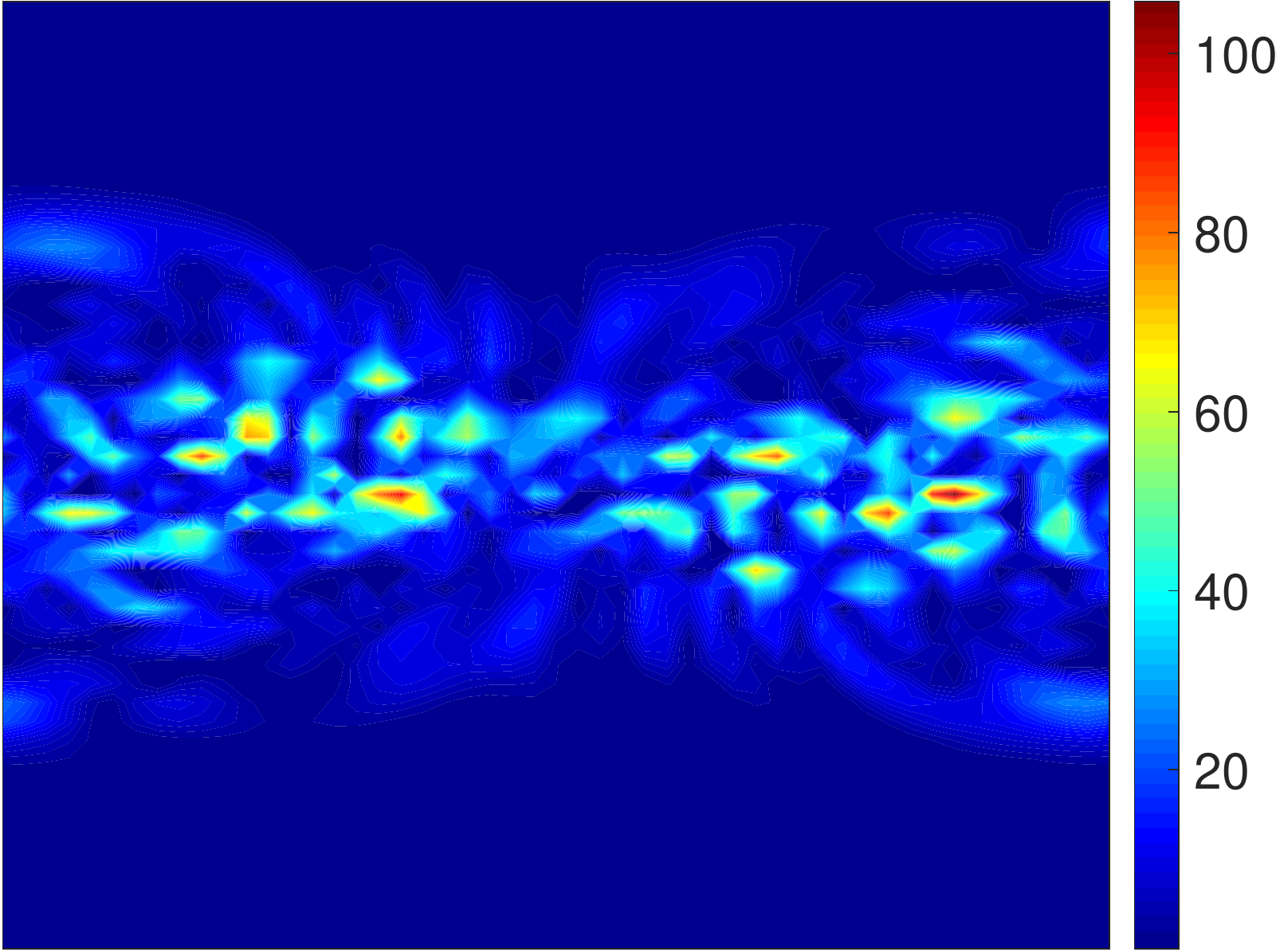}
\includegraphics[width=.24\textwidth, height=.12\textwidth,viewport=0 0 520 400, clip]{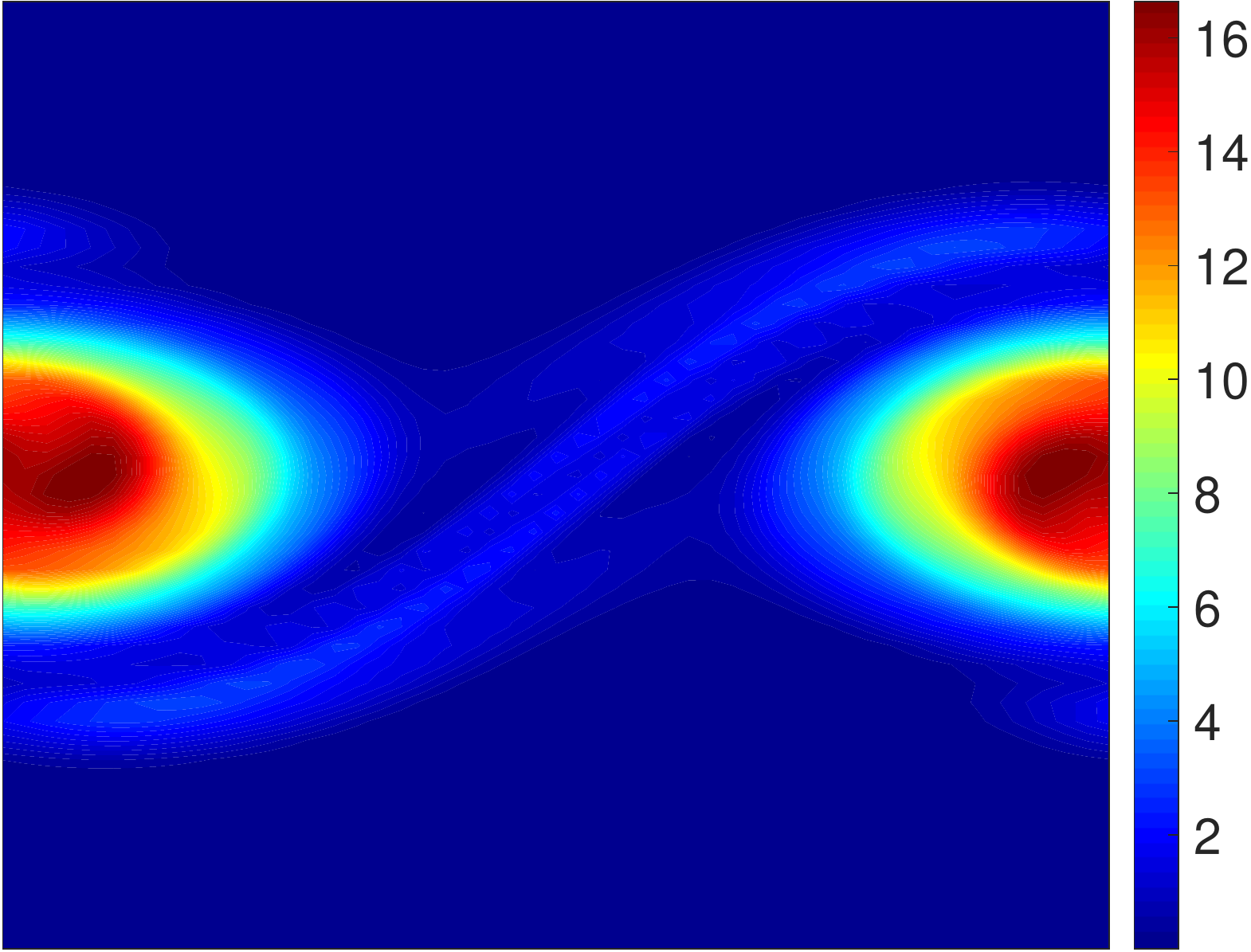}
\end{center}
\caption{\label{khemac}
Shown above are results (shown as vorticity contours) of ROM simulations that used EMAC-FEM as the FOM,  at t=1, 2, 3, 4, 5 (from top to bottom) for the $Re=100$ KH tests using 50 modes.}
\end{figure}

Results from using EMAC-FEM as the FOM are shown in Figures \ref{khemac} and \ref{emacEE}.  From Figure \ref{emacEE}, we observe that the FOM \LR{energy ($\frac12 \| u(t) \|^2$) and enstrophy ($\frac12 \| \nabla \times u(t) \|^2$)} match 
those of the DNS, at least up the precision of the plots.  Both the SKEW-ROM and CONV-ROM appear to become unstable, as their enstrophies grow dramatically after $t=2$ and 
stay high.  In turn, their energies decay much more quickly than the \LR{ FOM.}  The EMAC-ROM has the same energy and enstrophy as the other ROMs before $t=2$, but when the other ROMs become (essentially) unstable, the EMAC-ROM remains stable and gives reasonably accurate answers.  Figure \ref{khemac} shows vorticity contours of solutions are various time steps for the ROMs and the FOM, and this plot makes is clear how the EMAC-ROM gives a good approximation of the FOM (EMAC-FEM), while SKEW-ROM and CONV-ROM give terrible approximations.  

\subsubsection{Comparison of ROMs using SKEW-FEM as FOM}
    \label{sec:kh-skew-fem}

\begin{figure}[!h]
\begin{center}
\includegraphics[width=.45\textwidth, height=.32\textwidth,viewport=0 0 520 400, clip]{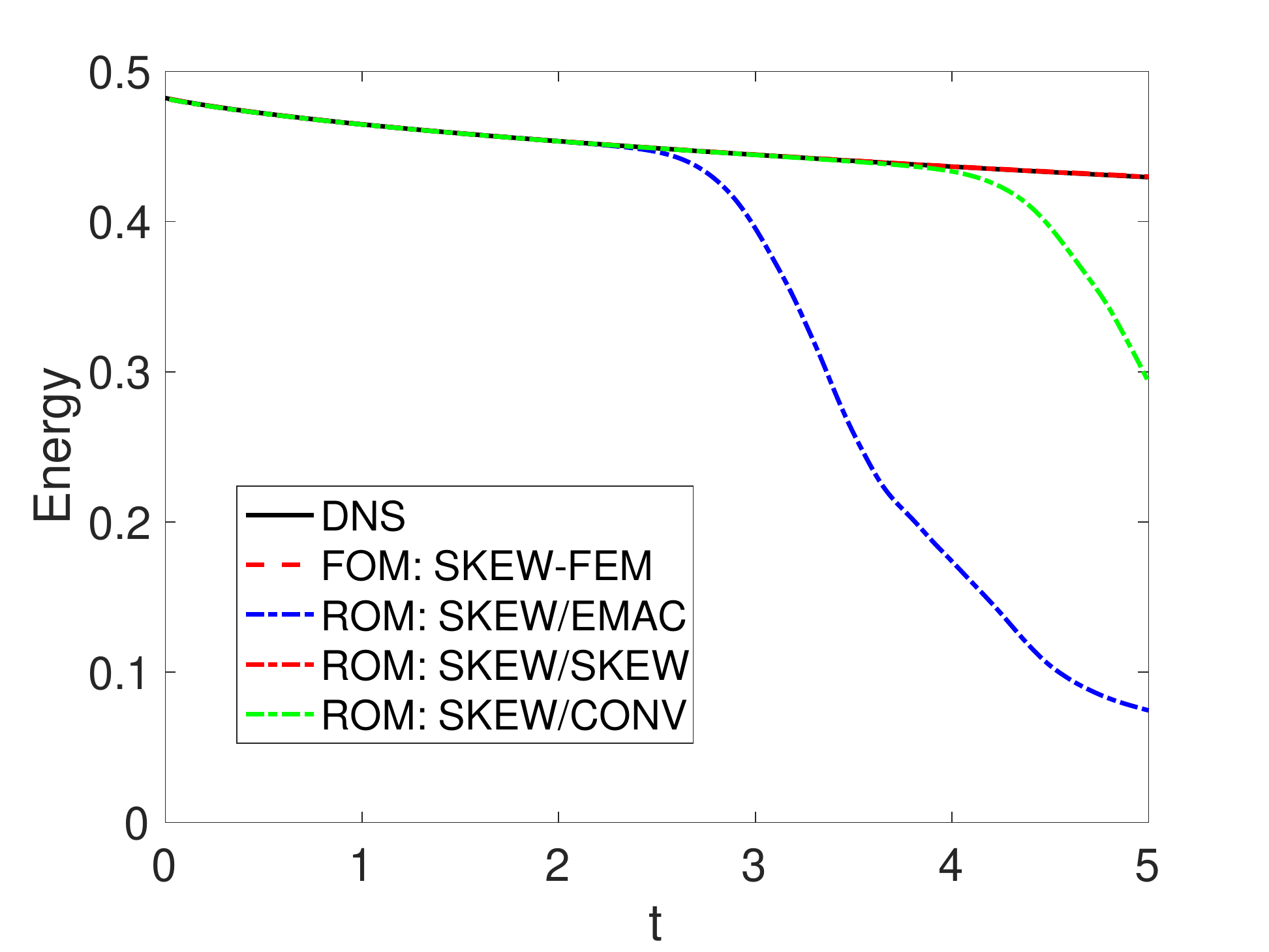}
\includegraphics[width=.45\textwidth, height=.32\textwidth,viewport=0 0 520 400, clip]{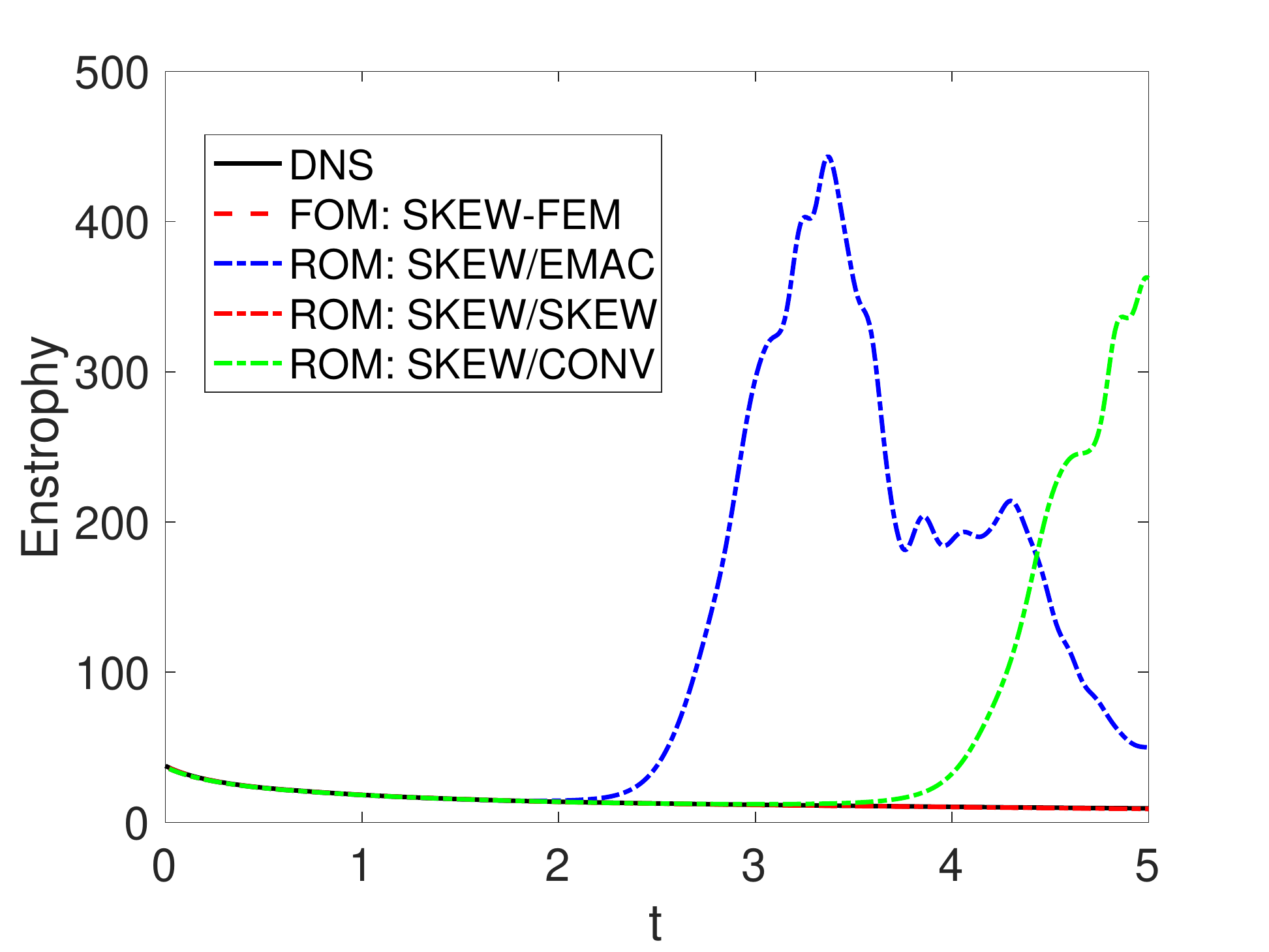} \\
\end{center}
\caption{\label{skewEE}
Shown above are the energy and enstrophy plots versus time for the $Re=100$ KH tests for varying ROMs (that are built from SKEW-FEM), DNS, the FOM (SKEW-FEM), and a DNS from  \cite{SJLLLS18}.
}
\end{figure}

\begin{figure}[!h]
\hspace{.2in} SKEW-ROM \hspace{.6in} EMAC-ROM \hspace{.7in} CONV-ROM \hspace{.4in} FOM (SKEW-FEM) 
\begin{center}
\includegraphics[width=.24\textwidth, height=.12\textwidth,viewport=0 0 520 400, clip]{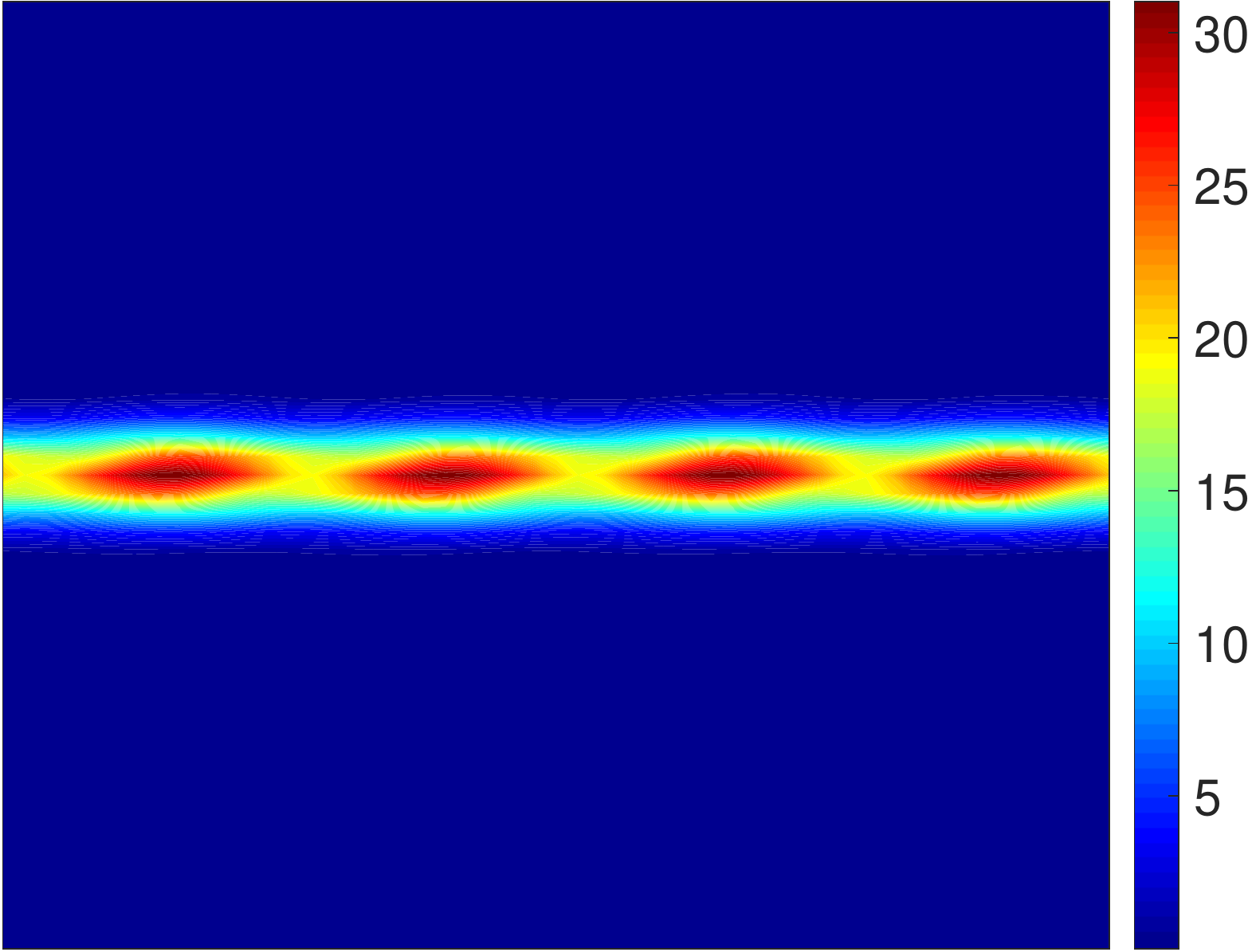}
\includegraphics[width=.24\textwidth, height=.12\textwidth,viewport=0 0 520 400, clip]{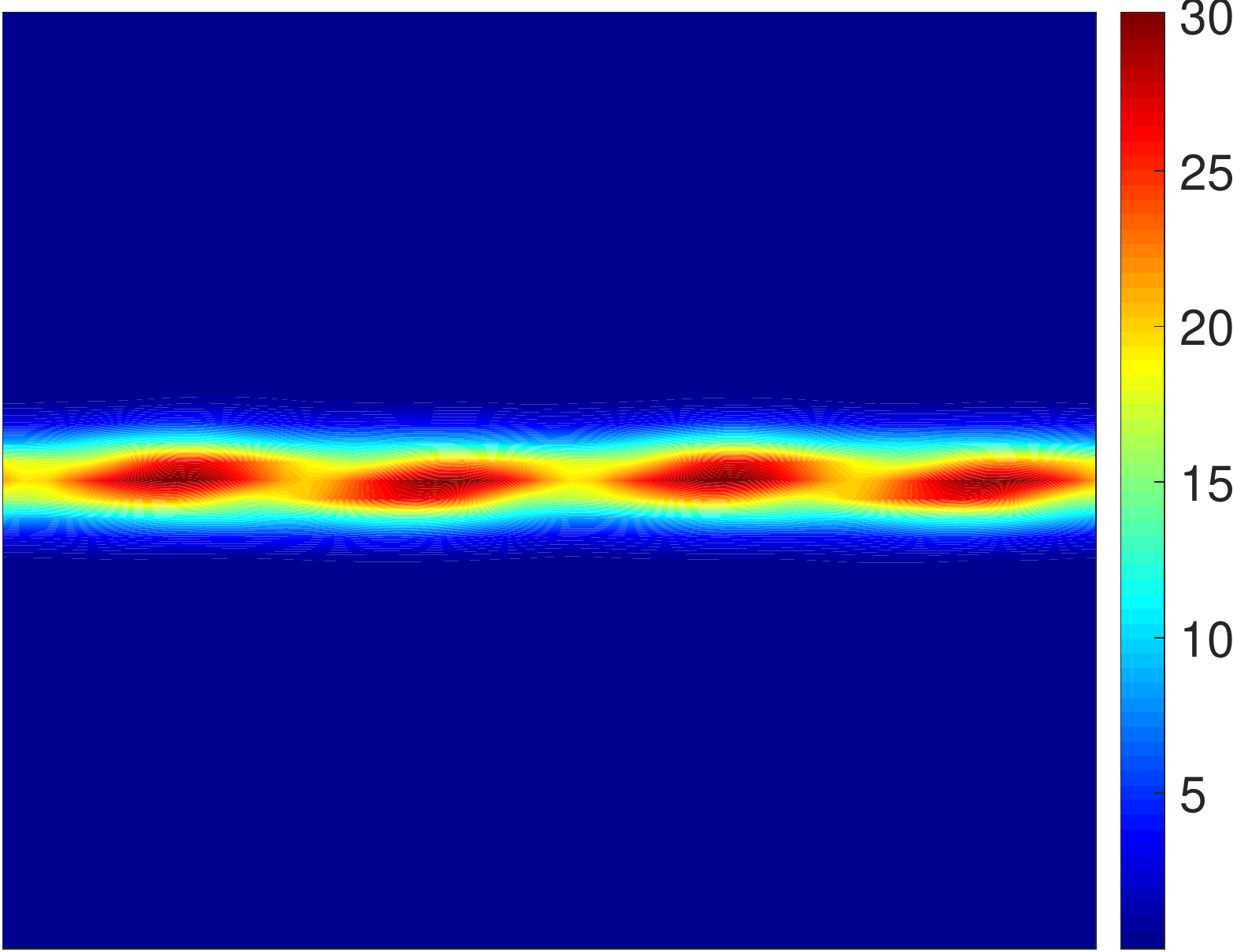} 
\includegraphics[width=.24\textwidth, height=.12\textwidth,viewport=0 0 520 400, clip]{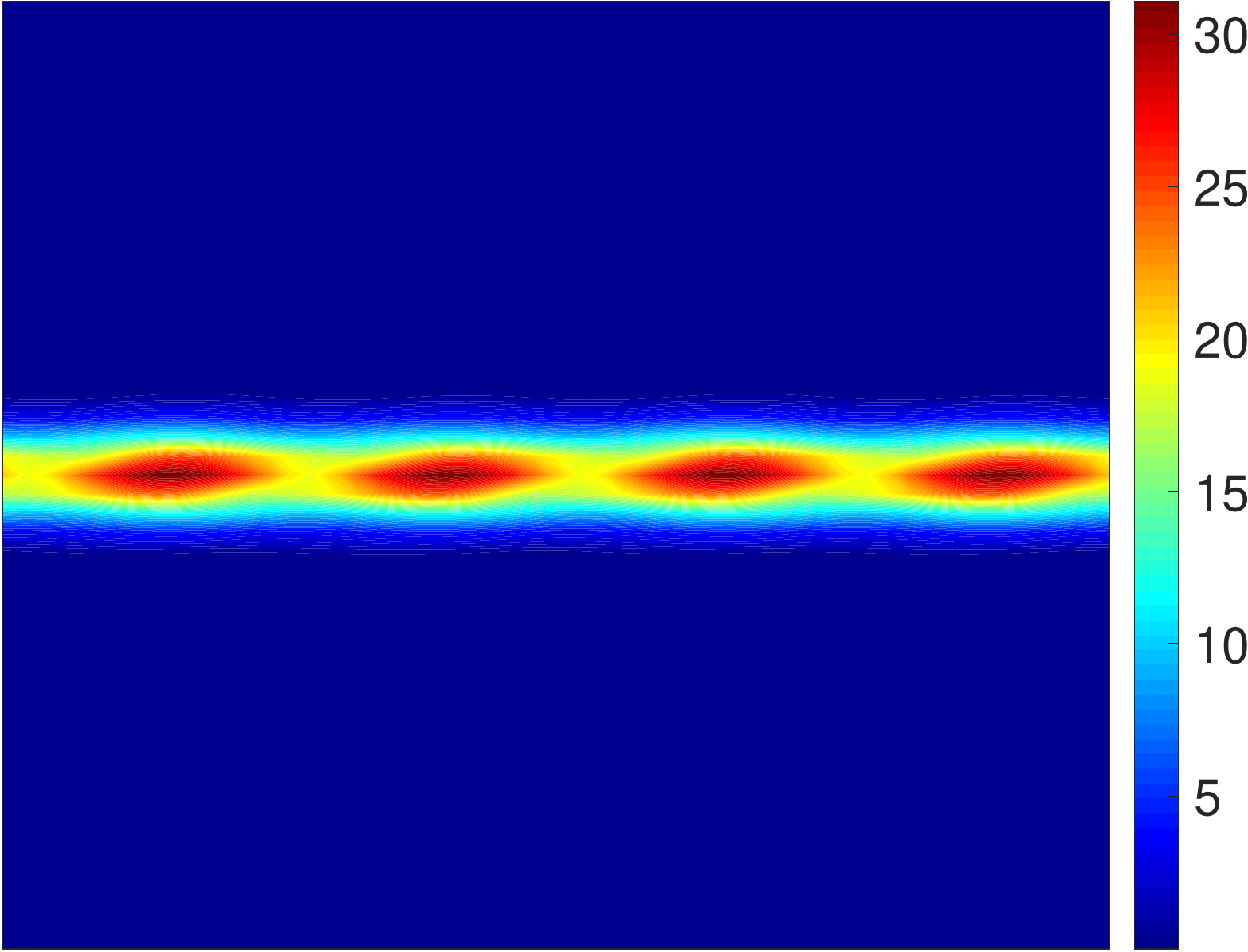}
\includegraphics[width=.24\textwidth, height=.12\textwidth,viewport=0 0 520 400, clip]{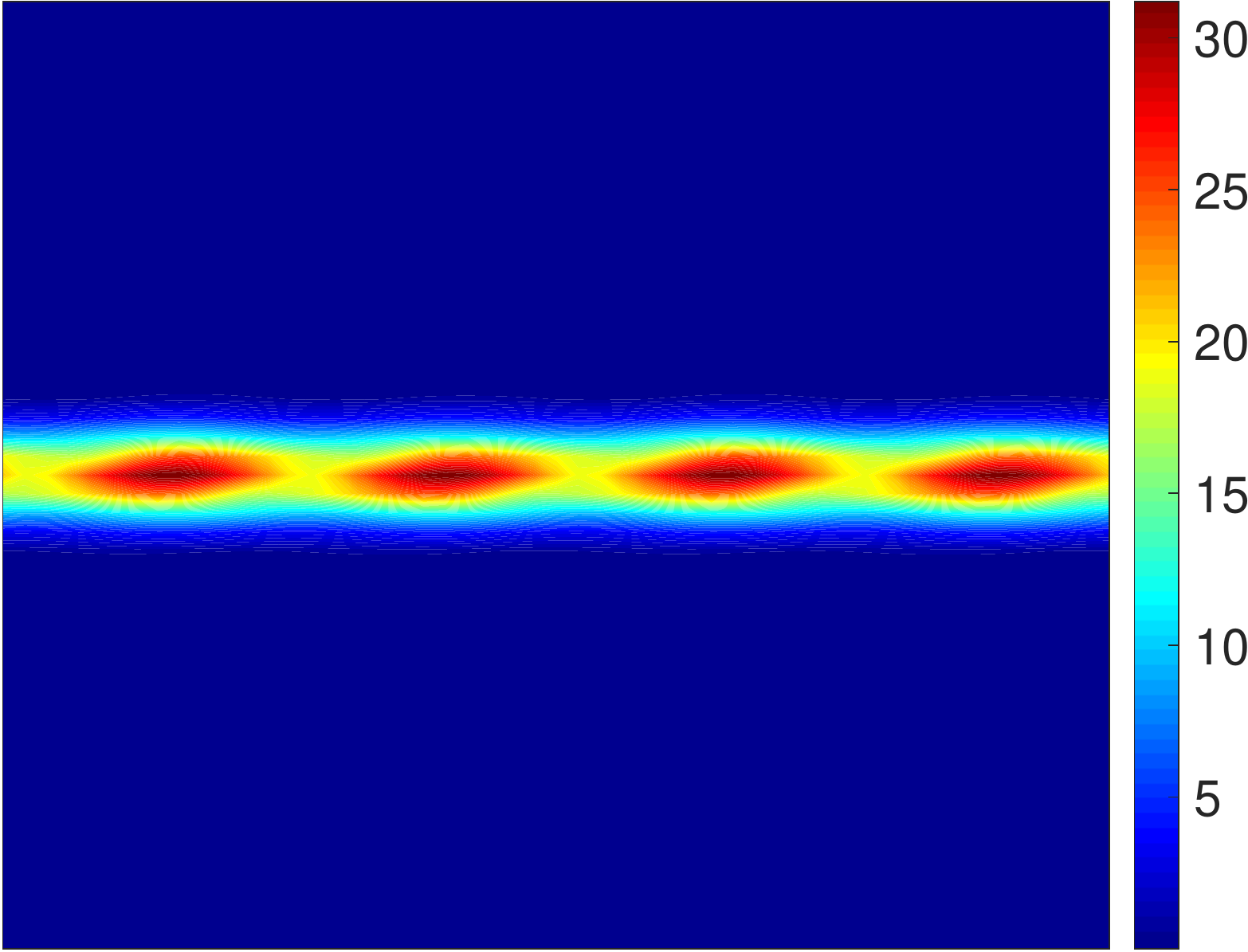}\\
\includegraphics[width=.24\textwidth, height=.12\textwidth,viewport=0 0 520 400, clip]{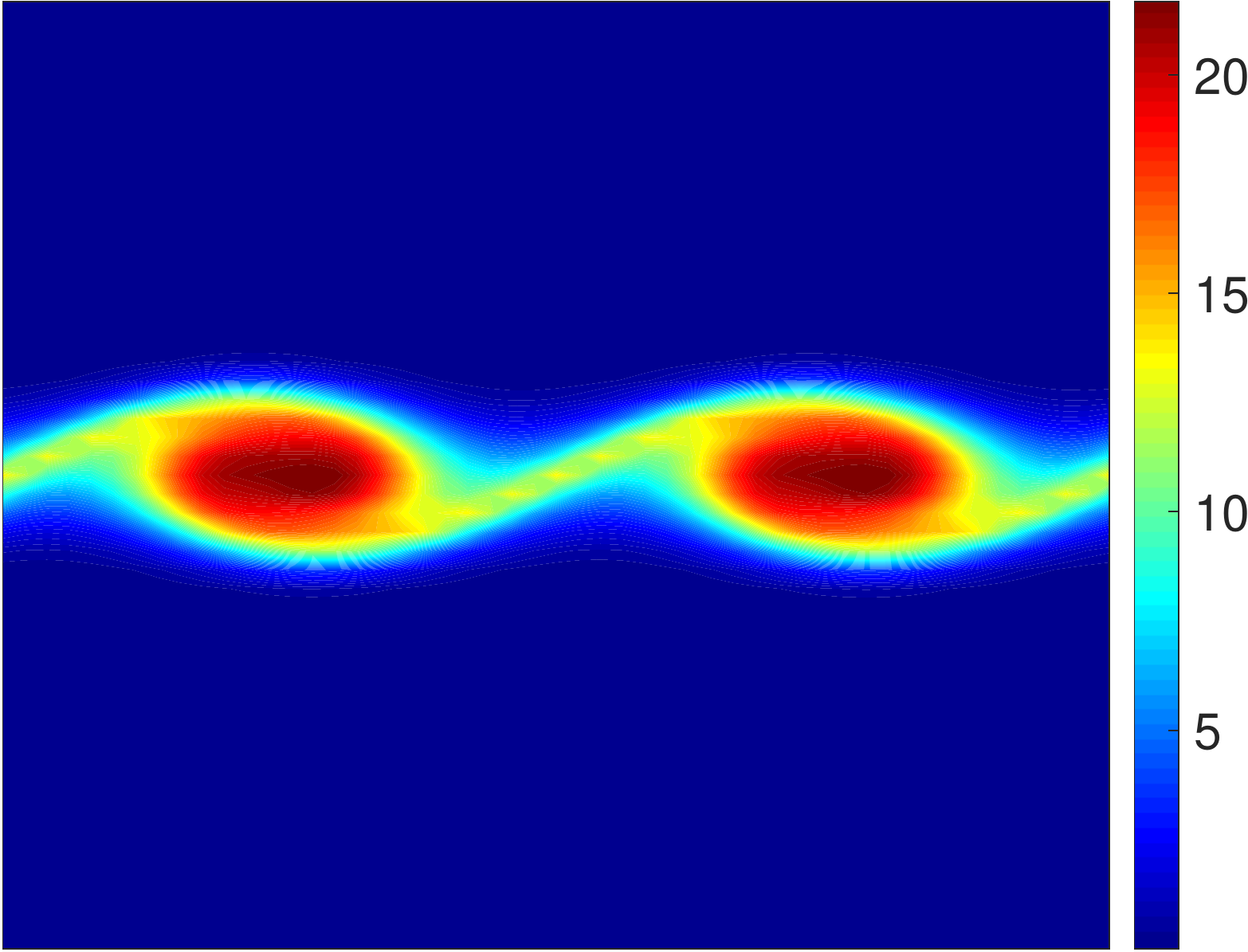}
\includegraphics[width=.24\textwidth, height=.12\textwidth,viewport=0 0 520 400, clip]{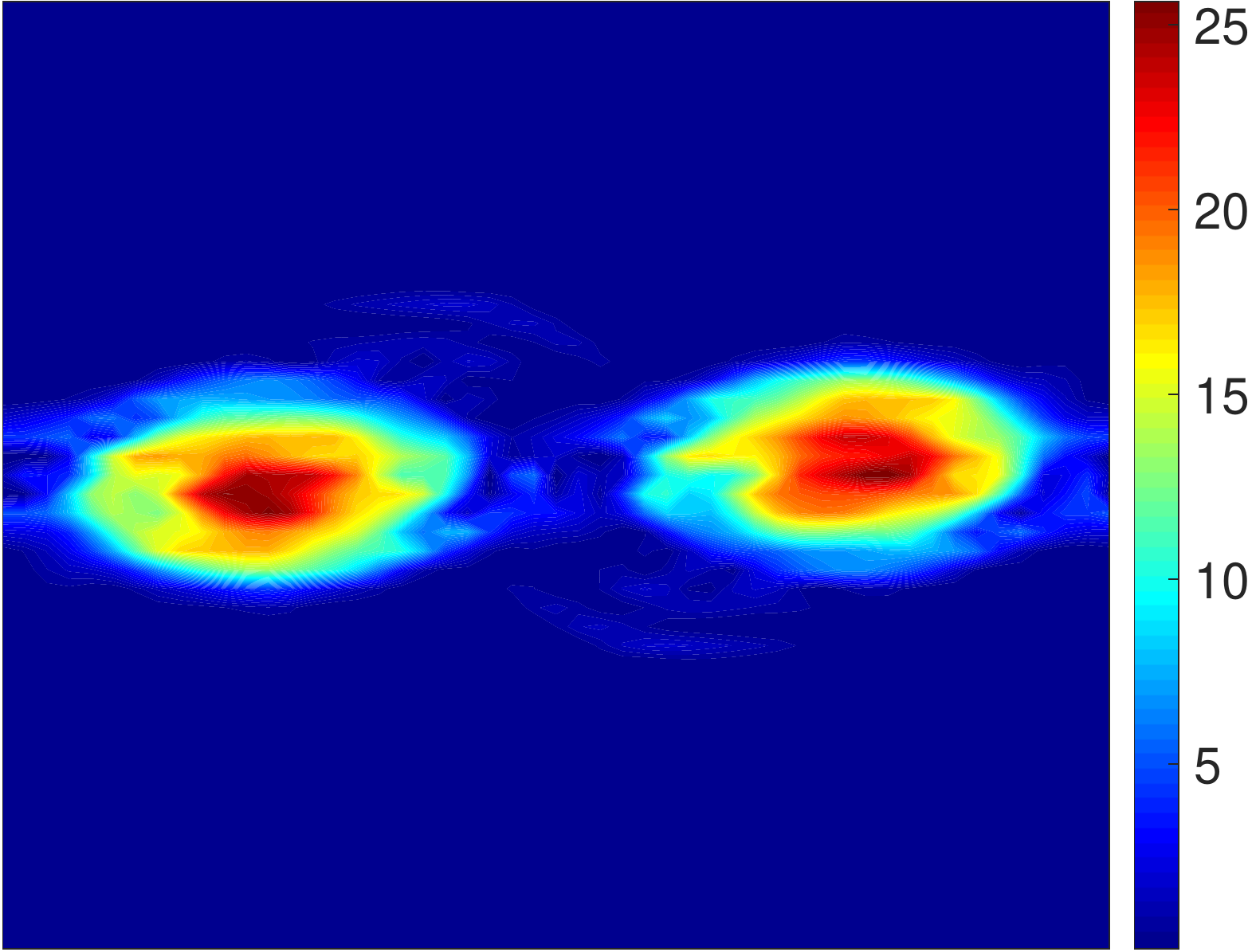} 
\includegraphics[width=.24\textwidth, height=.12\textwidth,viewport=0 0 520 400, clip]{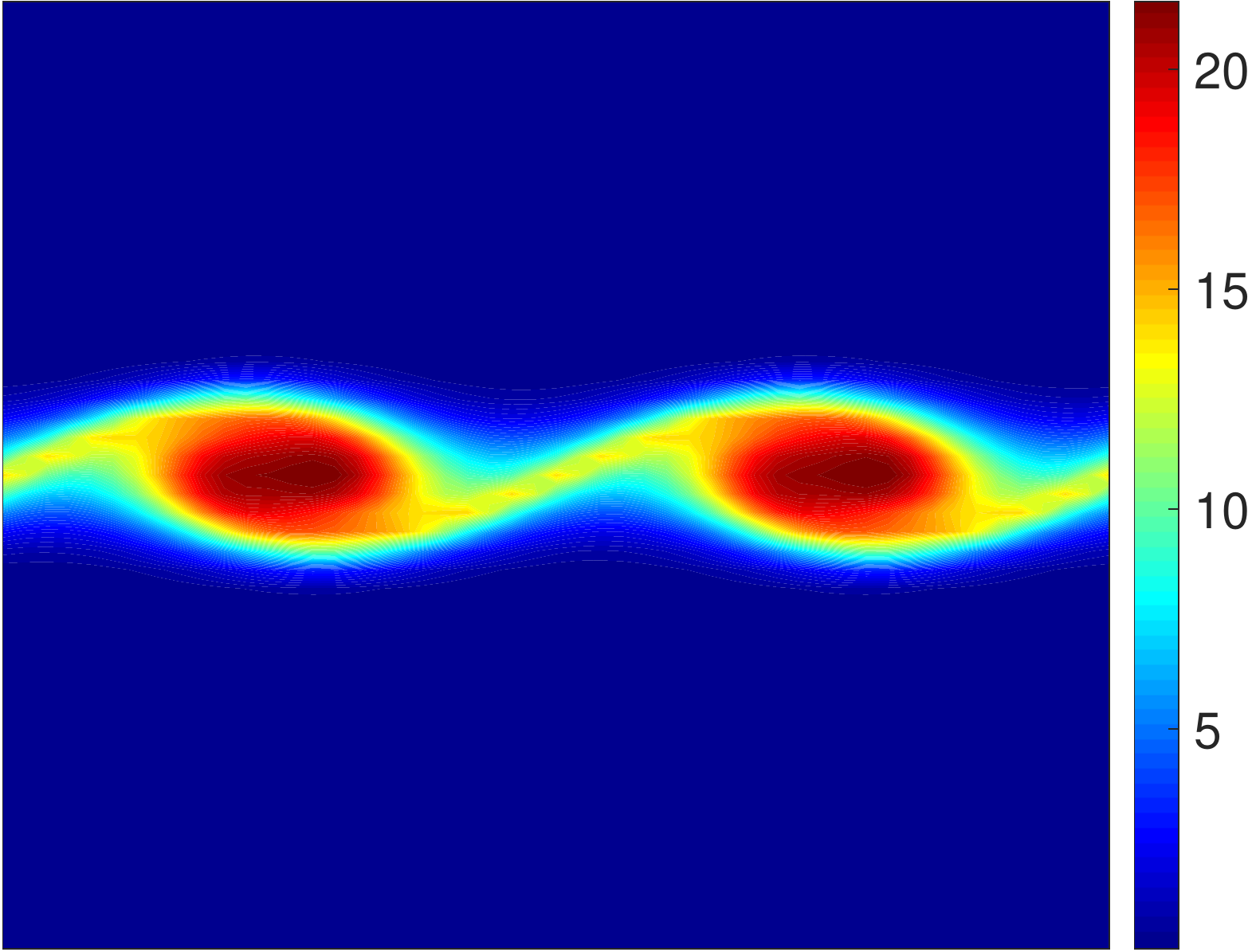}
\includegraphics[width=.24\textwidth, height=.12\textwidth,viewport=0 0 520 400, clip]{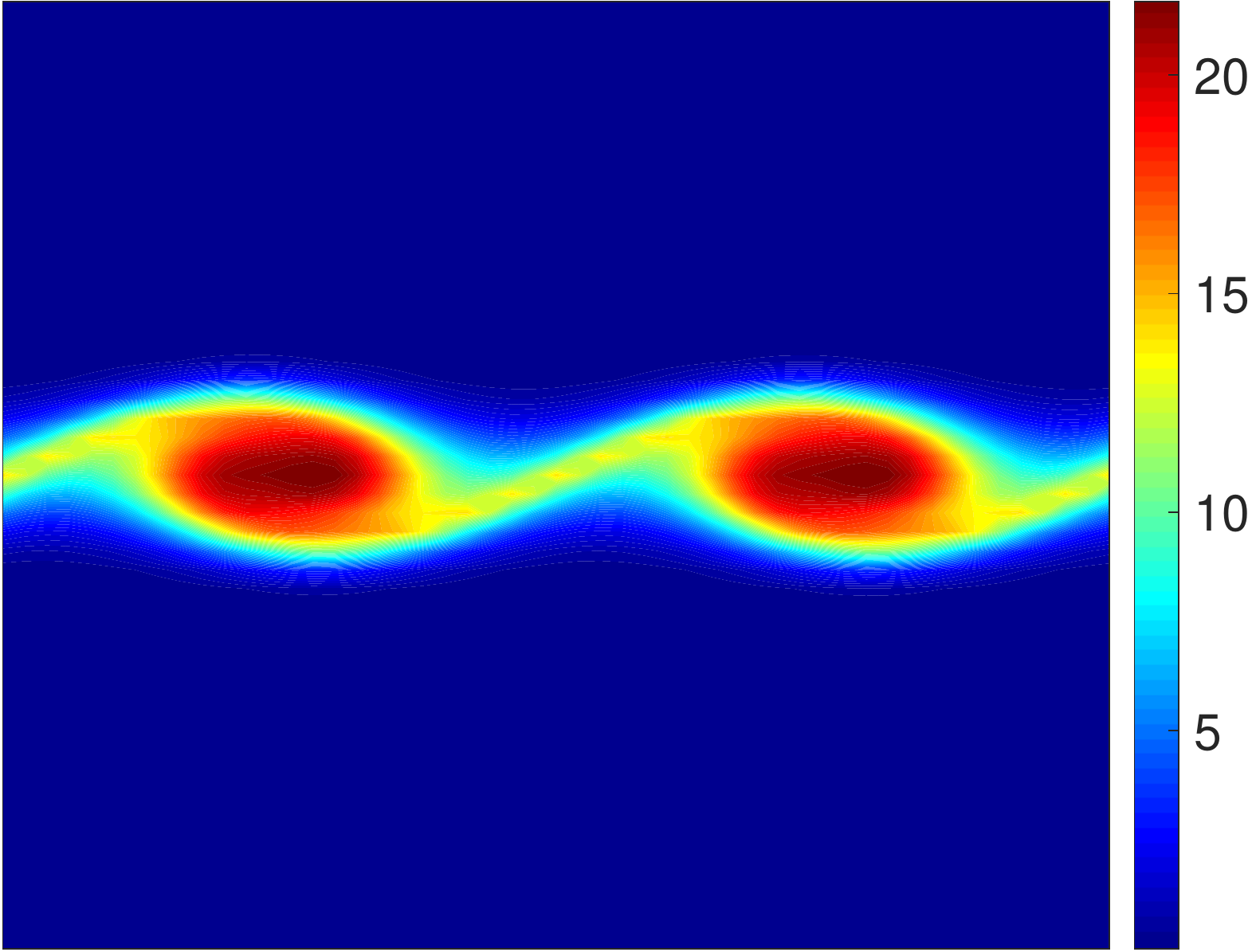}\\
\includegraphics[width=.24\textwidth, height=.12\textwidth,viewport=0 0 520 400, clip]{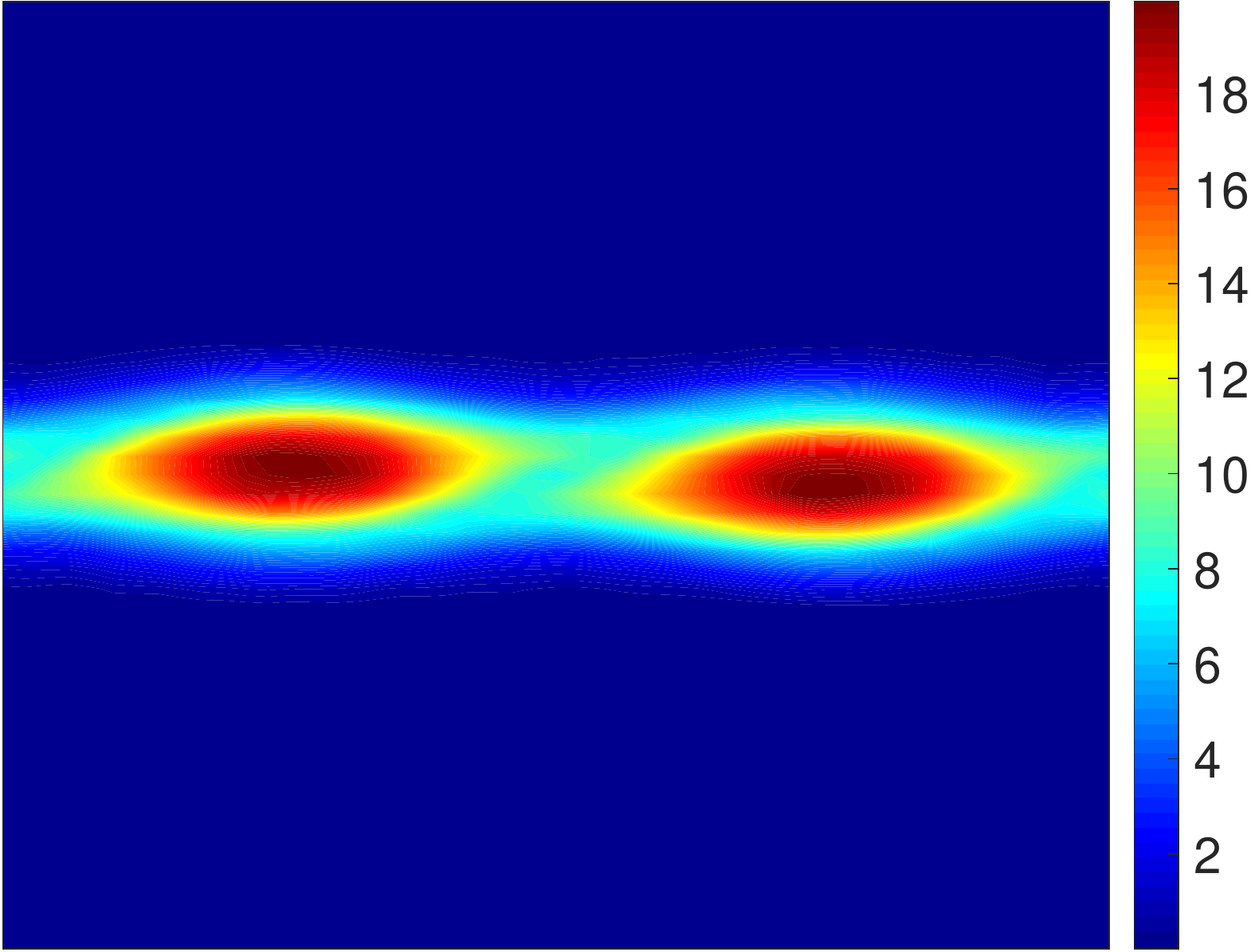}
\includegraphics[width=.24\textwidth, height=.12\textwidth,viewport=0 0 520 400, clip]{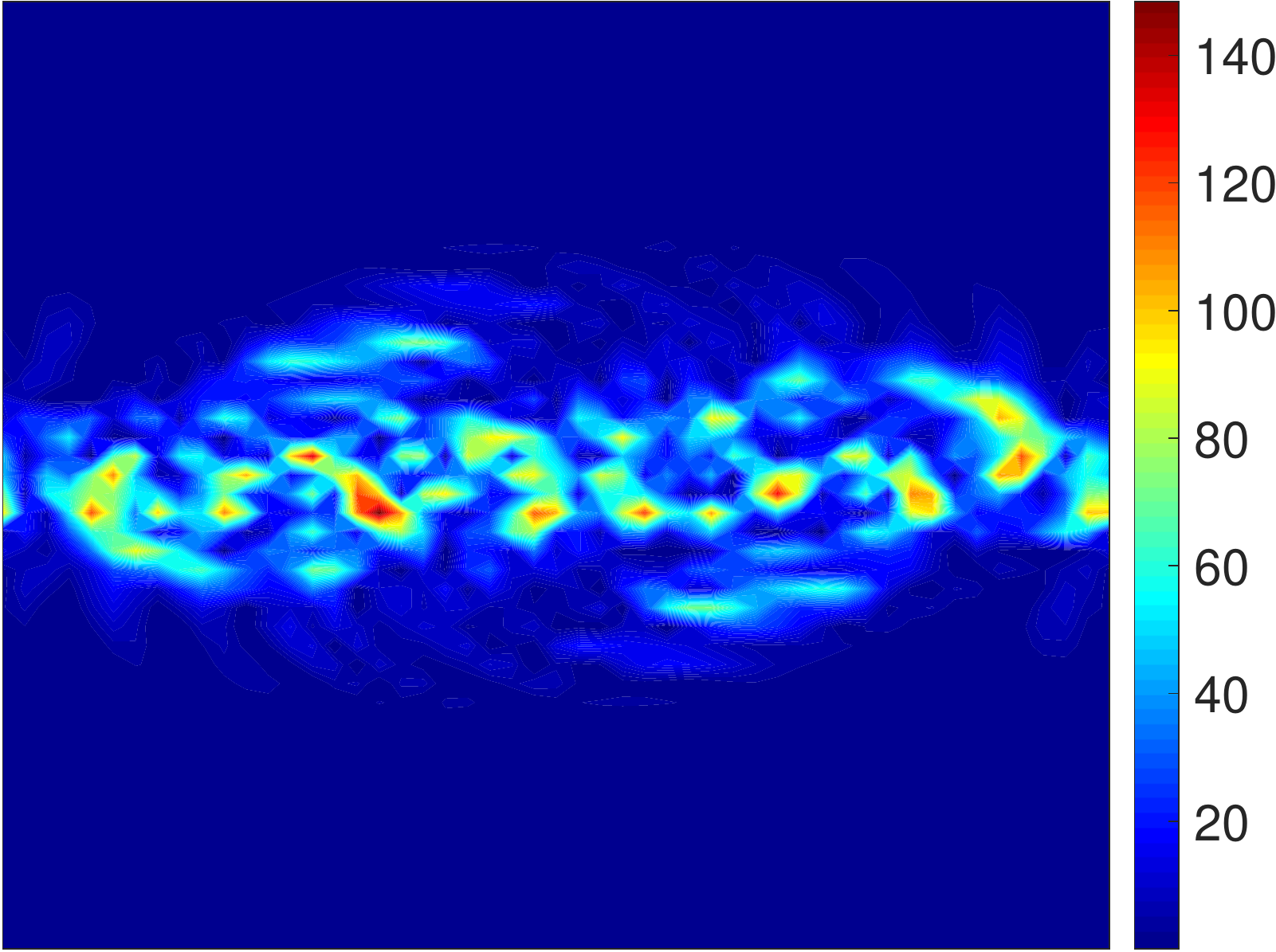} 
\includegraphics[width=.24\textwidth, height=.12\textwidth,viewport=0 0 520 400, clip]{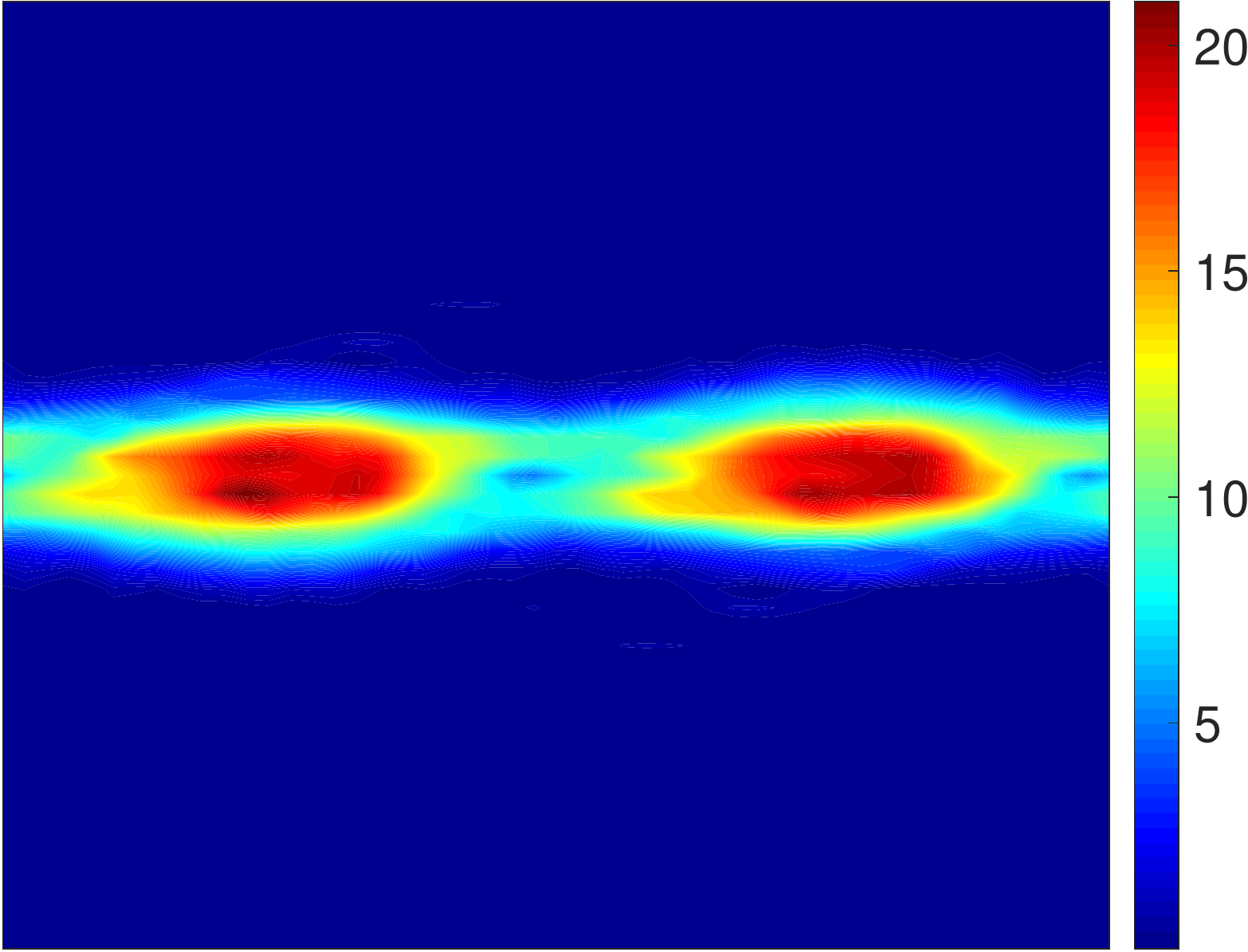}
\includegraphics[width=.24\textwidth, height=.12\textwidth,viewport=0 0 520 400, clip]{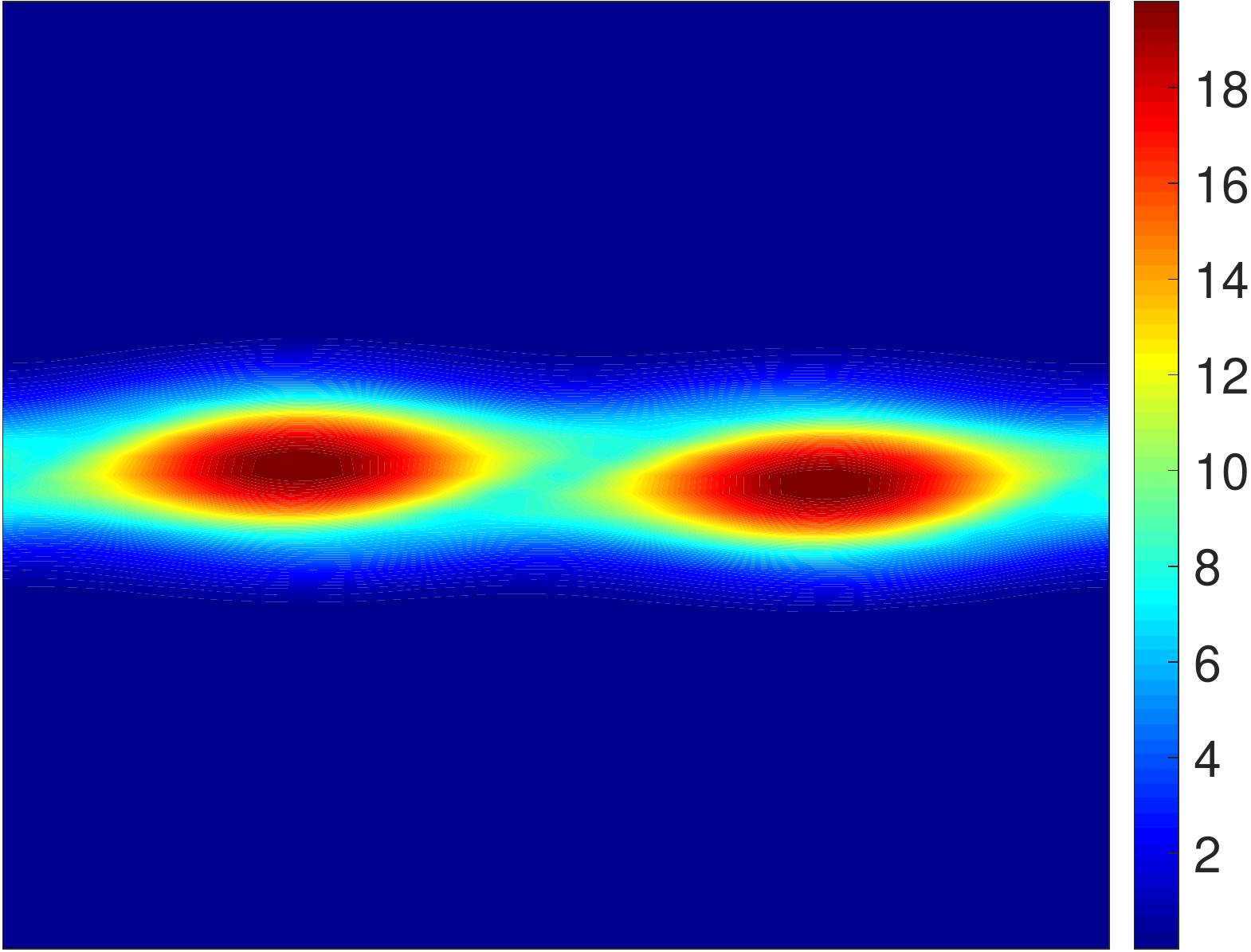}\\
\includegraphics[width=.24\textwidth, height=.12\textwidth,viewport=0 0 520 400, clip]{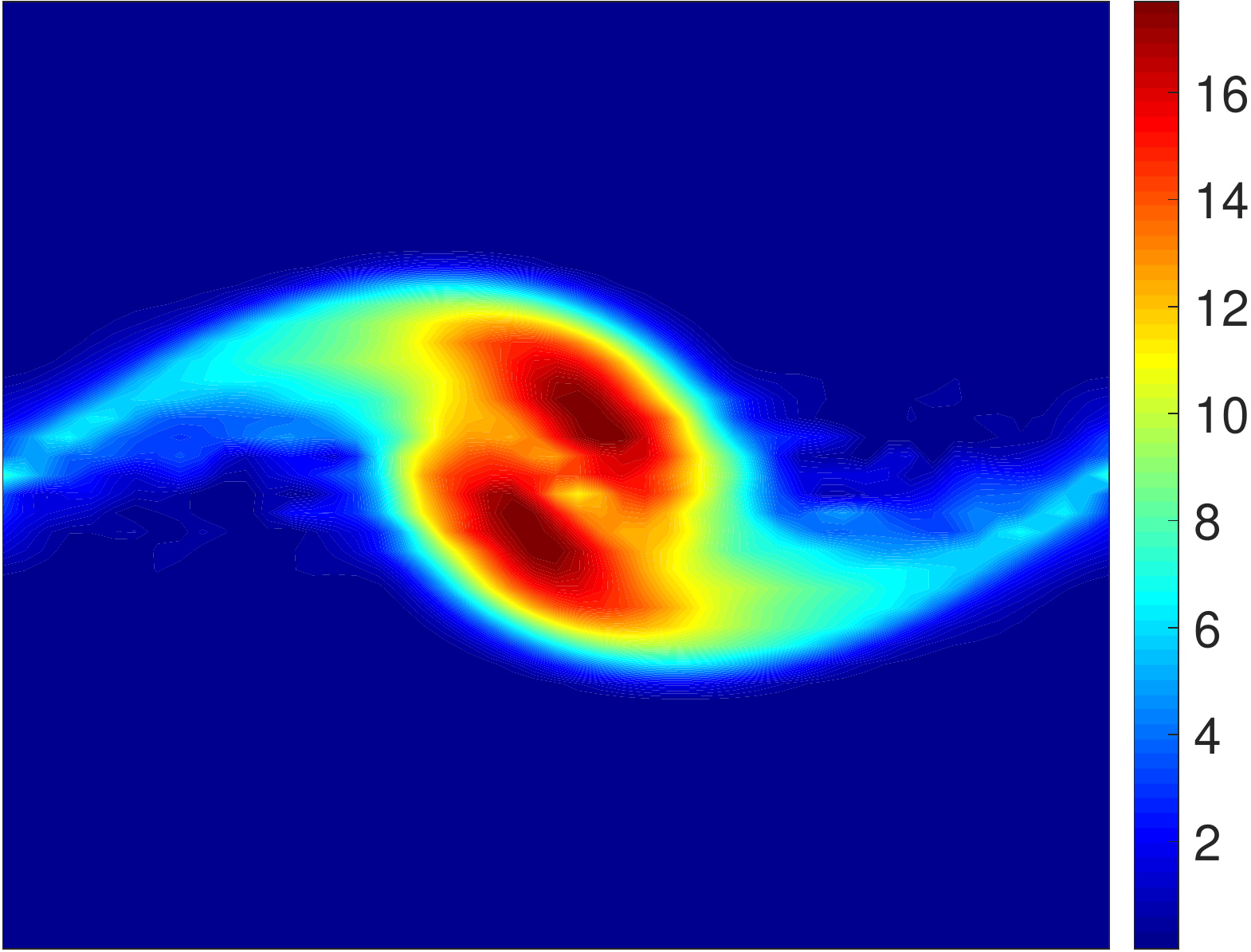}
\includegraphics[width=.24\textwidth, height=.12\textwidth,viewport=0 0 520 400, clip]{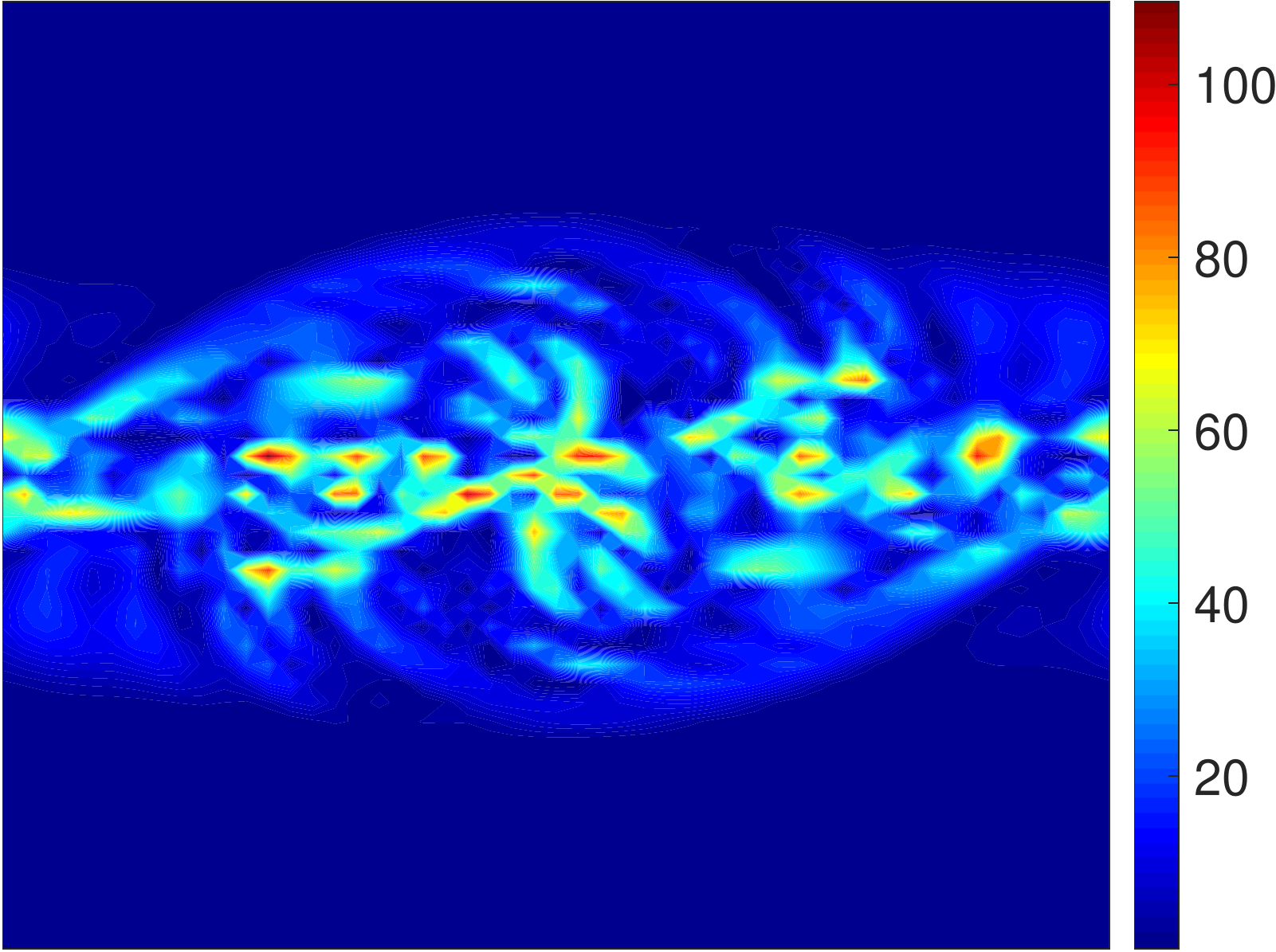} 
\includegraphics[width=.24\textwidth, height=.12\textwidth,viewport=0 0 520 400, clip]{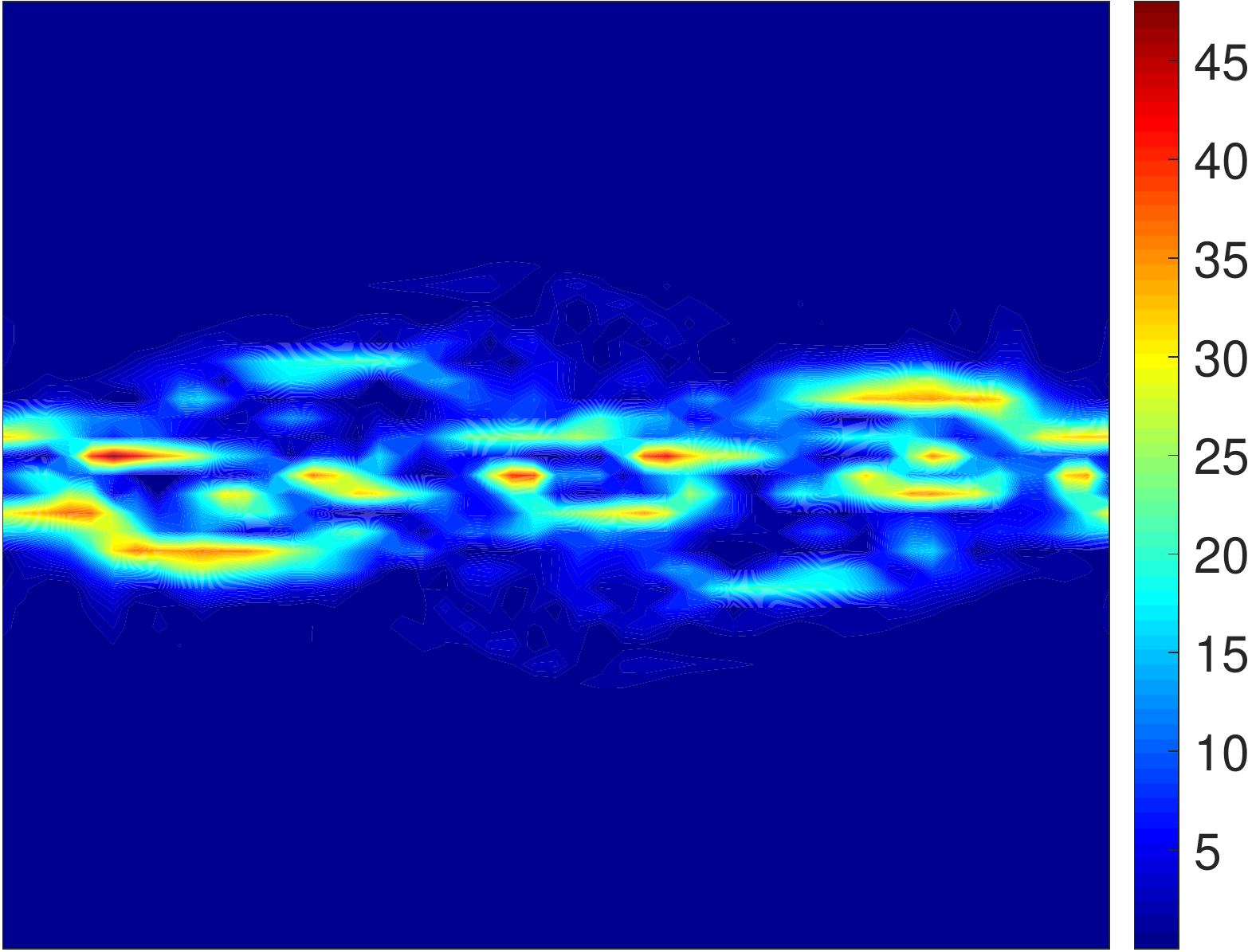}
\includegraphics[width=.24\textwidth, height=.12\textwidth,viewport=0 0 520 400, clip]{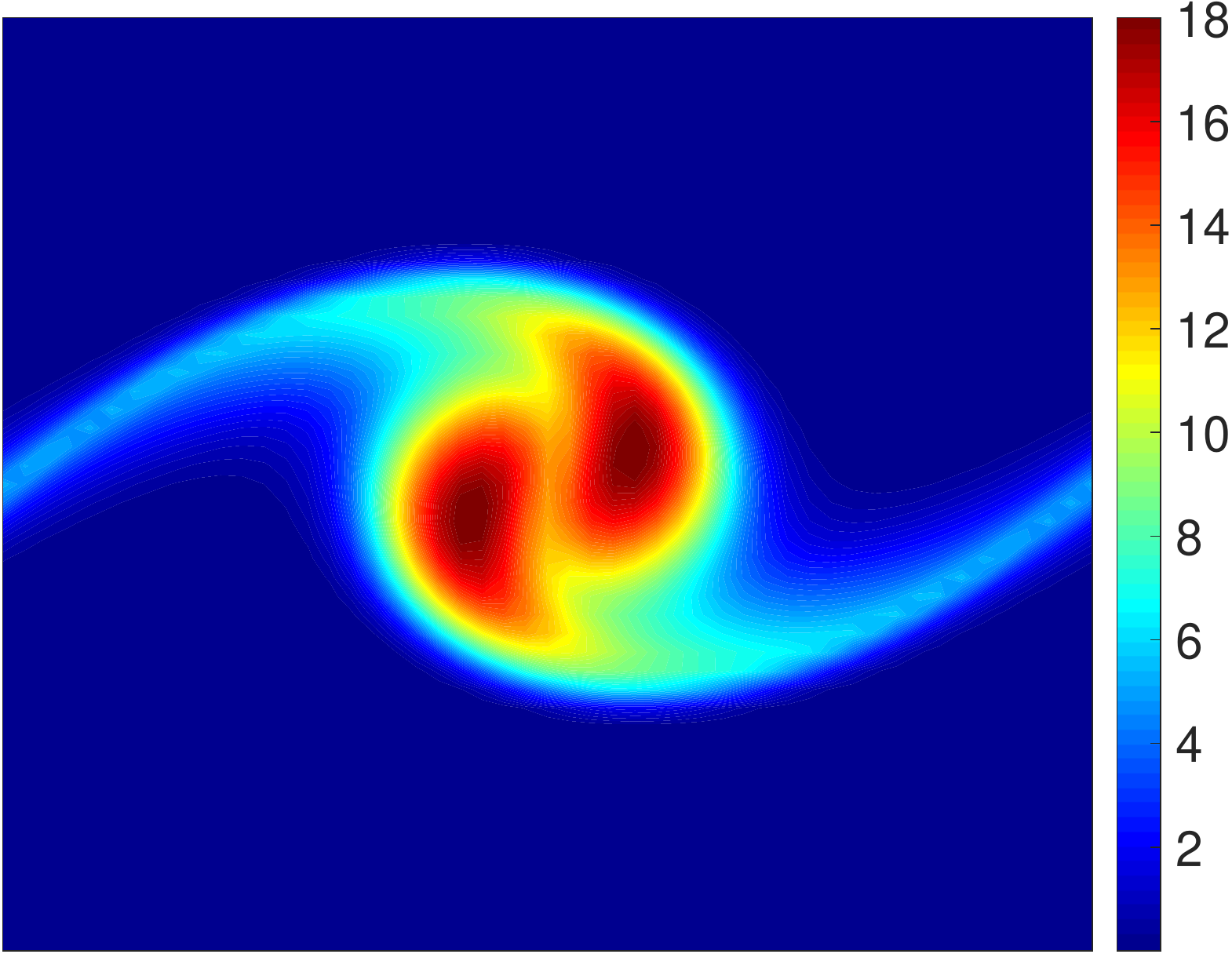}\\
\includegraphics[width=.24\textwidth, height=.12\textwidth,viewport=0 0 520 400, clip]{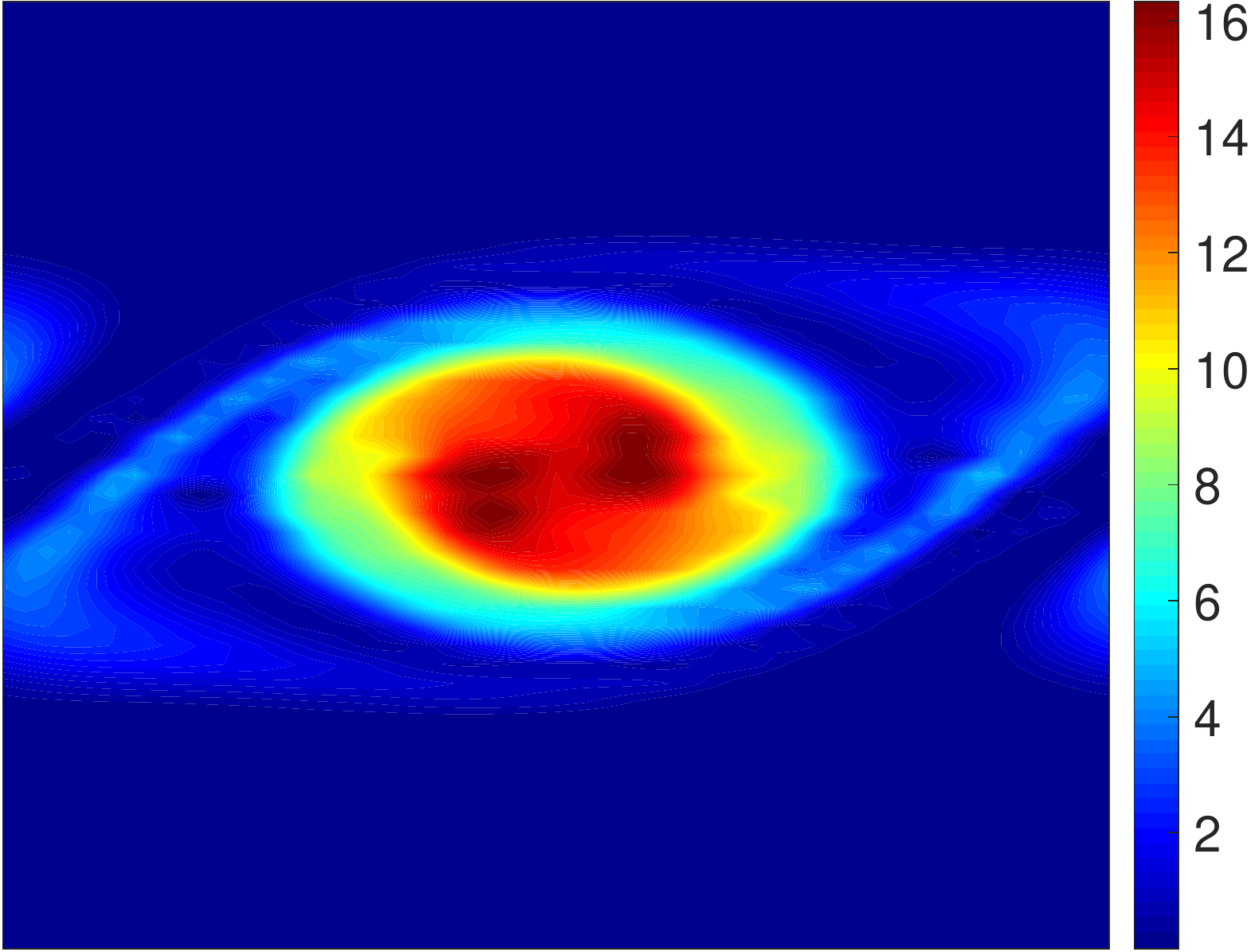}
\includegraphics[width=.24\textwidth, height=.12\textwidth,viewport=0 0 520 400, clip]{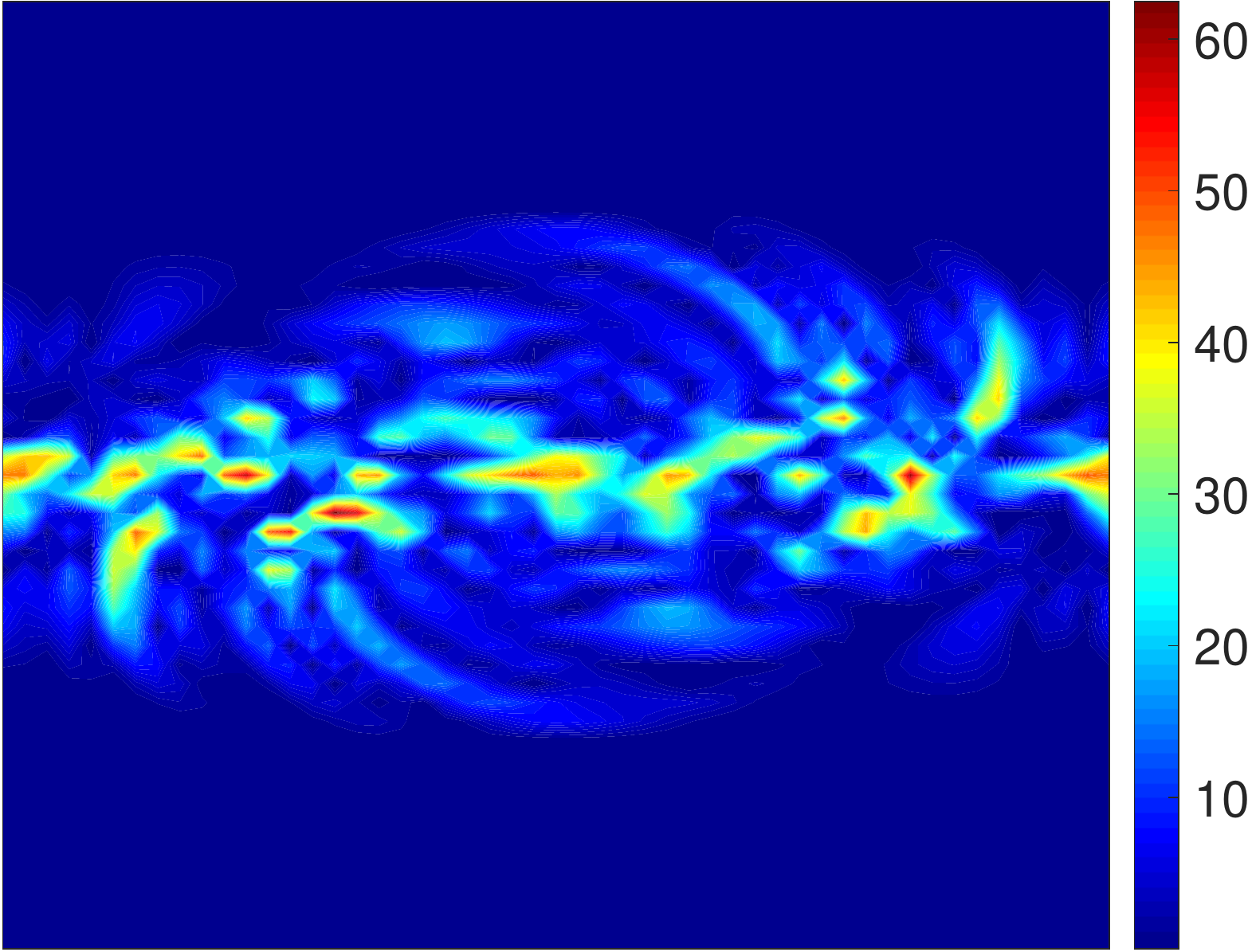} 
\includegraphics[width=.24\textwidth, height=.12\textwidth,viewport=0 0 520 400, clip]{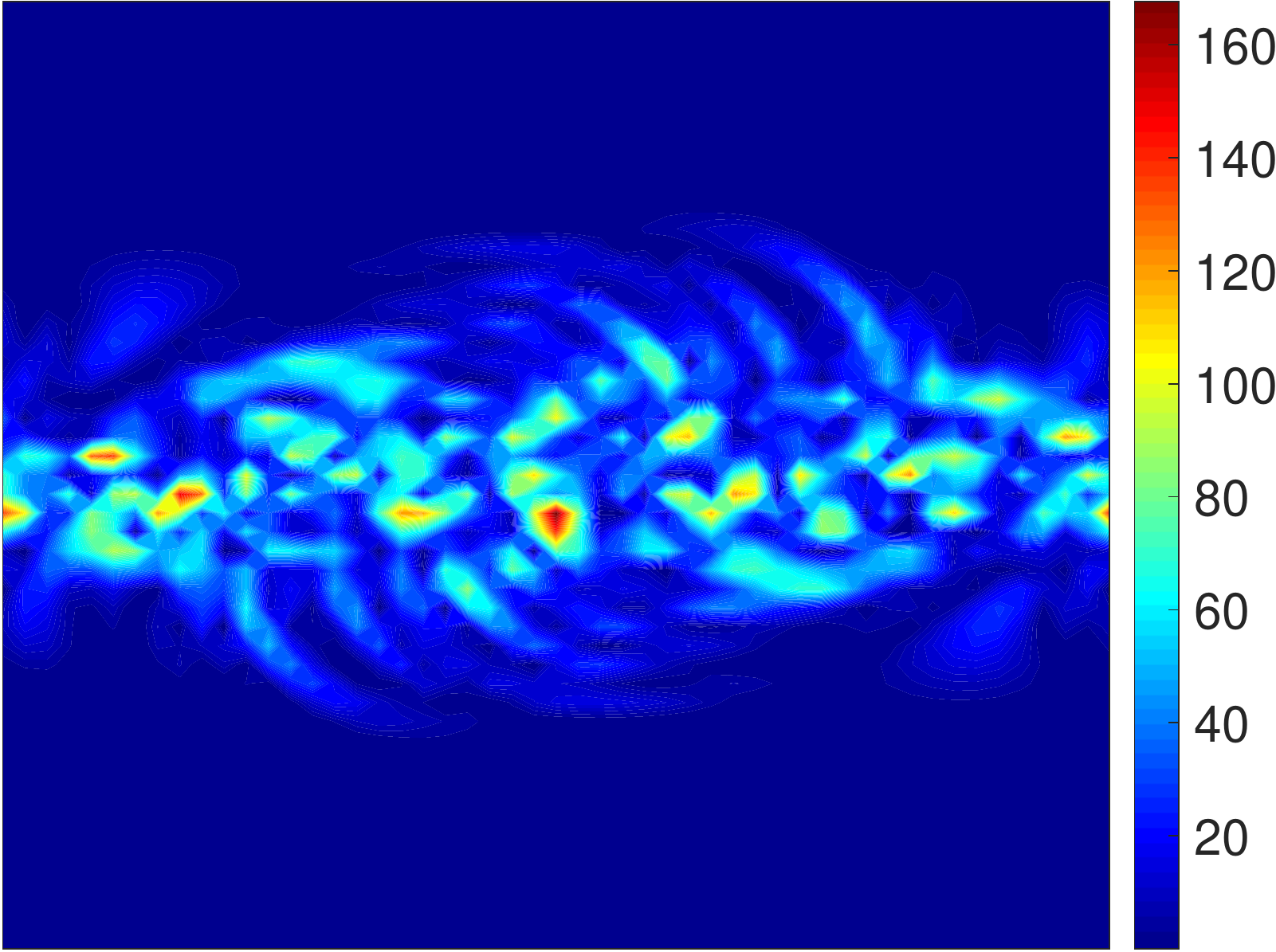}
\includegraphics[width=.24\textwidth, height=.12\textwidth,viewport=0 0 520 400, clip]{Emacsnapt5.pdf}
\end{center}
\caption{\label{khskew}
Shown above are results (shown as vorticity contours) of ROM simulations that used SKEW-FEM as the FOM,  at t=1, 2, 3, 4, 5 (from top to bottom) for the $Re=100$ KH tests using 50 modes.}
\end{figure}

We next repeat this same test, but now using SKEW-FEM as the FOM.  Results are shown in Figures \ref{khskew} and \ref{skewEE}.  Results are analogous to those 
in section~\ref{sec:kh-emac-fem}, where EMAC-FEM is used as the FOM.  The consistent ROM in this case (SKEW-ROM) dramatically outperforms EMAC-ROM and CONV-ROM.  While the inconsistent ROMs become unstable, the SKEW-ROM remains stable, matches the FOM energy and enstrophy, and matches the vorticity contours of the FOM although its final eddy location is shifted.
Hence these tests once again illustrate that best results come from FOM-ROM consistency.

These results also show that FOM-ROM consistency is not related to the particular formulation used at the FOM level.
Indeed, for both the EMAC-FEM used in section~\ref{sec:kh-emac-fem} and the SKEW-FEM used in this section, the best ROM results were obtained with the corresponding consistent formulation, i.e., EMAC-EMAC in section~\ref{sec:kh-emac-fem} and SKEW-SKEW in this section.

\subsubsection{Convergence of ROM to FOM depends on FOM/ROM consistency}
    \label{sec:kh-consistency}

\begin{figure}[!h]
\begin{center}
\includegraphics[width=.45\textwidth, height=.32\textwidth,viewport=0 0 520 400, clip]{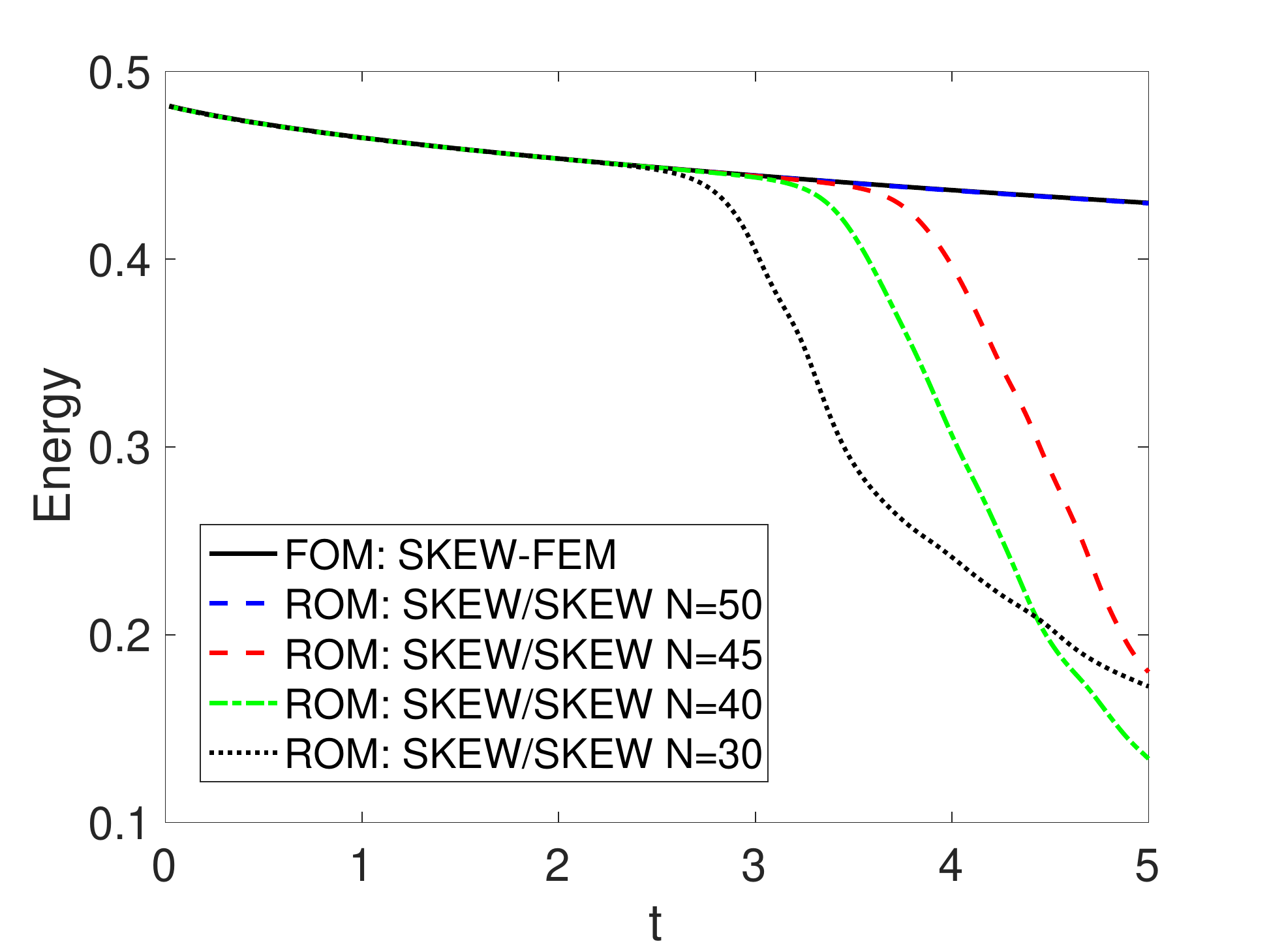}
\includegraphics[width=.45\textwidth, height=.32\textwidth,viewport=0 0 520 400, clip]{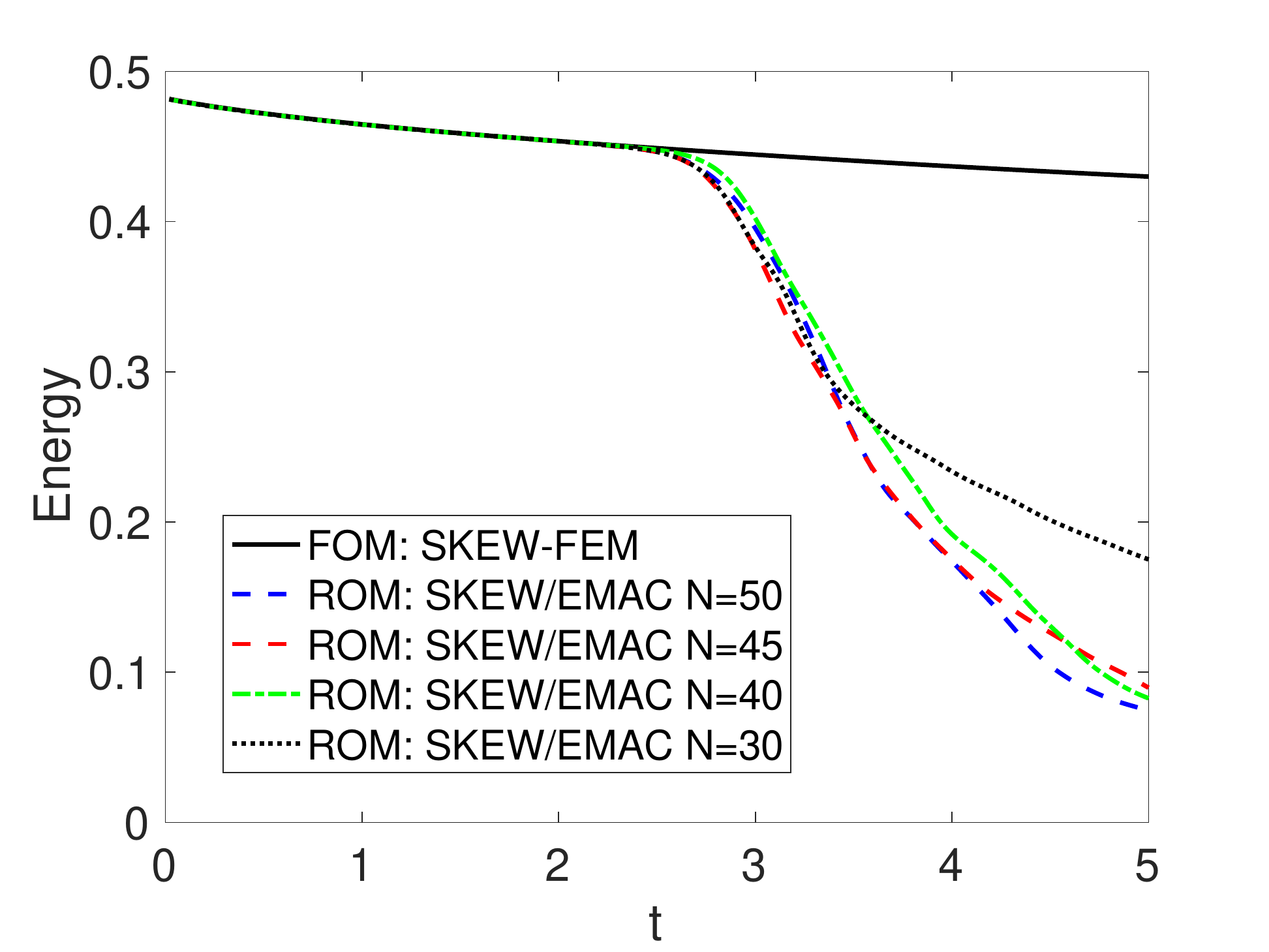} 
\end{center}
\caption{\label{skewconv}
Shown above is the energy plots for SKEW-FEM FOM together with SKEW-ROM with varying $N$ (left) and EMAC-ROM with varying $N$ (right).
}
\end{figure}

Finally, we illustrate \LR{ that the inconsistency of formulations causes an error that increasing the number of modes $N$ cannot fix (i.e., locking).  First, we compute the SKEW-FOM and SKEW-ROM (the consistent case) with N=30, 40, 45, and 50,} and plot their energy prediction along with the SKEW-FEM energy prediction in Figure \ref{skewconv} on the left (recall from Figure \ref{skewEE} on the left that SKEW-FEM energy is a relatively close prediction of the DNS energy prediction).  We observe that as $N$ increases, the SKEW-ROM energy converges to the SKEW-FEM energy, with the $N=50$ energy curve 
lying on top of the SKEW-FEM energy curve.  

\LR{ 
Figure \ref{skewconv} on the right shows the inconsistent case.  Here, SKEW-FOM is again used but now EMAC-ROM (the inconsistent case) with N=30, 40, 45, and 50 is computed.  The plot shows energy vs. time for each of these tests, and we observe that the EMAC-ROM does not converge 
to the FOM as $N$ increases; in fact, EMAC-ROM is very inaccurate after $t=2.5$ for all $N$ values tested.}

\section{Conclusions}
    \label{sec:conclusions}

In this paper, we \ti{\textit{analyzed} and} compared 
numerically ROMs that are FOM consistent (i.e., use the same nonlinearity discretization as the FOM) with ROMs that are FOM inconsistent (i.e., use a  nonlinearity discretization that is different from that used in the FOM).
Our first significant contribution is the rigorous proof that the ROMs that are not consistent with the FOM yield error bounds with additional terms 
dependent on the divergence error of the FOM.  Thus, the ROM solution can only recover the FOM solution to this level as the number of modes increases, i.e., can lock.  
Our second significant contribution is a  numerical investigation that shows the beneficial role of the FOM-ROM consistency in practical simulations.  Specifically, 
we considered channel flow around a cylinder and Kelvin-Helmholtz instability problems, and showed that the ROMs that are FOM consistent yield significantly more accurate results than the ROMs that are FOM inconsistent, and also illustrated the locking phenomena arising from FOM/ROM inconsistency. 
\ti{To our knowledge, this represents the first theoretical and numerical investigation of the FOM-ROM inconsistency with respect to the numerical discretization.}
\ti{This study brings new insight to the FOM-ROM consistency issue, \textit{proving}} 
that FOM-ROM consistency is significant not only when the computational model is modified (e.g., when different closures and stabilizations are used)~\ti{~\cite{ali2020stabilized,giere2015supg,girfoglio2021pod,pacciarini2014stabilized,rebollo2017certified,stabile2019reduced,strazzullo2021consistency}, } 
but also when fundamental discretization choices (e.g., the nonlinearity discretization) are made.

\ti{This first investigation of the FOM-ROM consistency with respect to the numerical discretization can be extended in several directions.}
An important research direction is the effect of nonlinear discretization and corresponding FOM-ROM consistency on conservation of physical properties at a ROM level.
At a FOM (FEM) level, the effect of nonlinear discretization on conservation of physical properties has been investigated, e.g., in~\cite{CHOR17,OR20}.
At a ROM level, the conservation of physical properties has been investigated, e.g., in~\cite{barone2009stable,chan2020entropy,chen2021bamcafe,karasozen2021structure,loiseau2018constrained,majda2012physics,mohebujjaman2017energy,mohebujjaman2019physically,sanderse2020non,yano2019discontinuous}.
To our knowledge, the effect of nonlinear discretization and corresponding FOM-ROM consistency on conservation of physical properties at a ROM level is still an open problem.

\ti{Probably the most important research direction is to \textit{prove} and illustrate numerically that FOM-ROM consistent approaches are more accurate than their inconsistent counterparts with respect to both computational modeling and numerical discretization for other settings.}
\ti{For example, we believe that the theoretical framework developed in this paper can be extended to other nonlinear problems, such as the Boussinesq equations and magnetohydrodynamics (MHD), \LR{ and with more effort to other PDEs that 
use formulations 
that are consistent at the continuous level but different at the discrete level.}
Another interesting research direction is the investigation of FOM-ROM consistency for 
}
different time stepping schemes, particularly ones that alter solution spaces such as projection and penalty methods. 
\ti{
An important research avenue is the theoretical and computational investigation of the FOM-ROM consistency when different spatial discretizations are used at the FOM and ROM levels.
For example, one could use a high order finite volume discretization at the FOM level and a low order finite element discretization at the ROM level.
Of course, the theoretical investigation of this important practical setting would pose additional challenges, e.g., the potential lack of Galerkin orthogonality, which was an essential tool used in the proof of Theorem~\ref{FOMROMerr}.  
}



\section*{Acknowledgements}
\ti{We thank the reviewers for the insightful comments and suggestions, which have significantly improved the paper.}


\end{document}